\documentclass[a4paper,11pt,oneside,reqno]{amsart}
\usepackage[all,cmtip]{xy}
\usepackage{amsmath}
\usepackage{amsthm}
\usepackage{subcaption}
\usepackage{amssymb}
\usepackage{graphicx}
\usepackage{mathrsfs}
\usepackage{tikz}
\usetikzlibrary{graphs}
\usepackage{enumerate}
\theoremstyle{plain}
\newtheorem{teore}{Theorem}[section]
\newtheorem{defin}[teore]{Definition}
\newtheorem{lem}[teore]{Lemma}
\newtheorem{coro}[teore]{Corollary}
\newtheorem{propo}[teore]{Proposition}

\newtheorem{pro}{Problem}
\newtheorem{claim}{Claim}[teore]
\newtheorem*{claim*}{Claim}
\newtheorem*{thm*}{Theorem}
\newtheorem*{defi*}{Definition}
\theoremstyle{remark}
\newtheorem{ejemplo}[teore]{{\sc Example}}
\newtheorem{ejemplos}[teore]{{\sc Examples}}
\newtheorem{notas}[teore]{{\sc Remark}}
\newcommand{\nrm}[1]{\|#1\|}

\newcommand{\tnrm}[1]{\|\hspace{-0.4mm}|#1\|\hspace{-0.4mm}|}

\newcommand{\prop}{\begin{propo}}
\newcommand{\fprop}{\end{propo}}
\newcommand{\cor}{\begin{coro}}
\newcommand{\fcor}{\end{coro}}

\newcommand{\IS}[2]{\mr{Is}_{#1}(#2)}

\newcommand{\defi}{\begin{defin}}
\newcommand{\fdefi}{\end{defin}}
\newcommand{\eje}{\begin{ejemplo}}
\newcommand{\feje}{\end{ejemplo}}
\newcommand{\ejes}{\begin{ejemplos}}
\newcommand{\fejes}{\end{ejemplos}}
\newcommand{\lema}{\begin{lem}}
\newcommand{\flema}{\end{lem}}
\newcommand{\teor}{\begin{teore}}
\newcommand{\fteor}{\end{teore}}
\newcommand{\nota}{\begin{notas}}
\newcommand{\fnota}{ \end{notas}}
\newcommand{\clam}{\begin{claim}}
\newcommand{\fclam}{\end{claim}}
\newcommand{\clams}{\begin{claim*}}
\newcommand{\fclams}{\end{claim*}}

\newcommand{\lclam}{\begin{lclaim}}
\newcommand{\flclam}{\end{lclaim}}
\newcommand{\prucl}{\prue[Proof of Claim:]}
\newcommand{\fprucl}{\fprue}
\newcommand{\ben}{\begin{enumerate}}
\newcommand{\een}{\end{enumerate}}
\newcommand{\bit}{\begin{itemize}}

\newcommand{\eit}{\end{itemize}}

\newcommand{\mc}[1]{\mathcal{#1}}
\newcommand{\mbb}[1]{\mathcal{#1}}
\newcommand{\mr}[1]{\mathrm{#1}}

\newcommand{\mk}[1]{\mathfrak{#1}}

\newcommand{\we}{\wedge}

\renewcommand{\bar}[1]{\overline{#1}}

\newcommand{\casos}{\begin{itemize}}
\newcommand{\fcasos}{\end{itemize}\setcounter{cs}{1}}

\newcommand{\fin}{\textsc{FIN}}

\newcommand{\srk}[2]{\mr{srk}_{#2}(#1)}
\newcommand{\rk}[1]{\mr{rk}(#1)}

\newcommand{\ro}{\varrho}

\newcommand{\conj}[2]{ \{ {#1}\,:\,{#2} \} }

\newcommand{\ou}{\omega_{1}}

\newcommand{\om}{\omega}

\newcommand{\ip}{\sqsubseteq}

\newcommand{\ka}{\kappa}

\newcommand{\buit}{\emptyset}

\newcommand{\ga}{\gamma}

\DeclareMathOperator{\conc}{^\smallfrown}

\newcommand{\al}{\alpha}
\newcommand{\be}{\beta}
\newcommand{\de}{\delta}
\newcommand{\De}{\Delta}
\newcommand{\la}{\lambda}

\newcommand{\sig}{\sigma}

\newcommand{\vphi}{\varphi}

\newcommand{\vep}{\varepsilon}

\newcommand{\eqs}{{\mathfrak X}}

\newcommand{\N}{{\mathbb N}}

\newcommand{\caso}[1]{\noindent {\sc Case  {#1}.  } }

\renewcommand{\ker}{\mathrm{Ker}}
\newcommand{\rest}{\upharpoonright}

\newcommand{\supp}{\mathrm{supp\, }}

\newcommand{\ran}{\mathrm{ran\, }}

\newcommand{\con}{\subseteq}

\newcommand{\cones}{\varsubsetneq }

\newcommand{\prue}{\begin{proof}}
\newcommand{\fprue}{\end{proof}}
\makeindex
\title[Homogeneous families on trees and subsymmetric basic sequences]{Homogeneous families on trees\\and subsymmetric basic sequences}
\author{C. Brech}
\address{Departamento de Matem\'atica, Instituto de Matem\'atica e Estat\'\i stica, Universidade de S\~ao Paulo,
Caixa Postal 66281, 05314-970, S\~ao Paulo, Brazil}
\email{brech@ime.usp.br}
\author{J. Lopez-Abad}
\address{Instituto de Ciencias Matem\'{a}ticas (ICMAT), CSIC-UAM-UC3M-UCM, C/ Nicol\'as Cabrera 13-15,
Campus Cantoblanco, UAM, 28049 Madrid, Spain}
\email{abad@icmat.es}
\author{S. Todorcevic}
\address{Institut de Math\'{e}matiques de Jussieu, UMR 7586, 2 place Jussieu - Case 247, 75222 Paris Cedex 05, France and Department of Mathematics, University of Toronto, Toronto, Canada, M5S 2E4}
\email{stevo@math.toronto.edu}
\thanks{The first author was partially supported by FAPESP grants (2012/24463-7 and 2015/26654-2) and by
  CNPq grants (307942/2012-0 and 454112/2015-7). The second author was partially supported by FAPESP grant (2013/24827-1),
    by the Ministerio de Econom\'\i a y Competitividad under grant MTM2012-31286 (Spain). 
    The first and the third author were partially supported by USP-COFECUB grant (31466UC). The third author was also supported by grants from
CNRS and NSERC (455916)}
    \thanks{This research was partially done whilst the authors were   visiting fellows at the Isaac Newton
Institute for Mathematical Sciences in the programme ``Mathematical, Foundational
and Computational Aspects of the Higher Infinite'' (HIF)}

\subjclass[2010]{03E05,\,05D10,\,46B15,\,46B06,\,46B26}

\keywords{Families of finite sets, Spreading models, 
Nonseparable Banach spaces}
\begin{document}

\begin{abstract}
We study  density requirements  on a given Banach space  that  guarantee the existence of subsymmetric
 basic sequences by extending Tsirelson's well-known space to larger index sets.  We
prove that for every cardinal $\kappa$ smaller than the first Mahlo cardinal there is a reflexive Banach space of
density $\kappa$ without subsymmetric basic sequences.   As for Tsirelson's space,  our construction is
based on the existence of a  rich collection of homogeneous families on large index sets for which one can estimate  the complexity on any given infinite set. This is used to describe detailedly the  asymptotic structure of the spaces. The collections of families  are of independent
interest and their existence is proved inductively. The fundamental stepping up argument is the analysis of
such collections of families on trees.

\end{abstract}
\maketitle

\section{Introduction}
Recall that a set of \emph{indiscernibles} for a given structure $\mathcal{M}$ is a subset $X$ with a total
ordering $<$ such that for every positive integer $n$ every two increasing $n$-tuples
$x_1<x_2<\cdot\cdot\cdot <x_n$ and $y_1<y_2<\cdot\cdot\cdot <y_n$ of elements of $X$ have the same properties
in $\mathcal{M}.$ A simple way of finding an extended structure on $\ka$ without an infinite set of
indiscernibles is as follows. Suppose that $\mc F$ is a family of finite subsets of $\ka$ that is compact, as
a natural subset of the product space $2^\ka$, and \emph{large}, that is,  every infinite subset of $\ka$ has
arbitrarily large subsets in $\mc F$. Let $\mc M_\mc F$ be the structure  $(\ka, (\mc F\cap [\ka]^n)_n)$ that
has $\ka$ as universe and that has infinitely many $n$-ary relations $\mc F\cap [\ka]^n\con [\ka]^n$. It is
easily seen that $\mc M_\mc F$ does not have infinite indiscernible sets.

While  in set theory and model theory indiscernibility is a well-studied and unambiguous notion, in the
context of the Banach space theory it has several versions, the most natural one being the notion of a
\emph{subsymmetric sequence} or a \emph{subsymmetric set}.    In a normed space $(X,\nrm{\cdot})$, a sequence
$(x_\al)_{\al\in I}$ indexed in an ordered set $(I,<)$ is called $C$-\emph{subsymmetric} when
\begin{equation*}
\nrm{\sum_{j=1}^n a_j x_{\al_j}}\le C\nrm{\sum_{j=1}^n a_j x_{\be_j}}
\end{equation*}
for every sequence of scalars $(a_j)_{j=1}^n$ and every $\al_1<\cdots < \al_n$ and $\be_1<\cdots <\be_n$ in
$I$.   When $C=1$ this corresponds exactly to the notion of an indiscernible set and it is easily seen that
this can always be assumed by renorming $X$ with an appropriate equivalent norm.   Another closely related
notion is \emph{unconditionality}. Recall that a sequence $(x_i)_{i\in I}$ in some Banach
space is \emph{$C$-unconditional} whenever
\begin{equation*}
\nrm{\sum_{i\in I} \theta_i a_i x_{i}}\le C\nrm{\sum_{i\in I} a_i x_{i}}
\end{equation*}
for every sequence of scalars $(a_i)_{i\in I}$ and every sequence $(\theta_i)_{i\in I}$ of signs. Of
particular interest are the indiscernible coordinate systems, such as the Schauder basic sequences. The unit
bases of the classical sequence spaces $\ell_p$, $p\ge 1$ or $c_0$ (in any density) are subsymmetric and
unconditional (in fact, symmetric, i.e. indiscernible by permutations) bases. Moreover, every
basic sequence in one of these spaces has a symmetric subsequence.  But this is not true in general:
there are basic sequences without unconditional
subsequences, the simplest example being the summing basis of $c_0$.  However, it is more difficult
  to find a weakly-null basis without unconditional subsequences (B. Maurey
and H. P. Rosenthal \cite{MaRo}). Now we know  that there are Banach spaces without unconditional basic
sequences. The first such example was given by  W. T. Gowers and B. Maurey \cite{GoMa}, a space which was
moreover reflexive. Concerning subsymmetric sequences, we mention that  the unit basis of the Schreier space
\cite{Schr} does not have subsymmetric subsequences, and the Tsirelson space \cite{Tsi} is the first example of
a reflexive space without subsymmetric basic sequences.

All these are separable spaces so it is natural to ask if large spaces must contain infinite   unconditional
or subsymmetric sequences, since from the theory of large cardinals we know that infinite indiscernible sets
exist in large structures. In general, this is a consequence of certain Ramsey principles (i.e.
higher-dimensional versions of the pigeonhole principles). Indeed, it was proved by Ketonen \cite{Ke} that
Banach spaces with density bigger than  the first $\om$-Erd\H{o}s cardinal have subsymmetric sequences.  Recall
that a cardinal number $\ka$ is called $\om$-Erd\H{o}s  when every  countable coloring of the collection of
finite subsets of $\ka$ has an infinite subset $A$ of $\ka$ where  the color of  a given finite subset $F$ of
$A$ depends only on the cardinality of $F$. Such cardinals are large cardinals, and their existence cannot be
proved on the basis of the standard axioms of set theory. It is therefore natural to ask what is the minimal
cardinal number $\mk{nc}$ ($\mk{ns}$) such that every Banach space of density at least $\mk{nc}$ (resp.
$\mk{ns}$) has an unconditional (respectively,  subsymmetric) basic sequence. It is natural to consider also the
relative versions of these cardinals restricted to various classes of Banach spaces  like, for example, the
class of reflexive spaces where we use the notations $\mk{nc}_\mr{refl}$ and $\mk{ns}_\mr{refl}$,
respectively. Note that it follows from  E. Odell's partial unconditionality result \cite{Od1}  that every
weakly-null subsymmetric basic sequence  is unconditional, hence $\mk{nc}_{\mr{refl}}\le \mk{ns}_{\mr{refl}}$.
Moreover, an easy application of Odell's result and Rosenthal's $\ell_1$-dichotomy gives that $\mk{nc}\le
\mk{ns}$.

Concerning lower bounds for these cardinal numbers, it was proved by S. A. Argyros and A. Tolias \cite{ArTol}
that $\mk{nc}>2^{\aleph_0}$, and by E. Odell  \cite{Od}  that  $\mk{ns}>2^{\aleph_0}$.  For the reflexive
case, we know that $\mk{nc}_\mr{refl}>\aleph_1$ (\cite{ArLoTo}), and in the recent paper \cite{ArMo},
S. A. Argyros and P. Motakis proved that
$\mk{ns}_\mr{refl}>2^{\aleph_0}$.  Finally we mention that in
\cite{LoTo} the Erd\H{o}s cardinals are characterized in terms of the existence of compact and large families and the
sequential version of $\mk{ns}_\mr{refl}$. More precisely it is proved that the first Erd\H{o}s cardinal
$\mk{e}_\om$ is the minimal cardinal $\ka$ such that every long weakly-null basis of length $\ka$  has a
subsymmetric basic sequence, or equivalently the minimal cardinal $\ka$ such that there is no compact and large
family of finite subsets of $\ka$.

The study of upper bounds is of different nature and seems to involve  more advanced set-theoretic
considerations connected to large cardinal principles.   This can be seen, for example,  from the
aforementioned result of Ketonen or from results of P. Dodos, J. Lopez-Abad and S.
Todorcevic  who proved in \cite{DoLoTo} that $\mk{nc}\le \aleph_\om$ holds consistently relative to the
existence of certain large cardinals and who proved  in \cite{LoTo}  that Banach's  Lebesgue measure
extension axiom implies that $\mk{nc}_\mr{refl}\le 2^{\aleph_0}$.

In this paper we continue the research on the existence of subsymmetric sequences in a normed space of large
density,    and we prove that  $\mk{ns}_\mr{refl}$ is rather large, distinguishing thus the cardinals
$\mk{ns}_\mr{refl}$  and $\mk{nc}_\mr{refl}$. In contrast to the   sequential version of $\mk{ns}_\mr{refl}$,
that is closely linked to indiscernibles of \emph{relational} structures $\mc M_\mc F$ for compact and large
families $\mc F$, the full version of $\mk{ns}_\mr{refl}$ is more related to the existence of indiscernibles
in structures that are not just relational but also have operations, suggesting that not only we need to
understand families on finite sets but also ``operations" with them. In the separable context, this is well-known
and can be observed in the construction of the Tsirelson space, where finite products of the \emph{Schreier
family} are used in a crucial way. The natural approach in the non-separable setting would be to generalize
Tsirelson's construction
 using analogues of the Schreier family, certain large compact families, on larger index
sets. However, in the uncountable level these families cannot be spreading and therefore, if one just copies Tsirelson's
construction on the basis of them, the corresponding non-separable Tsirelson-like spaces will always
contain almost isometric copies of $\ell_1$ (\cite[Theorem 8.2]{LoTo}).  This lead us to change our
perspective and use the well-known interpolation technique
 \cite[Example 3.b.10]{LiTz}, an approach that appeared  recently in the work of Argyros and Motakis mentioned above.
In this perspective, a key tool is a suitable operation $\times$, that we call \emph{multiplication},  of
compact families of finite sets. In fact, the multiplication is an operation which associates
to a family $\mc F$ on the fixed index set $I$ and family $\mc H$ on $\om$, a family $\mc F \times \mc H$ on $I$
which has, in a precise sense, many elements of the form $\bigcup_{n\in x} s_n$, where $x\in \mc H$  and $(s_n)_{n<\om}$ is an
arbitrary sequence of elements of $\mc F$.   It is well known that such multiplication exists in $\om$
and it models in some way the ordinal multiplication on uniform families. It is also the main tool
to define the generalized Schreier families on $\om$, vastly used in modern Banach space theory to study
ranks of compact notions (e.g. summability methods), or of asymptotic notions (e.g. spreading models). These
are \emph{uniform} families, so that any restriction of them looks like the entire family. We generalize
this property to the uncountable level by defining \emph{homogeneous} families, that
despite being uncountable families on large index sets, have countable Cantor-Bendixson rank which moreover does
not change substantially
when passing to restrictions. In particular, if $\mc F$ is homogeneous, then the
structure $\mc M_\mc F$ does not have infinite sets of indiscernibles, but we also get lower and upper bounds for the
rank of the collection of their (finite) sets of indiscernibles.

We then introduce the notion of a \emph{basis} of families, which is a rich collection of homogeneous families
admitting a multiplication, and we prove that they exist on quite large cardinal numbers. The existence of
such bases is proved inductively. For example, we prove that if $\ka$ has  a basis then $2^\ka$ has also a
basis.  This is done by representing $2^\ka$ as the complete binary tree $2^{\le \ka}$, and observing that we can use the
\emph{height} function $\mr{ht}:2^{\le \ka}\to \ka+1$ to pull back a basis on $\ka$ to a restricted version of
basis  on $2^{\le k}$,   consisting of homogeneous families of finite chains of $2^{\le \ka}$. Actually, we
prove the following more general equivalence (Theorem \ref{4jtierjteoirtertrtr344}).

\begin{thm*}
For an infinite  rooted tree $T$  the following are equivalent.
\begin{enumerate}[(a)]
\item There is a basis of families  on $T$.
\item  There is a basis of families  consisting of chains of $T$ and there is a basis consisting of  antichains of $T$.
\end{enumerate}
\end{thm*}
In particular, we obtain a basis on $2^\om$ that can be used to build a reflexive space of density
$2^{\om}$ without subsymmetric basic sequences, giving another proof of the result in \cite{ArMo}. Also,
one proves inductively that for every cardinal number $\kappa$ smaller than the first inaccessible cardinal,
there is a basis on $\kappa$ and a corresponding Banach space of density $\kappa$ with similar properties.
We then use Todorcevic's method of walks on ordinals \cite{To} to build trees on cardinals
up to the first Mahlo cardinal number and find examples of reflexive Banach spaces of large densities
without subsymmetric basic subsequences. Moreover,
as observed above for the structure $\mc M_\mc F$, we can bound the complexity of the (finite) subsymmetric
basic sequences and obtain the following.

\begin{thm*}
Every  cardinal  $\ka$ below the first Mahlo cardinal has a basis. Consequently, for every such cardinal
$\ka$ and every $\al<\ou$, there is a reflexive Banach space $\mk X$ of density $\ka$ with a long
unconditional basis and such that every bounded sequence in   $\mk X$  has an $\ell_1^\al$-spreading model
subsequence but the space $\mk X$ does not have $\ell_1^{\be}$-spreading model subsequences for $\be$ large
enough,  depending only on $\al$. In particular, $\mk X$ contains no infinite subsymmetric basic sequence.
\end{thm*}

The paper is organized as follows. In Section 2 we introduce some basic topological, combinatorial and
algebraic facts on families of finite chains of a  given partial ordering.  We then define homogeneous
families and bases of them. We present some upper bounds for the topological rank of a family
that uses the well-known Ramsey property of barriers on $\om$ and that will be used several times along the paper. In Subsection \ref{transfer_bases} we  give useful methods to transfer bases between partial orderings.    Section \ref{bases_on_trees} is the main part of this paper. The
main object we study is, given a tree, the collection $\mc A\odot_T \mc C$ of all finite subtrees of $T$
whose chains are in a fixed family $\mc C$ and such that the family of immediate successors of a given node
is in another fixed family $\mc A$. We study $\mc A\odot_T \mc C$ both combinatorially and topologically. The
combinatorial part  is based on the canonical form of a sequence of finite subtrees, and allow us to define a
natural multiplication. The topological one consists in finding upper bounds of the rank of the family $\mc
A\odot_T \mc C$ in terms of the corresponding ranks of the families $\mc A$ and $\mc C$, much in the spirit
of how one easily bounds the size of a finite tree from its height and splitting number. This operation
allows to lift bases on chains and of immediate successors to bases on the whole tree, our main result of
this work done in Theorem \ref{4jtierjteoirtertrtr344}. We apply this in Section 4 to prove that cardinal
numbers smaller than the first Mahlo cardinal have a basis. To do this, we represent such cardinals as nodes
of a tree having bases on chains and on immediate successors. We achieve this last part by proving several
principles of transference of basis. Finally, we use bases to build reflexive Banach spaces without subsymmetric
basic sequences.

\section{Basic definitions}

Let $I$ be a set. A set $\mc F$ is called a \emph{family on $I$} when the
elements of $\mc F$ are finite subsets of $I$.  Let $\mc P=(P,<)$ be  a partial ordering. A {\em family on
chains of $\mc P$} is a family on $P$  consisting of chains of $\mc P$, ie. totally ordered subsets of $P$. Let $\mr{Ch}_<$ be the collection of
all \emph{finite} chains of $\mc P$.   Given $k\le \om$, let
\begin{align*}
[I]^k:=& \conj{s\con I}{\#s=k}, & [P]_{<}^k:=&[P]^k \cap \mr{Ch}_<, \\
[I]^{\le k}:=& \conj{s\con I}{\#s\le k}, & [P]_<^{\le k}:=& [P]^{\le k} \cap \mr{Ch}_<.
\end{align*}
For a family $\mc F$ on     $I$ and $A\con I$, let $\mc F\rest A:=\mc F\cap \mc P(A)$. Recall that a family
$\mc F$ on $I$ is \emph{hereditary} when it is closed under subsets and it is \emph{compact} when it is a
closed subset of $2^I:=\{0,1\}^I$, after identifying each set of $\mc F$ with its characteristic function. In
this case, $\mc F$ is a scattered compact space.  Since each element of $\mc F$ is finite, it is not
difficult to see that $\mc F$ is compact if and only if every sequence $(s_n)_{n\in\om}$ in $\mc F$ has a
subsequence  $(t_n)_{n\in \om}$ forming a $ \De$-system with root in $\mc F$, that is, such that
$$t_{k_0}\cap t_{k_1}=t_{l_0}\cap t_{l_1} \in \mc F \text{ for every $k_0\neq k_1$ and $l_0\neq l_1$.}$$
The intersection $t_k\cap t_l$, $k\neq l$ is called the \emph{root} of $(t_n)_n$. By weakening the notion of  compactness, we say that $\mc F$ is \emph{pre-compact} if every sequence in $\mc F$ has a $\De$-subsequence (with root not  necessarily in $\mc F$).  It is easy to see that $\mc F$ is pre-compact if and only if its $\con$-closure $\conj{s\con I}{ s\con t  \text{ for some $t\in \mc F$}}$ is compact.

Recall the \emph{Cantor-Bendixson derivatives} of a topological space $X$:
$$   X^{(0)}:=X, \quad X^{(\al)}=\bigcap_{\be<\al} (X^{(\be)})'$$
 where $Y'$ denotes the collection of accumulation  points of $Y$, that is, those points $p\in X$ such that
each of its open neighborhoods has infinitely many points in $Y$.  The minimal ordinal $\al$ such that
$X^{(\al+1)}=X^{(\al)}$ is called the \emph{Cantor-Bendixson rank} $\mr{rk}_\mr{CB}(X)$ of $X$.  In the case
of a  compact family $\mc F$ on an index set $I$,  being scattered, its Cantor-Bendixson index is the first
$\al$ such that $\mc F^{(\al)}=\buit$, and therefore it must be a successor ordinal.

\defi[Rank, small rank and homogeneous families]
Given a compact family $\mc F$  on some index set $I$, let
$$\mr{rk}(\mc F):=\mr{rk}_\mr{CB}(\mc F)^{-}$$
where   $(\al+1)^-=\al$. We say that a compact family $\mc F$ is countably ranked when $\mr{rk}(\mc F)$ is
countable.  Let $\mc P$  be a partial ordering. For a given family $\mc F$ on chains of $\mc P$, the \emph{small rank
relative to $\mc P$} of $\mc F$ is
$$\mr{srk}_\mc P(\mc F):= \inf\conj{\mr{rk}(\mc F\rest C)}{C \text{ is an infinite chain of $\mc P$}}.$$
A compact and hereditary family $\mc F$ on chains of  $\mc P$ is called \emph{$(\al,\mc P)$-homogeneous} if $\mc F=\{\buit\}$ if $\al=0$, and if  $1\le \al<\ou$, then  
$[ P]^{\le 1}\con \mc F$ and
$$  \al= \mr{srk}_\mc P(\mc F) \le \mr{rk}(\mc F) <\iota(\al), $$
where $\iota(\al)$ is defined below in Definition  \ref{kyuyu39766}.  $\mc F$ is \emph{$\mc P$-homogeneous}
if it is $(\al,\mc P)$-homogeneous for some $\al<\ou$. \fdefi

Observe that these notions coincide for any two total orders on the fixed index set, as being a chain with respect to a total order does not depend on the total order itself. Hence, when $\mc F$ is a family on some index set $I$, we use $\mr{srk}$, $\al$-homogeneous and homogeneous for the corresponding notions with respect to any total order on $I$. If $\mc F$ is a compact family on a countable index set $I$, then $\mr{rk}(\mc F)$ is countable and, therefore, the small rank of $\mc F$ with respect to any partial order on $I$ is countable. In general, $\#\mr{rk}(\mc F)\le \# I$ and the extreme case can be  achieved.

\defi \label{kjo347755565}
The \emph{normal Cantor form} of an ordinal $\al$ is the unique expression $\al=\om^{\al[0]}\cdot n_0[\al]
+\cdots + \om^{\al[k]}\cdot n_k[\al]$ where $\al\ge\al[0]>\al[1]>\cdots >\al[k]\ge 0$ and $n_i[\al]<\om$ for
every $i\le k$.

Suppose that $*$  is an operation on countable ordinals and suppose that $\al>0$ is a countable ordinal.  We
say that $\al$  is  \emph{$*$-indecomposable} when $\be * \ga <\al$ for every $\be,\ga<\al$.
\fdefi
\nota
It it well-known that
\begin{enumerate}[(i)]
\item $\al$ is sum-indecomposable if and only $\al= \om^\be$.
\item $\al>1$ is product-indecomposable if and only if $\al= \om^\be$ for some sum-indecomposable  $\be$.
\item  For $\al>\om$, $\al$ is exponential-indecomposable if and only if $\al= \om^\al$.
\item   product-indecomposability imply sum-indecomposability, and exponential-indecomposability imply
product and sum indecomposability.
\end{enumerate}
\fnota
So, $1,\om,\om^2$ and $1,\om,\om^\om$  are the first 3 sum-indecomposable, and product-indecomposable
ordinals, respectively. If we define, given $\al<\ou$, $\bar{\al}_0:=\al$,
$\bar{\al}_{n+1}:=\omega^{\bar{\al}_n}$ and $\bar{\al}_\om:=\sup_n \bar{\al}_n$, then $\om, \bar{\om}_\om,
\bar{(\bar{\om}_\om)}_\om$ are the first 3 exponential-indecomposable ordinals.
 We will use \emph{exp-indecomposable} to refer to exponential-indecomposable
ordinals.
\defi\label{kyuyu39766}
Given a countable ordinal $\al$, let
$$\iota(\al)=\min\conj{\la> \al}{\la \text{ is exp-indecomposable}}.$$
Let $\mr{Fn}(\ou,\om)$ be the collection of all functions $f:\ou\to \om$ such that $\supp
f:=\conj{\ga<\ou}{f(\ga)\neq 0}$ is finite. Endowed with the pointwise sum $+$, $(\mr{Fn}(\ou,\om),+)$ is
an ordered commutative monoid.  Let $\nu: \ou\to \mr{Fn}(\ou,\om)$ be defined by $\nu(\al)(\ga)=n_i[\al]$ if
and only if $\ga=\al[i]$. Let $\sig: \mr{Fn}(\ou,\om)\to \ou$ be defined by $\sig(f)=\sum_{i\le k}\om^{\al_i}
\cdot f(\al_i)$, where $\{\al_0>\cdots >\al_n\ge 0\}=\supp f$. In other words, $\sig$ is the inverse of
$\nu$. Given $\al,\be<\ou$, the \emph{Hessenberg sum} (see e.g. \cite{Si}) is defined by
\begin{equation*}
\al \dot{+} \be := \sig(\nu(\al)+\nu(\be)).
\end{equation*}
 \fdefi
It is easy to see that if $\alpha$ is exp-indecomposable, then it is $\dot{+}$-indecomposable.

\defi
Let $\mc F$ and $\mc G$ be families on  chains of a partial ordering $\mc P$. Define
\begin{align*}
\mc F\cup \mc G:=& \conj{s\con P}{s\in \mc F \text{ or }s\in \mc G},\\
\mc F\sqcup_\mc P \mc G:=& \conj{s\cup t}{s\cup t \text{ is a chain and } s\in \mc F,\, t\in \mc G}, \\
\mc F\sqcup \mc G:=& \conj{s\cup t}{s\in \mc F,\, t\in \mc G}, \\
\mc F\boxtimes_\mc P (n+1) :=& (\mc F \boxtimes_\mc P n) \sqcup_\mc P \mc F; \quad \mc F \boxtimes_\mc P 1 := \mc F, \\
\mc F\boxtimes (n+1) :=& (\mc F \boxtimes n) \sqcup \mc F; \quad \mc F \boxtimes 1 := \mc F.
\end{align*}
\fdefi
Observe that when $\mc P$ is a total ordering  the operations $\sqcup_\mc P$  and $\sqcup$ are the same.

\prop \label{uygfuyhgu1} The operations $\cup$, $\sqcup_\mc P$ and $\sqcup$  preserve pre-compactness and hereditariness.  Moreover, if  $\mc F$ and $\mc G$ are countably ranked families on chains of $\mc P$, then
\begin{enumerate}[(i)]
\item  $\mr{rk}(\mc F\cup \mc G)=\max\{\mr{rk}(\mc F), \mr{rk}(\mc
G)\}$,
\item  $\mr{rk}(\mc
F\sqcup \mc G) = \mr{rk}(\mc F)\dot{+}\mr{rk}(\mc G)$,

\item  $ \mr{rk}(\mc F\sqcup_\mc P \mc G)  \le  \mr{rk}(\mc F)\dot{+}\mr{rk}(\mc G)$.
\end{enumerate}
Consequently, if $\mc F$  and $\mc G$   are $\mc P$-homogeneous, then $\mc F\cup \mc G$, $\mc F\sqcup_\mc P
\mc G$ and $\mc F\sqcup \mc G$ are $(\ga,\mc P)$-homogeneous with $\ga\ge \max\{\mr{srk}_\mc P(\mc
F),\mr{srk}_\mc P(\mc G)\}$ .

\fprop
\prue
It is easy to see that if $\mc F$ and $\mc G$ are pre-compact, hereditary, then $\mc F * \mc G$ is
pre-compact, hereditary, for $*\in \{\cup,\sqcup_\mc P,\sqcup\}$.    Let us see (i):   An easy inductive
argument shows that $(\mc F\cup \mc G)^{(\al)}=\mc F^{(\al)}\cup \mc G^{(\al)}$ for every countable $\al$. To prove 
(ii) and (iii), notice that by  a general fact, for every compact spaces $K$ and $L$ and every $\al$ one has that
\begin{equation}\label{iuuir78yug48394326yr122}
(K\times L )^{(\al)}= \bigcup_{
  \be \dot{+}\ga=\al
 } (K^{(\be)} \times L^{(\ga)}).
\end{equation}
When $\mr{rk}(K)$ and $\mr{rk}(L)$ are countable, we have that $\mr{rk}(K\times L)=\mr{rk}(K) \dot+ \mr{rk}(L)$.  The proof of
\eqref{iuuir78yug48394326yr122} is done by induction on $\al$ and by considering the case when $\al$ is
sum-indecomposable or not.  Now let $\mc F$ and $\mc G$ be countably ranked. Suppose that $\mc P$ is a
total ordering, and let $\mc F\times \mc G\to \mc F\sqcup \mc G$, $(s,t)\mapsto s\cup t$. This is clearly
continuous, onto and finite-to-one, so $\mr{rk}(\mc F\sqcup\mc G)= \mr{rk}(\mc F\times \mc G)=\mr{rk}(\mc F)
\dot+\mr{rk}(\mc G)$, which concludes the proof of (ii). If $\mc P$ is in general a partial ordering, then it follows from this that
$\mr{rk}(\mc F \sqcup_\mc P \mc G)\le \mr{rk}(\mc F) \dot+\mr{rk}(\mc G)$, proving (iii). Now suppose that
$\mc F$ and $\mc G$ are $\mc P$-homogeneous. We have clearly that
\begin{equation}\label{iojoi43oij43ioj43oij34}
\max\{\mr{srk}_\mc P(\mc F),\mr{srk}_\mc P(\mc G)\}\le \min\{ \mr{srk}_\mc P(\mc F\cup \mc G),\mr{srk}_\mc
P(\mc F\sqcup_\mc P \mc G), \mr{srk}_\mc P(\mc F\sqcup \mc G)\}.
\end{equation}
 On the other hand, $\max\{\mr{rk}(\mc F\cup \mc G),\mr{rk}(\mc F\sqcup_\mc P \mc G), \mr{rk}(\mc F\sqcup
\mc G) \}\le \mr{rk}(\mc F)\dot+\mr{rk}(\mc G)$. Since $\mr{rk}(\mc F),\mr{rk}(\mc G)<
\la:=\max\{\iota(\mr{srk}_\mc P(\mc F)), \iota(\mr{srk}_\mc P(\mc G))$, it follows by the indecomposability
of $\la$ that $\mr{rk}(\mc F)\dot+\mr{rk}(\mc G)<\la$. This, together with \eqref{iojoi43oij43ioj43oij34}
gives the desired result.
\fprue

\subsection{Bases of homogeneous families}
We recall  a  well-known generalization of Schreier families on $\om$, called uniform families. We  are going to use them
mainly as a tool to compute upper bounds of ranks  of operations of compact families. We use the following
standard notation: given $M,s,t\con \om$,  we write $s<t$ to denote that $\max s< \min t$ and let $M/s:=
\conj{m\in M}{s<m}$. Notice that a family $\mc F$ on $\om$ is pre-compact if and only if every sequence in $\mc F$ has a \emph{block} $\De$-subsequence $(s_n)_{n\in \om}$, that is, such that $s<s_{m}\setminus s< s_{n}\setminus s$ for every $m<n$, where $s$ is the root of $(s_n)_n$.
We write $s\ip t$ to denote that $s$ is an initial part of $t$, that is,   $s\con t$  and $t \cap (\max s +1)
=s$, and $s\sqsubset t$ to denote that $s\ip t$ and $s\neq t$.

 \defi
Given a family $\mc F$ on $\om$ and $n<\om$, let
$$\mc F_{\{n\}}:=\conj{s\con \om}{n<s\text{ and }\{n\}\cup
s\in \mc F}.$$ Let $\al$ be a countable ordinal number, and let $\mc F$ be a family on an infinite subset
$M\con \om$. The family $\mc F$ is called an \emph{$\al$-uniform family on $M$} when $\buit\in \mc F$ and
\begin{enumerate}[(a)]
\item  $\mc F=\{\buit\}$ if $\al=0$;
\item  $\mc F_{\{n\}}$ is $\be$-uniform on $M/ n$ for every $n\in M$, if $\al=\be+1$;
\item  $\mc F_{\{n\}}$ is $\al_n$-uniform on $M/ n$ for every $n\in M$ and $(\al_n)_{n\in M}$ is an  increasing sequence such that  $\sup_{n\in M} \al_n= \al$, if $\al$
is limit.
\end{enumerate}
\fdefi
It is important to remark that uniform families are not uniform fronts, which were introduced by
P. Pudlak and V. R\"{o}dl in \cite{PR} following previous works of C. Nash-Williams. Recall that a family $\mc B$ on $M$
is called an $\al$-uniform front on $M$ when $\mc B=\{\buit\}$ if $\al=0$, and if $\al>0$ then $\buit\notin \mc
B$, and $\mc B_{\{n\}}$ is a $\ga$-uniform front on $M/n$ for every $n\in M$, if $\be=\ga+1$, and    $\mc
B_{\{n\}}$ is a $\al_n$-uniform front on $M/n$ for every $n\in M$ and $(\al_n)_{n\in M}$ is increasing with
$\sup_{n\in M}\al_n$.  In fact, given a uniform family $\mc F$, the collection of its $\con$-\emph{maximal}
elements $\mc F^\mr{max}$ is a uniform front, and we can recover a uniform family out of a uniform front $\mc B$ by taking its closure under initial parts  $\mc B^{\ip}$, that is, the collection of initial parts of elements of $\mc B$:

\prop
\begin{enumerate}[(a)]
\item Every uniform family is compact.
\item The following are equivalent:
\begin{enumerate}[(b.1)]
\item  $\mc F$ is an $\al$-uniform family on $M$.
\item $\mc F^\mr{max}$ is an $\al$-uniform front on $M$ such that $\mc F=\overline{\mc F^\mr{max}}=(\mc
F^\mr{max})^{\ip}$.
\end{enumerate}
\end{enumerate}
\fprop
\prue
(a) is proved by a simple inductive argument on $\al$.  To prove that (b.1) implies (b.2), one observes first
that $(\mc F^\mr{max})^\ip=\mc F$, because $\mc F$ is compact, and then again use an inductive argument. The
proof of that (b.2) implies (b.1) one uses the well-known fact that if $\mc B$ is a uniform front, then
$\overline{\mc B}=\mc B^\ip$  (see for example \cite{ArTo}).
\fprue
\defi
Given two families $\mc F$ and $\mc G$ on $\om$ their sum and  product are defined by
\begin{align*}
\mc F\oplus \mc G:=& \conj{ s\cup t }{s<t,\, s\in \mc G\text{ and }t\in \mc F},\\
\mc F\otimes \mc G:=& \conj{  \bigcup_{i<n} s_i  }{   \{s_i\}_i\con \mc F,\, \max s_{i}<\min s_{i+1}, \, i<n,  \text{ and } \{\min s_i\}_i \in \mc
G}.
\end{align*}
\fdefi
The following are well-known facts of uniform fronts, and that are extended to uniform families by using the
previous proposition.  For more information on uniform fronts, we refer to  \cite{Lo}, \cite{LoTo3} and
\cite{ArTo}.
\prop\label{ioufhdfhvvbhjss}
\begin{enumerate}[(a)]
\item The rank of  an $\al$-uniform family is $\al$.
\item The unique $n$-uniform family on $M$, $n<\om$, is $[M]^{\le n}$.
\item If $\mc F$ is uniform and $t$ is a finite set, the $\mc F_t:=\conj{s\con \om}{t<s\text{ and }t\cup s\in \mc F}$ is uniform on $\om/t$. 
\item If $\mc F$ is an $\al$-uniform family  on $M$, then $\mc F\rest N$ is an $\al$-uniform family on $N$ for every $N\con M$
infinite. Consequently,  if $\mc F$ is an $\al$-uniform family on $\om$, then $\mc F$ is $\al$-homogeneous with $\mr{srk}(\mc F)=\mr{rk}(\mc F)=\al$. 
\item  If $\mc F$ is an $\al$-uniform family on $M$, and $\theta: M\to N$ is an order-preserving
bijection, then  $\conj{\theta"(s)}{s\in \mc F}$ is an $\al$-uniform family on $N$.
\item  Suppose that $\mc F$ and $\mc G$ are $\al$ and $\be$ uniform families on $M$, respectively. Then $\mc F\cup \mc G$,   $\mc F\oplus \mc
G$, $\mc F\sqcup \mc G$ and $\mc F\otimes \mc G$ are $\max\{\al,\be\}$, $\al+\be$, $\al \dot{+} \be$  and
$\al\cdot \be$-uniform families on $M$, respectively.
\item  Uniform fronts have the \emph{Ramsey property}: if $c:\mc F\to n$ is a coloring of a uniform front on $M$, then there is $N\con M$ infinite such that $c$ is constant on $\mc F\rest N$.
\item  Suppose that $\mc F$ and $\mc G$ are uniform families on $M$. Then there is some $N\con M$ such that either $\mc F\rest N\con \mc G$ or $\mc G\rest N\con \mc F$.  Moreover, when $\mr{rk}(\mc F)<\mr{rk}(\mc G)$,
 the first alternative
must hold and in addition $(\mc F\rest N)^\mr{max}\cap (\mc G\rest N)^{\mr{max}}=\buit$.
\item  If $\mc F$ is a uniform family on $M$, then there is $N\con M$ infinite such that $\mc F\rest N$ is
hereditary.
\item  If $\mc F$ is compact and $\ip$-hereditary family  on $\om$, then there is $M\con \om$ infinite such that $\mc F\rest M$ is a uniform family on $M$.   \qed
\end{enumerate}
\fprop
\nota \label{jkjuy3432}
\begin{enumerate}[(i)]
\item  The only new observation in the previous proposition is the fact in (e) that states that unions and square unions of uniform
families is a uniform family, and that can be easily proved by induction on the maximum of the ranks using, for example, that
$(\mc F \otimes \mc G)_{\{n\}} = (\mc F \otimes \mc G_{\{n\}}) \oplus \mc F_{\{n\}}$.

\item  A simple inductive argument shows that for every countable $\al$ and every  infinite $M\con \om$ there
is an $\al$-uniform family $\mc F$ on $M$, and, although uniform families are not necessarily
hereditary using (d) and (i) one can build them being hereditary.
\end{enumerate}
\fnota

We obtain the following consequence for families on an arbitrary partial ordering.
\cor\label{i89poopopwq92}
Suppose that $\mc P=(P,<_P)$ is a partial ordering, and suppose that $\mc F$ and $\mc G$ are compact and
hereditary families with $\mr{rk}(\mc F)<\mr{srk}(\mc G)$. Then every infinite chain $C$  of $\mc P$ has an
infinite subchain $D\con C$ such that $\mc F\rest D\con \mc G$.
\fcor
\prue
Fix $\mc F$, $\mc G$ and $C$ as in the statement. By going to a subchain of $C$, we may assume that $C$ is countable. Since $(C,<_P)$, so we can fix a bijection $\theta:C\to \om$, and we set $\mc F_0:=\conj{\theta"(s)}{s\in \mc F\rest C}$ and $\mc G_0:=\conj{\theta"(s)}{s\in \mc G\rest C}$. Observe that $\mc G_0\rest M$ has rank at least $\mr{srk}(\mc G)>\mr{rk}(\mc F_0)$ for every infinite $M\con \om$. Then by Proposition \ref{ioufhdfhvvbhjss} (j), (h) there is some infinite  $M\con \om$  such that $\mc F_0\rest M\con \mc G_0$, so $\mc F\rest \theta^{-1}(M) \con \mc G$.  
\fprue

Among uniform families, the \emph{generalized Schreier families} have been widely studied and used
particularly in Banach space theory. They have an algebraic definition and have extra properties, as for
example being spreading. Also, they have a sum-indecomposable rank. We recall the definition now.
\defi
\label{ioioijjweew} The \emph{Schreier family}  is
 \begin{equation*} \mc S:=\conj{s\con \om}{\#s\le \min s}.
\end{equation*}
A \emph{Schreier sequence} is defined inductively for $\al<\om_1$ by
\begin{enumerate}[(a)]
\item $\mc S_0:=[\om]^{\le 1}$,
\item  $\mc S_{\al+1}:=\mc S_\al \otimes \mc S$ and
\item $\mc S_\al:= \bigcup_{n<\om}( \mc S_{\al_n}\rest \om \setminus n)$ where $(\al_n)_n$ is such that  $\sup_n \al_n=\al$,  if $\al$ is limit.
\end{enumerate}
\fdefi

Note that the family $\mc S_\al$ depends on the choice of the sequences converging to limit ordinals.   
\defi[Spreading families]
A family $\mc F$ on $\om$ is  \emph{spreading} when for every $s=\{m_0<\dots<m_k\}\in \mc F$ and
$t=\{n_0<\cdots <n_k\}$ with $m_i\le n_i$ for every $i\le k$ one has that $t\in \mc F$.
\fdefi
\prop \label{iu4uieruieew44}
\begin{enumerate}[(a)]
\item 
Suppose that $\mc F$ and $\mc G$ are spreading families on $\om$. Then $\mc F\cup \mc G$, $\mc F\oplus \mc
G$ and $\mc F\otimes \mc G$ are spreading. If in addition $\mc F$ or $\mc G$ is hereditary, then $\mc F\sqcup \mc G$ is also spreading. 
\item If $\mc F$ is compact, hereditary and spreading, then $\mc F$ is $\mr{rk}(\mc F)$-homogeneous.
\item Suppose that $\mc F$ is a compact, hereditary and spreading family on $\om$ of finite rank $m$. Then there is some $n<\om$ such that $[\mc \om\setminus n]^{\le m}\con \mc F\con [ \om]^{\le m}$.
\end{enumerate}
\fprop
\prue
(a) is easy to verify. For (b): Given an infinite set $M\con \om$, let $\theta:\om\to M$ be the unique increasing enumeration of $M$. Then $s\mapsto \theta" s$ is a 1-1 and continuous mapping from $\mc F$ into $\mc F\rest M$, so $\mr{rk}(\mc F\rest M)\ge \mr{rk}(\mc F)\ge \mr{rk}(\mc F\rest M)$.     (c):  If there is some $s\in \mc F$ of cardinality at least $m+1$, then  $[\om\setminus \max s]^{\le m+1}\con \mc F$ because $\mc F$ is hereditary and spreading. And this is impossible, as it implies that the rank of $\mc F$ is at least $m+1$.  On the other hand, since the rank of $\mc F$ is $m$, there must be some element of it of cardinality $m$. Pick an element $s$ of $\mc F$ of such cardinality. Then $[\om\setminus \max s]^{\le m}\con \mc F$.  
\fprue
 The generalized Schreier families  are uniform families and they have extra
properties.
\prop \label{oi32o8980u233232}
\begin{enumerate}[(a)]
\item $\mc S_\al$ is hereditary, spreading and   $\om^\al$-uniform.
\item  For every  $\al\le \be$ there is $n<\om$ such that $\mc S_\al \rest ( \om \setminus n)\con \mc
S_\be$.
\end{enumerate}
\fprop
\prue
(a) The first two properties are well-known. The proof of that $\mc S_\al$ is a $\om^\al$-uniform family is
done by induction on $\al$. The case $\al=0$ is trivial, while it is easy to verify that $\mc S$ is an
$\om$-uniform family, so $\mc S_{\al+1}=\mc S_\al \otimes \mc S$ is a $\om^\al \cdot \om=
\om^{\al+1}$-uniform family by Proposition \ref{ioufhdfhvvbhjss} and inductive hypothesis. Suppose that $\al$
is limit. For a given $m,n\in \om$, let $\al_m^n<\om^{\al_n}$ be such that $(\al_m^n)_m$ is increasing,
$\sup_m \al_m^n =\al_n$ and $(\mc S_{\al_n})_{\{m\}}$ is a $\al_m^n$ uniform family. Since for every $m\in
\om$ we have that
$$(\mc S_\al)_{\{m\}}= \bigcup_{n\le m} (\mc S_{\al_n})_{\{m\}} \rest( \om  \setminus m)$$
it follows from Proposition \ref{ioufhdfhvvbhjss} (e) that $(\mc S_\al)_{\{m\}}$ is a $\be_m:=\max_{n\le
m}\al_m^n$-uniform family on $\om/m$. It is easy to see that $(\be_m)_m$ is increasing and satisfies that
$\sup_m \be_m =\om^\al$. (b) is proved by a simple inductive argument.
\fprue

\defi
Let $\mk S$ be the collection of  all hereditary, spreading uniform families on $\om$.
\fdefi

\prop\label{alpha_uniform_omega}
For every $\al<\ou$ there is a  hereditary, spreading $\al$-uniform family on $\om$.
\fprop
\prue
Let $(\mc S_\al)_{\al<\om}$ be a  Schreier sequence, and given a countable ordinal $\al$ with normal Cantor
form  $\al=\sum_{i\le k} \om^{\al_i}\cdot n_i$  we define
$$\mc F_\al:= (\mc S_{\al_0} \otimes [\om]^{\le n_0}) \oplus \cdots \oplus( \mc S_{\al_k}\otimes [\om]^{\le n_k}).$$
Then each $\mc F_\al$ is a hereditary, spreading $\al$-uniform family on $\om$.
\fprue

We present now  the concept of basis, which is a collection of families that intends to generalize the collection of uniform family on $\om$, and the multiplication $\otimes$ between the families in the basis.  It seems that there is no canonical definition for the multiplication $\mc F\times \mc G$ of two families on an index set $I$. However, when $\mc G$ is a family on $\om$ we can define it quite naturally as follows.
\defi
Let $\mc F$ be a homogeneous family on chains of a partial ordering $\mc P$, and let $\mc H$ be a homogeneous
family on $\om$. We say that a family $\mc G$ on chains of $\mc P$ is a \emph{multiplication of
}$\mc F$ by $\mc H$ when
\begin{enumerate}[(M.1)]
\item $\mc G$ is homogeneous and $\iota(\mr{srk}_\mc P(\mc G))=\iota(\mr{srk}_\mc P(\mc F)\cdot \mr{srk}(\mc
H))$.
\item   Every sequence $(s_n)_{n<\om}$ in $\mc F$ such that $\bigcup_n s_n$ is a chain of $\mc P$ has an  infinite subsequence $(t_n)_n$ such that for every $x\in \mc H$ one has that  $\bigcup_{n\in x}t_n\in \mc G$.
\end{enumerate}
\fdefi

\eje
\begin{enumerate}[(i)]
\item $\mc F\boxtimes_\mc P n$ is a multiplication of any homogeneous family $\mc F$ by $[\om]^{\le n}$. In general, given a family $\mc H\in \mk S$ of finite rank $n$, $\mc F\boxtimes_\mc P n$ is a multiplication of $\mc F$ by $\mc H$ (see Proposition  \ref{iu4uieruieew44} (b)). 
\item  If $\mc P$ does not have any infinite chain, then any homogeneous family $\mc F$
on chains of $\mc P$ has finite rank and given any homogeneous family
$\mc H$ on $\om$, $\mc G = \mc F$ satisfies (M.2) .
\end{enumerate}  \feje
Notice that always $\mc F\con \mc G$ for every multiplication $\mc G$ of $\mc F$ by any family $\mc H\neq
\{\buit\}$,  and that when $\mc F,\mc H\neq \{\buit\}$, then {\it(M.1)} is equivalent to $\iota(\mr{srk}_\mc P(\mc G))=\max\{\iota(\mr{srk}_\mc P(\mc F)),\iota( \mr{srk}(\mc
H))\}$, because $\iota(\al\cdot \be)=\max\{\iota(\al),\iota(\be)\}$  if $\al,\be\ge 1$.  When the family $\mc F=[\ka]^{\le 1}$ and $\mc H$ is the Schreier family $\mc S$, then the
existence of a family $\mc G$ satisfying (M.2) is equivalent to $\ka$ not being $\om$-Erd\H{o}s (see
\cite{LoTo}, and the remarks after Theorem \ref{firstMahlo}). Let us use the following notation. Given a
collection $\mk C$ of families on  chains of $\mc P$ and $\al<\ou$ let $\mk C_\al:=\conj{\mc F\in \mk
C}{\srk{\mc F}{\mc P}=\al}$.

\defi[Basis]
 Let $\mc P=(P,<)$ be a partial ordering with an infinite chain. A \emph{basis (of homogeneous families) on chains of
$\mc P$} is a pair $(\mk B, \times)$  such that:
\begin{enumerate}[(B.1)]
\item  $\mk B$ consists of homogeneous families on chains of $\mc P$, it contains all cubes $[P]^{\leq n}_{\mc P}$, and      $\mk B_\al\neq \buit$ for
all $\om\le \al<\ou$.
\item  $\mk B$ is closed under $\cup$ and $\sqcup_\mc P$, and if $\mc F\con \mc G\in \mk B$ is homogeneous on chains of $\mc P$   and  $\iota(\srk{\mc F}{\mc P})= \iota(\srk{\mc G}{\mc P})$, then $\mc F\in \mk B$.

\item   $\times : \mk B \times \mk S\to \mk B$ is such that for every $\mc F\in \mk B$ and every $\mc H\in \mk S$ one has that $\mc F\times \mc H$  is a multiplication of $\mc F$ by $\mc H$.
\end{enumerate}
When $\mc P = (P,<)$ is a total ordering, we simply say that $\mk B$ is a \emph{basis of families on $P$}.
\fdefi
Multiplications $\mc F\times \mc H$ for families $\mc H$ of finite rank always exist, so $(B.3)$ is equivalent to the existence of multiplication  by families $\mc H\in \mk S$ of infinite rank.   Also, $\{0\}\times \mc H=\{\buit\}$ is always a multiplication. The condition on $(B.2)$ concerning subfamilies is a technical requirement, but is not essential (see Proposition \ref{iweiijoweee44}). 

\prop
There is  basis of families on $\om$.
\fprop
\prue
Let $\mk B$ be the collection of all homogeneous families $\mc F$ such that there is some uniform family $\mc G$ with $\mc F\con \mc G$ and $\iota(\mr{srk}(\mc F))=\iota(\mr{srk}(\mc G))$.   
Given $\{\buit\}\cones \mc F\in \mk B$ and $\mc H\in \mk S$,  choose some uniform family $\mc G$ with $\mc F\con \mc G$ and $\iota(\mr{srk}(\mc F))=\iota(\mr{srk}(\mc G))$, and define
$$\mc F\times_\om \mc H:= (\mc G\otimes \mc H) \oplus \mc G.$$ 
{\it (B.1)}:   $\mk B_\al\neq \buit$  follows from Proposition \ref{alpha_uniform_omega}. {\it(B.2)}:  $\mk B$ is closed under $\cup$ and $\sqcup$ by Proposition \ref{uygfuyhgu1} and Proposition \ref{ioufhdfhvvbhjss}.  {\it(B.3)}: Fix $\{\buit\}\cones \mc F\in \mk B$, and $\mc H\in \mk G$. Choose an homogeneous family $\mc G$ defining  $\mc F\times_\om \mc H=(\mc G\otimes \mc H)\oplus \mc G$. We verify first (M.1): Let $\al,\be<\ou$ be such that $\mc G$ and $\mc H$ are $\al$-uniform and $\be$-uniform, respectively. Observe that $\iota(\mr{skr}(\mc F))=\iota(\mr{srk}(\mc G))=\iota(\al)$ and $\mr{srk}(\mc H)=\be$. Then $\mc F\times_\om H$ is $\al\cdot (\be+1)$-uniform, $\mr{srk}(\mc F\times_\om \mc H)=\al(\be+1)$, so $\iota(\mr{srk}(\mc F\times_\om \mc H))=\max\{\iota(\al),\iota(\be)\}=\max\{\iota(\mr{srk}(\mc F)),\iota(\mr{srk}(\mc H))\}=\iota(\mr{srk}(\mc F)\cdot \mr{srk}(\mc H))$.  This finishes the proof of (M.1). We check now 
(M.2):  Suppose that $(s_k)_k$ is a sequence in $\mc F$. Let $(t_k)_{k<\om}$ be a $\De$-subsequence with root $t\in
\mc F$ such that $t<t_k\setminus t <t_{k+1}\setminus t$ for every $k$. Suppose that $x\in \mc H$.   Then
$\{\min t_k \setminus t\}_{k\in x}\in \mc H$, because $\mc H$ is spreading. Hence, $\bigcup_{k\in x}
(t_k\setminus t) \in  \mc G\otimes \mc H$.  Since $t<\bigcup_{k\in x} (t_k\setminus t)$, it follows that
$\bigcup_{k\in x} t_k = t \cup \bigcup_{k\in x} (t_k\setminus t)\in (\mc G\otimes \mc H)\oplus \mc G$.
\fprue

Our main result is the following:

\teor\label{firstMahlo}
Every cardinal $\theta$ strictly smaller  than the first Mahlo cardinal  has a basis of families on $\theta$.
\fteor

The proof, in Section \ref{bskdjgbsekjbgcnamnc}, is done inductively on $\ka$ and using trees with small height and levels (to be precised later)
in order to step up. For example, 
when $\ka$ is not strong limit, there
must be $\la<\ka$ such that $2^\la\ge \ka$, so, by the inductive hypothesis, there must be a basis on $\la$, and
there is a rather natural way to lift it up to a basis on the nodes of the complete binary tree, via the height mapping. The case when
$\ka$ is not regular is similar. When on the contrary $\ka$ is inaccessible, the tree $T$ of cardinality $\kappa$ is
substantially more complicated, and in fact relies on the method of walks in ordinals. In any of these cases, the
main difficulty is to pass from a basis on chains and a basis ``on the antichains'' of a tree $T$ to a basis of families
on the nodes of $T$ with respect to any total ordering on them, thus getting a basis on $|T|$ (see Theorem \ref{4jtierjteoirtertrtr344}).

Recall that $\ka$ is $\om$-Erd\H{o}s when for every coloring $c:[\ka]^{<\om}\to 2$ there is an infinite
$c$-homogeneous subset $A\con \ka$, that is, for $s,t\in A$, $c(s)=c(t)$ when $\#s=\#t$. A compact and
hereditary family on $\ka$ is called \emph{large} when $\mr{srk}(\mc F)\ge \om$, or, equivalently, when $\mc
F$ satisfies (M.2) for $[\ka]^{\le 1}$ and the Schreier family $\mc S$. It is proved in \cite{LoTo} that the
existence of such families in $\ka$ is equivalent to $\ka$ not being $\om$-Erd\H{o}s.

\begin{pro}
Characterize when $\ka$ has $\ge \om$-homogeneous families.
\end{pro}

\begin{pro}
Characterize the cardinal numbers $\ka$  such that that there exists $c:[\ka]^{<\om}\to 2$ and an $\ge
\om$-homogeneous family on $\ka$ such that every $s\in \mc F$ is $c$-homogeneous and such that for every
sequence
 $(s_n)_{n<\om}$ in $\mc F$ and every $l<\om$ there are $n_1<\cdots< n_l$ such that $\bigcup_{i=1}^ls_{n_i}$
 is $c$-homogeneous.
\end{pro}

The first such $\ka$ not satisfying this coloring property is at least the first Mahlo cardinal and smaller
than the first $\om$-Erd\H{o}s cardinal. 

Before going into the particular case of the partial ordering being a tree, we present some operations of partial orderings and results guaranteeing
that we can transfer families or bases from some partial orderings to more complex ones (Subsection \ref{transfer_bases}). The following characterization of the existence of a basis will be useful.

\prop   \label{iweiijoweee44}
A partial ordering $\mc P$ with an infinite chain has a basis if and only if there is a pair $(\mk B,\times)$, called \emph{pseudo-basis}, with the following properties:
\begin{enumerate}[(B.1')]
\item  $\mk B$ consists of homogeneous families on chains of $\mc P$, it contains all cubes, and for every $\om \leq \alpha < \ou$ there is $\mc F\in \mk B$ such that $\alpha \leq \mr{srk}_\mc P(\mc F) \leq \iota(\alpha)$.
\item $\mk B$ is closed under $\cup$ and $\sqcup_\mc P$.
\item   $\times : \mk B \times \mk S_{\ge \om}\to \mk B$, $\mk S_{\ge\om}$ being the infinite ranked families of $\mk S$, is such that for every $\mc F\in \mk B$ and every $\mc H\in \mk S$ one has that $\mc F\times \mc H$  is a multiplication of $\mc F$ by $\mc H$.
\end{enumerate}
\fprop
\prue
 Suppose that $(\mk B,\times)$ satisfies (B.1'), (B.2') and (B.3'). Let $C=\{p_n\}_n$ be an infinite chain of $\mc P$, of order type $\om$. Fix a basis $(\mk B(\om),\times_\om)$ of families on $\om$.  For each $\mc G\in \mk B(\om)$, let $\bar{\mc G}:= \conj{\{p_n\}_{n\in x}}{x\in \mc G}$.  Then $\bar{\mc G}$ is homeomorphic to $\mc G$. Given $\mc F \in \mk B$, let $\widetilde{\mc F}: = \conj{s \in \mc F}{s \cap C = \emptyset}$. Now let $\mk B'$ be the collection of all unions
 $\widetilde{\mc F} \cup \bar{\mc G}$ such that $\mc F\in \mk B$, $\mc G\in \mk S$, and finally let $\mk B''$ be the collection of all $\mc P$-homogeneous families $\mc F$ such that there is some $\mc G\in \mk B'$ with $\mc F\con \mc G$ and $\iota(\mr{srk}_\mc P(\mc F))=\iota(\mr{srk}_\mc P(\mc G))$.  For each $\mc F\in \mk B''$ we choose $\mc G_\mc F\in \mk B$ and $\mc H_\mc F\in \mk B(\om)$ such that $\mc F\con  \widetilde{\mc G_\mc F}\cup \bar{\mc H_\mc F} \in \mk B'$ and $\iota(\mr{srk}_\mc P(\mc F))=\iota(\mr{srk}_\mc P( \widetilde{\mc G_\mc F}\cup \bar{\mc H_\mc F}))$. If $\mc H$ has infinite rank,  we define $\mc F\times' \mc H:= (\widetilde{\mc G_F \times  \mc H} ) \cup \bar{(\mc H_\mc F\times_\om \mc H)}$, and when $\mc H$ has finite rank $n$, let $\mc F\times' \mc H:= \mc F \boxtimes_\mc P n$.   It is easy to check that $(\mk B'',\times')$ is a basis on chains of $\mc P$.
\fprue

\subsection{Ranks and uniform families}

The objective of this part is to give technical tools for getting upper bounds for ranks with the use of uniform families.


\defi
Given families $\mc F$ on $\om$ and $\mc G$ on a partial ordering $\mc P$, we say that a mapping $f:\mc F\to \mc G$ between two families is $(\ip,\con)$-increasing when $s\ip t$
implies that $f(s)\con f(t)$.
\fdefi
The fact that a  point is in a certain derivative of a compact metrizable space $K$ can be witnessed by a
continuous and 1-1 mapping from a uniform family into $K$.

\prop[Parametrization of ranks] \label{njewirjiowejrew}
Suppose that $K$ is  a compact metrizable space and let $\al<\ou$. Then a point $p\in K$  is such that $p\in
K^{(\al)}$ if and only if for every  $\al$-uniform family $\mc B$ there is a 1-1 and continuous function
$\theta: \mc B\to K$ such that $\theta(\buit)=p$. In case $K=\mc F$ is  a compact family on $I$, $p\in \mc
F^{(\al)}$ if and only if for every $\al$-uniform family $\mc B$ there is a  1-1 and continuous  mapping
$\theta: \mc B\to  \mc F$  such that $p=\theta(\buit)$ and such that $\theta$ is $(\ip,\con)$-increasing.
\fprop

\prue
Given $p\in K$ and $\vep>0$, let $B(p,\vep)$ be the open ball around $p$ and radius $\vep$. The proof is by
induction on $\al$. Suppose that $p\in K^{(\al)}$ and let $\mc B$ be a $\al$-uniform family on $M$ and $\mc
C$ the collection of $\ip$-maximal subsets in $\mc B$.  Without loss of generality we assume that $M=\om$.
Let  $\al_n<\al$ be such that $\mc C_{\{n\}}$ is $\al_n$-uniform on $\om/n$.   Choose  $(p_n)_n$ in
$K^{(\al_n)}$ converging non-trivially to $p$ such that there are mutually disjoint closed balls $B_n$ around
$p_n$ with $\mr{diam}(B_n)\downarrow_n 0$.   Since each $p_n\in K^{(\al_n)}$ it follows by inductive
hypothesis that  for each $n$ there is a  1-1 and continuous function
$$\theta_n:  \overline{\mc B_{\{n\}}} =\mc C_{\{n\}}  \to B_n$$  with   $\theta_n(\buit)=p_n$. Let $\theta: \mc B\to K$ be defined by $\theta(\buit)=p$, $\theta(s):= \theta_{\min s}(s\setminus \{\min s\})$, for $s\neq \buit$.  By the choice of the balls $B_n$ it follows that $\theta$ is 1-1. We verify now that $\theta$ is continuous:
Suppose that $(s_k)_k$ tends to $s$.  Suppose first that $s\neq \buit$, let $n:=\min s$. Then there is $k_0$
such that for every $k\ge k_0$ one has that $\min s_k=n$. It follows that for every $k\ge k_0$,
$\theta(s_k)=\theta_n(t_k)$ and $\theta(s)=\theta_n(t)$, where  $t_k:=s_k\setminus \{n\}$ and $t:= s\setminus
\{n\}$. Hence, $\lim_{k\to \infty} \theta(s_k)=\lim_{k\to \infty}\theta_n(t_k)=\theta_n(t)=\theta(s)$.
Suppose now that $s=\buit$.  Fix $\ga>0$ and suppose that $d(p, \theta(s_k))\ge \ga$ for every $k$ belonging
to an infinite subset $M\con \om$. Without loss of generality, we may assume that $(s_k)_{k\in M}$ is a
$\De$-system with empty root such that $s_k <s_l$ if $k<l$ in $M$.      Since   $\theta(s_k)\in B_{n_k}$, for
$n_k:=\min s_k$   for every $k$, and since $(n_k)_{k\in M}$ tends to infinity, $(p_{n_k})_{k\in M}$ converges
to $p$, so that there is some $k$ such that $d(p, \theta(s_k))<\ga$, a contradiction. The reverse implication
is trivial.

Now, if $\mc F$ is a compact  family on $I$ and $\al$ is a countable ordinal, then $p\in \mc F^{(\al)}$ if
and only if $p\in (\mc F\rest I_0)^{\al}$ for some countable subset $I_0$. With a small modification of the proof we have just exposed  applied to the compact and metrizable space $K=\mc F\rest I_0$ one can find recursively on $\al$ a 1-1 and
continuous $\theta:\mc B\to \mc F\rest I_0$ which in addition is $(\ip,\con)$-increasing.
\fprue


\lema
\label{778tr5643367898}
 Suppose that $\mc B$ and $\mc C$ are uniform families,  $\mc F$ is a compact  family on  some index set $I$
  with $\mr{rk}(\mc F)<\mr{rk}(\mc C)$. Suppose that $\la: \mc B\otimes \mc C\to \mc F$ is $(\ip,\con)$-increasing. Then there is
  a  finite subset $x$ of $\om$ and  some infinite set $x<M$ such  that $\{x\}\sqcup \mc B\rest M\con \mc B\otimes \mc C$ and such that $\la  $ is constant on
   $\{x\}\sqcup (\mc B\rest M)^\mr{max} $.  If in addition $\la$ is continuous, then $\la$ is constant on $\{x\}\sqcup (\mc B\rest M) $.
 \flema
\prue
Let $\mc B_0:=\mc B^{\mr{max}}$, $\mc C_0:=\mc C^{\mr{max}}$, and for each $s\in \mc B_0\otimes \mc C_0$, let
$s=\bigcup_{i\le k_s} s_i$ be the canonical decomposition of $s$; i.e. $s_i<s_{i+1}$ are in $\mc B_0$ and
$\{\min s_i\}_{i\le k_s}\in \mc C_0$.
\clam
There is an  infinite subset $M\con \om$ such that one of the following conditions hold.
 \begin{enumerate}[(a)]
  \item For every $s_0<\cdots <s_{k-1}$ in $\mc B_0\rest M$ with $\{\min s_i\}_{i<k}\in \mc C\rest M\setminus \mc C_0$ and every $s_{k-1}<x<y\in \mc B_0\rest M$ one has that   $\la(\bigcup_{i<k}s_i\cup x) \neq\la(\bigcup_{i<k}s_i\cup y)$.
\item  For every $s=\bigcup_{i\le k_s}s_i\in \mc B_0\rest M\otimes \mc C_0\rest M$ there is $k<k_s$  such that for every  $s_{k-1}<x<y\in \mc B_0\rest M$ one has that   $\la(\bigcup_{i<k}s_i\cup x)=\la(\bigcup_{i<k}s_i\cup y)$.
\end{enumerate}
\fclam
\prucl
We use the   the Ramsey property of uniform fronts and a diagonal argument. We can find a decreasing sequence $(N_j)_j$ of infinite subsets of $\om$, such that $n_{j}:=\min N_j< n_{j+1}:=\min N_{j+1}$ and such that for every $j$ and every $s_0<\cdots<s_{k-1}$ in $\{n_i\}_{i\le  j}$ with $\{\min s_i\}_{i< k}\in \mc C\setminus \mc C_0$ one has that  either
\begin{enumerate}[(i)]
\item for every $x<y$ both in $\mc B_0\rest N_{j+1}$ one has that $\la(\bigcup_{i<k}s_i \cup x)\neq \la(\bigcup_{i\le j}s_i \cup y)$,  or  
\item for every $x<y$ both in $\mc B_0\rest N_{j+1}$ one has that  $\la(\bigcup_{i<kj}s_i \cup x)= \la(\bigcup_{i\le j}s_i \cup y)$.
 \end{enumerate} 
This is done inductively: suppose that $(N_i)_{i\le j}$ is already defined, $n_i:=\min N_i$. Given a sequence  $s_0<\cdots<s_{k-1}$ in $\{n_i\}_{i\le  j}$ with $\{\min s_i\}_{i< k}\in \mc C\setminus \mc C_0$, we can define the coloring $c: (\mc B\oplus\mc B)\rest (N_{j}/ s_{k-1})\to 2$,  for $x<y$  each in $\mc B\rest (N_{j}/ s_{k-1}$    by $c(x\cup y)=0$ if $\la(s_0\cup \cdots \cup s_{k-1}\cup x)\neq \la(s_0\cup \cdots \cup s_{k-1}\cup y)$ and 1 otherwise. Notice that for every $x\in \mc B\rest (N_{j}/ s_{k-1}$ one has that $s_0\cup \cdots \cup s_{k-1}\cup x\in \mc B\otimes \mc C$, since $\{\min s_{i}\}_{i<k}\in \mc C \setminus \mc C_0$ and so $C_{\{\min s_{i}\}_{i<k}}$ is a uniform family (Proposition \ref{ioufhdfhvvbhjss}) of rank $>0$, so it contains all singletons $>\min s_{k-1}$. 

Let now $d:\mc B\otimes \mc C\rest N\to 2$ be defined by $c(\bigcup_{i\le k_s}s_i)=0$ if there is $k<k_s$ such that $\la(\bigcup_{i\le k} s_i)=\la((\bigcup_{i<k} s_i\cup s_{k+1})$ and $c(\bigcup_{i\le k_s}s_i)=1$ otherwise.  Let $M\con N$ be such that $d$ is constant on $\mc B\otimes \mc C\rest M$ with value $\ro\in 2$.  Then if $\ro=1$, then  (b) holds on $M$. If $\ro=0$, the   (a) holds on $M$. 
\fprucl
\clam
(b) above holds. 
\fclam
\prucl
Otherwise,   (a) holds.  We claim that there is some $\la(\buit)\con z \in \mc F^{(\al)}$, $\al:=\mr{rk}(\mc C)$,   which is impossible by hypothesis. The proof is by induction on $\al$.  let $(s_i)_{i<\om}$ be a sequence  in $\mc B_0\rest M$ such that $s_i<s_{i+1}$ for every $i$, $m_i:=\min s_i$,  and $N:=\{m_i\}_i\con M$. For each $i$, let $\la_i: \mc B\otimes \mc C_{\{m_i\}}\rest N\to \mc F$, $\la_i(x):=\la(s_i\cup x)$ for every $ x\in \mc B\otimes \mc C_{\{m_i\}}\rest N$.  Since $\la_i$ also satisfies (a), we obtain that for every $i$ there is some  $\la(s_i)\con z_i\in (\mc F)^{(\al_i)}$, where $\al_i=\mr{rk}(\mc C_{\{m_i\}}$.  By (a), $(\la(s_i))_{i}$ are pairwise different, so $(z_i)_i$ is non-trivial. Let $(z_i)_{i\in A}$ be a non-trivial $\De$-sequence with root $z\in \mc F$ ($\mc F$ is compact). Since $\mc C$ is uniform, it follows that either $\al=\be+1$ and  $\al_i=\be$ for every $i$, or $\al$ is limit and $\sup_i \al_i=\al$. In any case, $z\in \mc F^{(\al)}$. Since $\la$ is $(\ip,\con)$-increasing, it follows that $\la(\buit)\con \la(s_i)$ for every $i$, so $\la(\buit)\con z\in \mc F^{(\al)}$, as required. 
\fprucl
%
%
So,  (b) above holds. Fix $s=\bigcup_{i\le k_s}s_i\in \mc B_0\rest M\otimes \mc C_0\rest M$. Let $k<k_s$ be such that, setting $x:=\bigcup_{i<k}s_i$, then $\la(x\cup y)=\la(x\cup z)$ for every $x<y<z$ for every $y,z\in \mc B_0\rest M$. We claim that $\la(x\cup y)=\la(x\cup z)$ for every $x<y,z\in \mc B_0\rest M$.  Find $y,z<w\in \mc B_0\rest M$. Then $\la(x\cup y)=\la(x\cup w)=\la(x\cup z)$.

 If we assume that $\la$ is in addition continuous, since $\mc B\rest M$ is scattered, the set of isolated points is dense. Hence  $\{x\}\sqcup (\mc B\rest M)^\mr{max}$ is dense
  in $\{x\}\sqcup (\mc B\rest M)$. Since $\la$ is constant on $\{x\}\sqcup (\mc B\rest M)^\mr{max}$, it is constant on  $\{x\}\sqcup (\mc B\rest M)$.
\fprue

Our first upper estimation on ranks is the following.  
\prop
\label{oi43jiro43ro34njrtguyy73}
 Suppose that $\mc F$ and $\mc G$ are countably ranked     families  and suppose that $\la:\mc F\to \mc G$ is  $\con$-increasing. Then
 $$\mr{rk}(\mc F)< \sup_{t\in \mc G}(\mr{rk}(\conj{s\in \mc F}{\la(s)\con t})+1)   \cdot (\mr{rk}(\mc G)+1).$$
 If in addition $\la$ is continuous, then we obtain that
  $$\mr{rk}(\mc F)< \sup_{t\in \mc G}(\mr{rk}(\conj{s\in \mc F}{\la(s)=t})+1)   \cdot (\mr{rk}(\mc G)+1).$$
\fprop
 \prue
 Let $\al:= \sup_{t\in \mc G}(\mr{rk}(\conj{s\in \mc F}{\la(s)\con t}) +1)  $, $\be:=\mr{rk}(\mc G)$,
  and  suppose that $\mc F^{(\al\cdot (\be+1))}\neq \buit$. Let $\mc B$ and $\mc C$ be $\al$ and $\be+1$ uniform families,
let $f: \mc B\otimes \mc C\to \mc F$ be a 1-1, continuous and $(\ip,\con)$-increasing function, and let
$\theta:=\la \circ f$. By hypothesis, $\theta$ is $(\ip,\con)$-increasing, so   it follows from   Lemma
\ref{778tr5643367898} that there are $x \subset \omega$ finite and $x<M$ infinite such that
  $\{x\}\sqcup \mc B\rest M\con \mc B\otimes \mc C$ and  such that   $\theta$ is constant on $\{x\}\sqcup (\mc B\rest M)^\mr{max}$ with value $t\in \mc G$. This implies that the mapping $\theta_0(y):=f(x\cup y)$ for every
   $y\in \mc B\rest M$  is a 1-1 and continuous mapping $\theta_0: \mc B\rest M\to \conj{s\in \mc F}{\la(s)\con t}$, hence $\mr{rk}(\conj{s\in \mc F}{\la(s)\con t})\ge \al$, and this is impossible.

 Suppose that in addition $\la$ is continuous. Then  the desired result is proved similarly by changing the definition of $\al$ with $\al:=  \sup_{t\in \mc G}(\mr{rk}(\conj{s\in \mc F}{\la(s)= t}) +1)$, and then using the  particular case of continuous functions in Lemma  \ref{778tr5643367898}.
 \fprue

\subsection{Transfering families and bases with the use of operations}\label{transfer_bases}
We present several procedures to define bases from other bases, focused on transfer methods that will be used in Section \ref{bskdjgbsekjbgcnamnc} to step up bases on  a cardinal to bases on bigger cardinal numbers.  In particular we will see that if an infinite cardinal $\ka$ has a basis, then the complete binary tree $2^{\le \ka}$ has a basis consisting on chains of $2^{\le \ka}$. We will also give more involved transfer methods that   will  be used to understand more complicated trees (see Subsection \ref{trees_walks}).  We start with some terminology.

\defi
Let $\mc P=(P,\le_P)$ and $\mc Q=(Q,\le_Q)$ be partial orderings and  $\la: P\to Q$.
\begin{enumerate}[(i)]
\item  $\la$ is
\emph{chain-preserving} when $p_0\le_P p_1$ implies that  $\la(p_0)\le_Q \la(p_1)$ or $\la(p_1)\le_Q
\la(p_0)$.
\item   $\la$ is \emph{1-1 on chains} when $\la \rest C$ is 1-1 for every chain $C$
of $\mc P$.
\item  $\la$ is \emph{adequate} when it is chain-preserving and 1-1 on chains.
\end{enumerate}

\fdefi

In other words, $\la$ is chain-preserving if it is a \emph{graph homomorphism} between the corresponding
\emph{comparability graphs}. Observe that $\la$ is chain-preserving if and only if $\la"(C)$ is a chain of
$\mc Q$ for every chain $C$ of $\mc P$. Observe also that when $\mc Q$ is a total ordering, every mapping
$\la:P\to Q$ is chain preserving. The previous proposition can be this generalized as follows:

\teor\label{ui32iuui32t744}
Let  $\mc P$ and $\mc Q$ be partial orderings which have infinite chains, and let $\la: \mc P\to \mc Q$  be adequate. If   $\mc Q$ has a
basis of homogeneous families, then  so has $\mc P$.
\fteor

In order to prove this, we introduce the following operation.

\defi[preimage]
Given partial orderings $\mc P$ and $\mc Q$, $\la:P\to Q$ and a family   $\mc G$ on chains of  $\mc Q$, let
 \begin{align*}
 \la^{-1}(\mc G):=&\conj{s\con P}{s\text{ is a chain of $\mc P$ and } \la"s\in \mc G}.
 \end{align*}
\fdefi


\lema \label{kloi89787855}
Suppose that $\mc P$ and $\mc Q$ are two partial orderings,  $\la:\mc P\to \mc Q$   is  adequate. Suppose also that $\mc G$ is a  family on chains of $\mc Q$.
\begin{enumerate}[(a)]
\item If  $\mc G$ is pre-compact, hereditary, then so is $\la^{-1}(\mc G)$.
\item   If $\mc G$ is  countably ranked, then
\begin{align}
\mr{rk}(\la^{-1}(\mc G))  < &  \om \cdot (\mr{rk}(\mc G)+1).\label{uiui7656}
\end{align}
Consequently, if $\mc P$ has infinite chains, and $\mc G$ is $(\al,\mc Q)$-homogeneous, $\al\ge \om$,   then $\la^{-1}\mc G$ is $(\be,\mc P)$-homogeneous with $\al\le \be <\iota(\al)$.
\end{enumerate}
\flema
\prue

 Set $\mc F:= \la^{-1}(\mc G)$. It is clear that $ \mc F$ is
hereditary when $\mc G$ is hereditary. Suppose that  $\mc G$ is pre-compact.  Let $(x_n)_n$ be a sequence in
$\mc F$. W.l.o.g. we assume that $(x_n)_n$ converges to $A\con P$. Limit of chains are chains, so $A$ is a
chain of $P$. The proof will be finished when we verify that $A$ is finite.  We  assume that $(\la"x_n)_n$ is
a $\De$-sequence with root $y$.  It is easy to see that  $\la"A \con y$. Since $\la$ is 1-1 on chains, it
follows that   $\#A \le \# y$, so $A$ is finite.   Suppose that $\mc G$ has countable rank. We apply
Proposition  \ref{oi43jiro43ro34njrtguyy73} to $\la: \mc F\to \mc G$ to conclude that
\begin{equation}\label{oi8987444}
\mr{rk}(\mc F)<\sup_{y\in \mc G} (\mr{rk}(\conj{x\in \mc F}{\la(x)\con y})+1) \cdot (\mr{rk}(\mc G)+1).
\end{equation}
Observe that given $y\in \mc G$, since $\la$ is 1-1 on chains, it follows that
$$\conj{x\in \mc F}{\la(x)\con y} \con [P]^{\le \# y},$$
 so from      \eqref{oi8987444} we obtain the desired inequality in \eqref{uiui7656}.

Suppose that $\mc P$ has infinite chains and let $\mc G$ be $(\al, \mc Q)$-homogeneous.  Let us see that $\mc F$ is $(\be, \mc P)$-homogeneous with $\be\ge \al$.   Let $X$ be an infinite chain of $\mc P$ such that $\mr{srk}_\mc P(\mc F)=\mr{rk}(\mc F\rest X)=\be$. Then $Y:=\la"X$ is an infinite chain of $\mc Q$. Since $h: \mc F\rest X\to \mc G\rest Y$, $h(s):= \la" s$ is an homeomorphism,  it follows that $\mr{rk}(\mc G\rest Y) =\be$, hence $\mr{srk}_\mc Q(\mc G) \le \be$. On the other hand, it  follows from \eqref{uiui7656} and the fact that $\mc G$ is homogeneous that
\[ \pushQED{\qed} 
\mr{srk}_\mc P(\mc F) \le \mr{rk}(\mc F)<\om \cdot (\mr{rk}(\mc G)+1) <\iota(\mr{srk}_\mc Q(\mc G))\le \iota(\mr{srk}_\mc P(\mc F)).\qedhere
\popQED
\] 
\let\qed\relax
\fprue

\prue[{\sc Proof of Theorem \ref{ui32iuui32t744}}]

Let $(\mk C,\times_\mc Q)$ be a basis of families on chains of $\mc Q$. Let
$\mk B$ be the collection of all $\mc P$-homogeneous families $\la^{-1}\mc G$ with $\mc G\in \mk C$.
For each $\mc F\in \mk B$, choose $\mc G_\mc F$ such that $\mc F= \la^{-1}(\mc G_\mc F)$, and for $\mc H\in \mk S$, let $\mc F\times \mc H:= \la^{-1}(\mc G_\mc F\times_\mc Q \mc H)$.  We check that $(\mk B,\times)$  satisfies (B.1'), (B.2') and (B.3), which is enough to guarantee the existence of a basis on $\mc P$, by Proposition \ref{iweiijoweee44}. Given $\mc G\in \mk C_\al$, $\la^{-1}\mc G\in \mk B$ and by Lemma \ref{kloi89787855} we know that $\la^{-1}\mc G$ is $\be$-uniform with $\al \leq \be < \iota(\al)$. We check now (B.2').
Suppose that $\mc G_0,\mc G_0\in \mk C$. Then $\la^{-1}\mc G_0 \sqcup_\mc P\la^{-1}(\mc G_1)= \la^{-1}(\mc
G_0 \sqcup_\mc Q \mc G_1)$,  so $\la^{-1}\mc G_0 \sqcup_\mc P\la^{-1}(\mc G_1)\in\mk B $, because $\mc
G_0\sqcup_\mc Q\mc G_1\in \mk C$. Similarly one shows that $\mk B$ is closed under $\cup$.

Finally, we verify (B.3) for $(\mk B,\times)$. Fix $\mc F\in \mk B$ and $\mc H\in \mk S$. Then $\mc F\times \mc H=\la^{-1}(\mc G_\mc F \times_\mc Q\mc H)$.
\begin{align}
\iota(\mr{srk}_\mc P(\mc F\times \mc H))= & \iota(\mr{srk}_\mc Q(\mc G_\mc F\times_\mc Q\mc H))=\max\{\iota(\mr{srk}_\mc Q(\mc G_\mc F)), \iota(\mr{srk}(\mc H))\} \nonumber \\
=&\iota(\mr{srk}_\mc P(\mc F) \cdot \mr{srk}(\mc H)).\nonumber\end{align}
Let now $(s_n)_n$ be a sequence in $\mc F$ such that $C:=\bigcup_n s_n$ is a chain. Then $\la"C= \bigcup_{n}\la"s_n $ is a $\mc Q$-chain, and $\la"s_n \in \mc G_\mc F$. By the property (M.2) of $\times_\mc Q$, we obtain that there is a subsequence $(t_n)_n$ of $(s_n)_n$ such that $\bigcup_{n\in x} \la" t_n\in \mc G_\mc F\times_\mc Q\mc H$ for every $x\in \mc H$. This means that $\bigcup_{n\in x} t_n\in \mc F\times \mc H$ for such $x\in \mc H$.
\fprue

Recall that given a sequence of partial orderings $(\mc P_i)_{i\in I}$ we denote by $\biguplus_{i\in I}\mc P_i$ its \emph{disjoint union}, which is the partial ordering on $\bigcup_{i\in I} P_i\times\{i\}$ defined by $(p,i)<(q,j)$ if and only if $i=j$ and $p<_{\mc P_i} q$.

\prop\label{hjgdty42932833}
Suppose that $\theta$ is a regular cardinal number such that $\ou^{\ou}<\theta$, and suppose that every
$\xi<\theta$ has a basis on families on $\xi$. Suppose that $(\theta_\xi)_{\xi<\theta}$ is a sequence of
infinite ordinals such that $\sup_\xi \theta_\xi=\theta$. Then  the disjoint union of
$\biguplus_{\xi<\theta}(\theta_\xi,<)$ has a basis of families on chains.
\fprop
\prue
The proof is a counting argument. Set $\mc P:=\biguplus_{\xi<\theta} (\theta_\xi,<)$. First of all, let $C\con \theta$ be such that $\#C=\theta$ and $(\theta_\xi)_{\xi\in C}$ is strictly increasing with supremum $\theta$.  For each $\xi\in C$ let $(\mk C^\xi, \times_\xi')$ be a basis on $\theta_\xi\times\{\xi\}$. Let $F: C\to \ou^{\ou}$ be the mapping that to $\xi\in C$ and $\al<\ou$ assigns
$$F(\xi)(\al):=\min \conj{\mr{rk}(\mc F)}{\mc F\in \mk C^\xi_\al}<\iota(\al)<\ou.$$
Since $\ou^{\ou}<\theta$ there must be $D\con C$ of cardinality $\theta$ and $f\in \ou^{\ou}$ such that
$F(\xi)=f$ for every $\xi\in D$.  Define now for each $\xi<\theta$, $\mu_\xi:=\min\conj{\ga\in D}{\theta_\xi\le \theta_\ga}$.
Fix $\xi<\theta$. Let $\mk B^\xi$ be equal to $\mk C^\xi$ if $\xi\in D$, and let $\mk B^\xi$ be the collection of families $\conj{x\times \{\xi\}}{ x\con \theta_\xi \text{ and } x\times \{\mu_\xi\}\in \mc F} $ for $\mc F\in \mk C^{\mu_\xi}$. For $\xi\in D$, let $\times_\xi=\times_\xi'$. Suppose that $\xi\notin D$. For each $\mc F\in \mk B_\xi$, let $\mc G_\mc F$ be such that $\mc F=\conj{x\times \{\xi\}}{ x\con \theta_\xi \text{ and } x\times \{\mu_\xi\}\in \mc G_\mc F}$, and define
$$\mc F\times_\xi \mc H:= \conj{x\times \{\xi\}}{ x\times\{\mu_\xi\}\in \mc G_\mc F \times_{\mu_\xi} \mc H}.$$
It is easy to see that $(\mk B^\xi,\times_\xi)$ is a basis on $\theta_\xi\times \{\xi\}$ for every $\xi<\theta$.
Let $\mk B$ be the collection of all $\mc P$-homogeneous families $\mc F$ such that $\mc F \rest \theta_\xi \times \{\xi\}\in \mk B^\xi$.  Define for $\mc F\in \mk B$ and $\mc H\in \mk S$
$$\mc F\times \mc H:=\bigcup_{\xi<\theta}  (\mc F\rest (\theta_\xi \times \{\xi\}) \times_\xi \mc H).$$
We check that $(\mk B,\times)$ is a pseudo-basis on chains of $\mc P$. It is easy to see that $\mk B$ contains all finite cubes. Now  let $\om\le \al<\ou$  and we prove that $\mk B_\al\neq \buit$.  For each $\xi\in D$, let $\mc F_\xi\in \mk C^\xi_\al$.  Define for each $\xi<\theta$ $\mc G_\xi=\mc F_\xi$ if $\xi\in D$, and $\mc G_\xi:= \conj{x\times \{\xi\}}{x\con \theta_\xi \text{ and } x\times \{\mu_\xi\}\in  \mc F_{\mu_\xi}}$.  Notice that    each $\mc F_\xi$ is $\al_\xi$-uniform with $\al\le \al_\xi<\iota(\al)$. Notice also that $\sup_{\xi<\theta} \mr{rk}(\mc F_\xi) =f(\al)<\iota(\al)$.  Now let $\mc G:=\bigcup_{\xi<\theta}\mc G_\xi$. Since
$$\mr{rk}(\mc G)\le \sup_{\xi<\theta}\mr{rk}(\mc G_\xi) +1 =f(\al)+1<\iota(\al),$$
and $\mr{srk}_{\mc P}(\mc G)=\al$, it follows that $\mc G$ is $(\al,\mc P)$-homogeneous.  It is easy to see that $\mk B$ is closed under $\cup$ and $\sqcup$.  Now we prove that $\times$ is a multiplication. Fix $\mc F\in \mk B$ and $\mc H\in \mk S$, and suppose that $\mc F$ is $(\al,\mc P)$-homogeneous.  By definition
\begin{align*}
\iota(\mr{srk}_\mc P(\mc F\times \mc H))=& \iota(\min_{\xi<\theta} (\mr{srk}((\mc F\rest (\theta_\xi\times \{\xi\}))\times_\xi \mc H)))=\\
=&\max\{\min_{\xi<\theta}\iota(\mr{srk}((\mc F\rest (\theta_\xi \times\{\xi\})), \iota(\mr{srk}(\mc H) )\}  \\
=& \max\{\iota(\al),\iota(\mr{srk}(\mc H))\}=\iota(\mr{srk}_{\mc P}(\mc F)\cdot \mr{srk}(\mc H)).
\end{align*}
It is easy to see that $\times$ satisfies (M.2).
\fprue

We need to analyze then the lexicographical
orderings, finite or infinite. This is the content of the next part. Given $I$ and $J$, let $\pi_I:I\times
J\to I$, $\pi_J:I\times J\to J$ be the canonical projections. Given $i\in I$, $j\in J$, and  $\De\con I\times
J$, let
\begin{align*}
(\De)_i:=& \pi_J"(\pi_I^{-1}(i) \cap \De), \\
(\De)^j:= &\pi_I"(\pi_J^{-1}(j)\cap \De).
\end{align*}
be the corresponding sections.

Recall that given two partial orderings $\mbb P=(P,\le_P)$ and $\mbb Q=(Q,\le_Q)$, let $\mbb P\times_\mr{lex}
\mbb Q:=(P\times Q,<_\mr{lex})$ be the lexicographical product of $\mbb P$ and $\mbb Q$ defined by
$(p_0,q_0)<_\mr{lex} (p_1,q_1)$ if and only if $p_0<_P p_1$, or if $p_0=p_1$ and $q_0<_Qq_1$.  This
generalizes easily to finite products $\mc P_1 \times_\mr{lex}\mc P_2 \times_\mr{lex}\cdots \times_\mr{lex}
\mc P_n$, and the corresponding finite powers $\mc P_\mr{lex}^n$.  One can also define infinite
lexicographical products, but they are not going to be used here. Instead we use \emph{quasi-lexicographical
power} $\mc P_\mr{qlex}^{<\om}$  on $P^{<\om}$ defined by $(p_i)_{i<m} <_\mr{qlex} (q_i)_{i<n}$ if and only
if $m<n$ or if $m=n$ and $(p_i)_{i<m} <_\mr{lex} (q_i)_{i<m}$. Finally, let $\mr{lh}:P^{<\om}\to \om$ be the
length function. The main result here is the following.

\teor\label{ui32iuui32t7442212}
Let  $\mc P$ and $\mc Q$ be partial orderings.
\begin{enumerate}[(a)]
\item If $\mc P$ and $\mc Q$ have bases of families on the corresponding chains, then $\mc P \times_\mr{lex} \mc Q$ also has a basis of families on its chains.
\item   If  there is a basis on chains of $\mc P$, then there is
also a basis of families on chains of each finite lexicographical power $\mc P^n_\mr{lex}$ and there is a
basis of families on chains of $\mc P^{<\om}_\mr{lex}$.
\end{enumerate}
\fteor

\defi[Fubini product of families]
Given $I$ and $J$, let $\pi_I:I\times J \to I$, and $\pi_J:I\times J\to J$ be the canonical coordinate
projections.
 Given families $\mc F$  and $\mc G$ on $I$ and  $J$ respectively, let
 \begin{align*}
\mc F\circledast_\mr{F} \mc G:= & \conj{x\con I \times J}{\pi_I"x\in \mc F \text{ and $(x)_i\in \mc G$ for every $i\in I$} }.
 \end{align*}
 We call $\mc F \circledast_\mr{F} \mc G$ the Fubini product of $\mc F$ and $\mc G$.  Given $n\ge 1$, let $\mc F_\mr{F}^{n+1}=\mc F^n_\mr{F}\circledast_\mr{F}\mc F$.
\fdefi
It is easy to see that if $\mc F$  and $\mc G$ are families on chains of $\mc P$ and $\mc Q$, respectively,
then $\mc F\circledast_F \mc G$ is a family on chains of $\mc P\times_\mr{lex} \mc G$.

\defi[Power operation]
For each $n<\om$, let $\mc F_n$ be a family on chains of $\mc P^n_\mr{lex}$, and let $\mc G$ be a family on $\om$. We define $((\mc F_n)_n)^\mc G$ as the collection of all $x\con P^{<\om}$ such that
\begin{enumerate}[(i)]
\item  $x\cap [P]^n \in \mc F_n$ for every $n<\om$.
\item  $\mr{lh}"x\in \mc G$.
\end{enumerate}
Given a family  $\mc F$  on chains of $\mc P$, let $\mc F^\mc G:= ((\mc F_\mr{F}^n)_n)^\mc G$.

\fdefi

\lema \label{kloi8978785522222}
Let $\mc P$ be $\mc Q$  two partial orderings. For each $1\le n<\om$, let $\mc F_n$ be a family on chains of $\mc P_\mr{lex}^n$,  and suppose also that  $\mc G$ and $\mc
H$ are families on chains of $\mc Q$ and $\om$, respectively. Set $\mc F:=\mc F_1$.
\begin{enumerate}[(a)]
\item  If $\mc F_n$, $n<\om,$ $\mc G$, $\mc H$ are pre-compact, hereditary, then so are  $\mc F\circledast_\mr{F} \mc
G$, $((\mc F_n)_n)^\mc H$, and  $\mc F^\mc H$.
\item   If $\mc F_n$, $n<\om,$  $\mc G$ and $\mc H$ are countably ranked families, then
\begin{enumerate}[(b.1)]
\item $\mr{rk}(\mc F\circledast_\mr{F} \mc G) < \left( \mr{rk}(\mc G)\cdot \om)\right) \cdot (\mr{rk}(\mc F) +1)$,
\item $\mr{rk}(((\mc F_n)_n)^\mc H) < \sup_{n<\om} (\mr{rk}(\mc F_n) +1) \cdot (\mr{rk}(\mc H)+1)$,
\item  $\mr{rk}(\mc F^\mc H) <  (\mr{rk}(\mc F) \cdot \om)^\om \cdot (\mr{rk}(\mc H)+1)$, if $\mr{rk}(\mc F)\ge 1$.
\end{enumerate}
\item When the corresponding families are countable ranked,
\begin{enumerate}[(c.1)]
\item  $\mr{srk}_\mr{lex}(\mc F\circledast_\mr{F} \mc G)= \min\{\mr{srk}_\mc P(\mc F), \mr{srk}_\mc Q(\mc
G)\}$,
\item  $\mr{srk}_{\mr{qlex}}(((\mc F_n)_n)^\mc H)=\min\{\min_{n<\om} \mr{srk}_{\mc
P_\mr{lex}^n}(\mc F_n), \mr{srk}(\mc H)\}$,     and
\item  $\mr{srk}_{qlex}(\mc F^{\mc
H})=\min\{\mr{srk}_\mc P(\mc F), \mr{srk}(\mc H)\}$.
\end{enumerate}
\item If each $\mc F_n$ is $(\al,\mc P^n_\mr{lex})$-homogeneous,  $\mc G$ is $(\al,\mc Q)$-homogeneous and $\mc H$ is $\al$-homogeneous, then $\mc F\circledast_\mr{F} \mc G$ is $(\al, \mr{lex})$-homogeneous, and both $((\mc F_n)_n)^\mc H$ and $\mc F^\mc H$ is $(\al, \mr{qlex})$-homogeneous.
\end{enumerate}
\flema
\prue

\noindent \emph{The operation $\circledast_\mr{F}$}: The hereditariness is easy to prove. Suppose that $\mc
F$ and $\mc G$ are pre-compact. Let $(x_n)_n$ be  a sequence in $\mc F\circledast_\mr{F} \mc G$. Let $M\con
\om$ be infinite and  such that   $(\pi_I" x_n)_{n\in M}$ is a $\De$-sequence with root $y$.   Let $N\con M$
be infinite such that $((x_n)_p)_{n\in N}$ is a $\De$-sequence with root $z_p$ for every $p\in y$. Let
$x:=\bigcup_{p\in y} z_p$. It is easy to see  that $(x_n)_{n\in N}$ is a $\De$-sequence with root $x$.
Suppose now that $\mc F$ and $\mc G$ have countable rank,   and set  $\mc H:=\mc F\circledast_\mr{F} \mc G$.
We apply Proposition \ref{oi43jiro43ro34njrtguyy73}  to the projection $\pi_P: P\times Q\to P$ to conclude
that
$$
\mr{rk}(\mc H)< \sup_{y\in \mc F} (\mr{rk}(\conj{x\in \mc H}{\pi_P"(x)\con y})+1) \cdot (\mr{rk}(\mc F)+1).
$$
Now observe that for a given $y\in \mc F$ one has that
$$\conj{x\in \mc H}{\pi_P"x\con y}\con \bigsqcup_{p\in y} \conj{ \{p\}\times z\con P\times Q}{ z\in  \mc G },$$
by definition of $\mc H$.  Clearly $\conj{ \{p\}\times z\con P\times Q}{ z\in  \mc G }$ is homeomorphic to
$\mc G$, so by Proposition \ref{uygfuyhgu1} (iii.3.),
$$\mr{rk}(\conj{x\in \mc H}{\pi_P"x\con y})   <\mr{rk}(\mc G) \cdot \om.$$

\noindent \emph{The power operation}:  The fact that this operation preserves hereditariness is trivial to
prove. Suppose that each $\mc F_n$ is a pre-compact family on chains of $\mc P_\mr{lex}^n$ and that $\mc H$
is a pre-compact family on $\om$. Set $\mc Z:= ((\mc F_n)_n)^\mc H$, and suppose that $(x_k)_k$ is a sequence
in $\mc Z$.  We assume that $(\mr{lh}(x_n) )_n$ is a $\De$-sequence in $\om$ with root $z$. Now for each
$n\in z$ and each $k<\om$, let $y_k^{n}:= \mr{lh}^{-1}(n)\cap x_k \in \mc  F_n$.  Let $(x_k)_{k\in M}$ be a
subsequence of $(x_k)_{k}$ such that for each $n\in z$ one has that $(y_k^n)_{k\in M}$ is a $\De$-sequence
with root $y_n$.  It is easy to verify that $(x_k)_{k\in M}$ is a $\De$-sequence with root $\bigcup_{n\in z}
y_n$. The inequality in  (b.2) follows  from Proposition
\ref{oi43jiro43ro34njrtguyy73}.  The properties of $\mc F^\mc H$ follow from the corresponding properties of
the Fubini product and the power operation.

(c) follows from the fact that if $C=\{(p_n,q_n)\}_{n<\om}$ is a chain of $\mc P \times_\mr{lex} \mc Q$ then
there is an infinite $M\con \om$ such that either $p_m\neq p_m$ for every $m<n\in M$, or $p_m=p_n$ and
$q_m\neq q_n$ for every $m<n$ in $M$.  Similarly, given a chain $C=\{\bar p^n\}_{n<\om}$, there is an
infinite $M\con \om$ such that either $\mr{lh}(\bar p^n)= l$ for every $n\in M$ and $\{\bar{p}^n\}_{n\in M}$
is an infinite chain of $\mc F_\mr{F}^l$, or $\mr{lh}(\bar p^m)\neq    \mr{lh}(\bar p^n)$ for every $m\neq n$
in $M$.  (d) is a consequence of (a) (b) and (c).
\fprue

\section{Bases  of families  on trees}\label{bases_on_trees}
The goal of this section is to give a method to step up bases to trees. The idea is very natural: 
A tree $T$ is determined by its chains and antichains. Given two families $\mc A$
and $\mc C$ on antichains and on chains of a tree $T$, respectively, one can define a third  family $\mc
A\odot_T \mc C$ consisting of all subsets of $T$  generating a subtree whose antichains are in $\mc A$ and
its chains are in $\mc C$. In general, the antichains of a tree are difficult to understand; on the contrary,
the particular antichains consisting of immediate successors of a node are typically simpler (e.g. in a
complete binary tree), and it makes more sense to define $\mc A\odot_T \mc C$  in terms of these particular
simpler antichains. This operation on families will allow us, for example, to step up from a basis of
families on a cardinal number $\ka$ to a basis on its cardinal exponential $2^\ka$, and more.

Recall that a (set-theoretical) \emph{tree} $T=(T,<)$ is a set  of nodes $T$ with a partial order $<$ such
that $\conj{u<t}{u\in T}$ is well ordered for every $t\in T$. A \emph{rooted} tree is a tree with a minimal
element $0$, called the \emph{root} of $T$. All trees we use are rooted, so that whenever we say tree, we
mean a rooted tree.   Trees are a sort of lexicographical product of two orderings, the one defining the tree
order $<$ and the following.  Given $t\le u$ in $T$, let $\mr{Is}_t(u)$ be the immediate successor of $t$
which is below $u$, that is, $\mr{Is}_t(u)$ be the smallest $v \leq u$ such that $t<v$.
 Then, given $t\in T$ and $x\con T$, let
$$ \mr{Is}_t"x:=\conj{\mr{Is}_t(u)}{t<u\in x}.$$
For simplicity, we write $\textrm{Is}_t$ for $\mr{Is}_t"(T)$, that is, the set of all immediate successors of
$t$ in $T$. For every $t\in T$, fix a total ordering $<_t$ of $\mr{Is}_t$.
 Let $<_a$ be the   partial ordering in $T$ defined by $t<_a u$ if and only if there is $v$ such that $t,u\in \mr{Is}_v$ and $t<_v
 u$.   Hence, a chain with respect to $<_a$ is a set of immediate successors of a fixed node.   Notice that both $<$ and $<_a$ can be extended to a total ordering $\prec$ on $T$ by defining
 $t_0\prec t_1$ if and only if $t_0 < t_1$, or if $t_0$ and $t_1$ are $<$-incomparable and $\mr{Is}_{t_0\we t_1}(t_0)<_a \mr{Is}_{t_0\we
 t_1}(t_1)$. However, as we already observed, none of the notions involved in this work depends on the particular \emph{total} order one can choose on some index set, since being a chain with respect to a total order does not depend on the total order itself. 
 
 We are now able to state the main result of this section.

\teor \label{4jtierjteoirtertrtr344}
The following are equivalent for an infinite tree $T$.
\begin{enumerate}[(a)]
\item There is a basis of families  on $T$.
 \item  There is a basis of families on chains of $(T,<)$, if there is an infinite $<$-chain, and there is a basis of families of families on chains of $<_a$, if there is an infinite $<_a$-chain.
\end{enumerate}
\fteor
We pass now to recall well-known combinatorial principles on trees. Let $T=(T,<)$ be a
complete rooted tree with root $0$. Recall that a \emph{chain} of a tree is a totally ordered subset of it.     Let 
\begin{align*}
\mr{Ch}_a:=&\conj{s\in[T]^{<\om} }{s \text{ is a $<_a$-chain}}\\
\mr{Ch}_c:=&\conj{s\in [T]^{<\om} }{s \text{ is a $<$-chain}}.
\end{align*}
Notice that it follows from K\"onig's Lemma that when $T$ is infinite, either $\mr{Ch}_c$ is not compact, that is, there is an infinite $<$-chain, or $\mr{Ch}_a$ is not compact, that is, there is an
infinite $<_a$-chain. 
 Given $t \le u\in T$, let
 $$[t,u]:= \conj{v\in T}{ t\le v\le u}.$$
 Similarly, one defines the corresponding (semi) open intervals.  Given $x\con T$, let $x_{\le t}:=x\cap [0,t]$,  $x_{\ge t}:=\conj{u\in x}{t\le u}$, and let $x_{<t}$ and $x_{>t}$ be their open analogues.

Given $t,u\in T$, let
$$t\wedge u:= \max(T_{\le t}\cap T_{\le u}),$$
which is well defined by the completeness of $T$.  We say that $s\con T$ is \emph{$\wedge$-closed} when
$t\wedge u\in s$ for every $t,u\in s$. Given $s\con T$, let $\langle s\rangle$ be the \emph{subtree
generated} by $s$, that is, the minimal $\wedge$-closed subset of $T$ containing $s$. We say that a subset
$\tau\con T$ is a \emph{subtree} of $T$ when $\langle \tau\rangle = \tau$.

We introduce another operation that will be very useful. Given $t,u\in T$, let
\begin{align*}
 t\wedge_\mr{is} u:=& \left\{\begin{array}{ll}
 \min\{t,u\} & \text{if $t,u$ are comparable}\\
 \mr{Is}_{t\wedge u}(t) & \text{if $t,u$ are incomparable}.
 \end{array}
 \right.
\end{align*}
Given $s\con T$, let  $\langle s\rangle_\mr{is}$ be the minimal subset of $T$ that contains $s$ and that is closed under $\wedge$ and $\wedge_\mr{is}$. It is clear that $\langle s\rangle_\mr{is}$ is a substree and that $\langle s\rangle\con  \langle s\rangle_\mr{is}$. We   will say that $\tau\con T$ is an $\mr{is}$-subtree if $\langle \tau\rangle_\mr{is}$. The following is easy to prove. 
\defi 
Given a family $\mc F$ on $T$, let
\begin{align*}
\langle \mc F \rangle:= &\conj{x\con \langle s \rangle}{s\in \mc F},\\
\langle \mc F \rangle^\mr{is}:= &\conj{x\con \langle s \rangle_\mr{is}}{s\in \mc F},\\
\langle \mc F \rangle_T:=& \conj{\langle s \rangle}{s\in \mc F},\\
\langle \mc F \rangle_{\mr{is}-T}:=& \conj{\langle s \rangle_\mr{is}}{s\in \mc F},\\
\mr{Is}(\mc F):= &\conj{\mr{Is}_t" s }{t\in T \text{and }s\in \mc F}.
\end{align*}
\fdefi
We give some basic properties. The first one is easy to prove, so we leave the details to the reader. 

\prop
For every finite set $s \subseteq T$ and every $t \in T$, we have that
\begin{align*}
\langle s\rangle =& \conj{t_0\we t_1}{t_0,t_1\in s},\\
\langle s \rangle_\mr{is} = &\langle s\rangle \cup \{t_0 \we_\mr{is} t_1: t_0,t_1 \in s, t_0 \perp t_1\},\\
\mr{Is}_t" \langle s\rangle=& \mr{Is}_t" s.  
\end{align*}
		
In particular, $\langle s \rangle$ is finite whenever $s$ is finite, and $\mr{Is}(\mc F)=\mr{Is}(\langle \mc F\rangle)$ for every hereditary family $\mc F$ on $T$.   In general, if $(s_i)_{i \in I}$ is a
family of subsets of $T$, then
\begin{equation*}
\pushQED{\qed} 
 \langle \bigcup_{i\in I} s_i\rangle=\bigcup_{\{i,j\}\in [I]^2} \langle s_i\cup s_j\rangle.  \qedhere
\end{equation*}
\fprop

\prop	Let $\mc F$ be  a family on $T$.
\begin{enumerate}[(a)]
\item $\mc F=\langle \mc F\rangle$ if and only if $\mc F$ is hereditary and closed under $\langle \cdot \rangle$. 
\item $\langle \mc F\rangle$ is compact if and only if  $\langle \mc F\rangle_T$ is compact.
\item $\langle \mc F\rangle^\mr{is}$ is compact if and only if  $\langle \mc F\rangle_{\mr{is}-T}$ is compact.
\item If $\langle \mc F \rangle^\mr{is}$ is compact, then so is $\langle \mc F\rangle$. 
\end{enumerate}
 \fprop
\prue
$(a)$ is trivial. $(b)$:   Suppose that $\langle \mc F\rangle$ is compact. 
Then $\langle \mc F\rangle_T$ is compact, because $\conj{\tau\con T}{\tau\text{ is a subtree}}$ is compact.   Reciprocally, suppose that $\langle\mc F\rangle_T$ is compact, and let $(s_n)_n$ be a sequence in $\langle\mc F\rangle$. If $(\langle s_n\rangle)_{n\in M}$ is a converging subsequence with limit $\tau$, then it is easy to extract (using for example the Ramsey Theorem) a further converging subsequence $(s_n)_{n\in N}$ whose limit must belong to $\langle \mc F \rangle$ because this is always a closed set.  Similarly one proves  $(c)$.  $(d)$: Since $\langle \mc F\rangle \con \langle \mc F\rangle^\mr{is}$ and $\langle \mc F\rangle$ is closed, the compactness of $ \langle \mc F\rangle^\mr{is}$ implies the compactness of $ \langle \mc F\rangle$.  
\fprue

Given $s \con T$, let
$$(s)_{\mr{max}}: = \{\text{maximal elements of }s\}.$$

\defi
Let $s=(t_k)_{k\in \om}$ be a sequence of nodes in $T$.
\begin{enumerate}[(a)]
\item $s$ is called a \emph{comb} if $s$ is an antichain such that
\begin{equation*}
\text{$t_k\wedge t_l =t_k\wedge t_m$ and $t_k\wedge t_l<t_l\wedge t_m$ for every $k<l<m$. }
\end{equation*}
The chain $(u_k)_k$, $u_k=t_k\we t_l$ ($k<l$)  is called the $\we$-chain of the comb $(t_k)_k$.  
\begin{figure}[h]
\begin{tikzpicture}[scale=0.8,grow=up]
\node[circle, label=right:$u_1$, draw, fill=black, thin, scale=0.5] {}
child {
node[circle, label=right:$u_2$, draw, fill=black, thin, scale=0.5] {}
child {node[circle, label=right:$u_3$, draw, fill=black, thin, scale=0.5] {}
child {node{} edge from parent[dotted]}
child {node[circle, label=right:$t_3$, draw, fill=black, thin, scale=0.5] {}    }
}
child{node[circle, label=right:$t_2$, draw, fill=black, thin, scale=0.5] {}}
}
child {node[circle, label=right:$t_1$, draw, fill=black, thin, scale=0.5] {}}
 ;
\end{tikzpicture}
\caption{A comb $(t_k)_{k\in \om}$ and its corresponding $\we$-chain $(u_k)_k$.}
\end{figure}
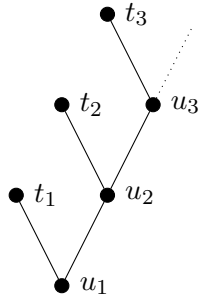

\item  $s$ is called a \emph{fan} if
\begin{equation*}
t_k\wedge t_l=t_{k'}\wedge t_{l'} \text{ for every $k\neq l$ and $k'\neq l'$}.
\end{equation*}
The node $u:=t_k\wedge t_l$ ($k\neq l$) is called the $\we$-root of the fan $(t_k)_k$.

\begin{figure}[h]

\begin{tikzpicture}[scale=0.8,grow=up]
\node[circle, label=below:$u$, draw, fill=black, thin, scale=0.5] {} child { node{$\cdots$}} child {
node[circle, label=above:$t_k$, draw, fill=black, thin, scale=0.5] {}     } child { node{$\cdots$} edge
from parent[draw=none]} child { node[circle, label=above:$t_2$, draw, fill=black, thin, scale=0.5] {}
} child { node[circle, label=above:$t_1$, draw, fill=black, thin, scale=0.5] {}       }
 ;
\end{tikzpicture}

\caption{A fan $(t_k)_{k\in \om}$ and its corresponding $\we$-root $u$.}
\end{figure}
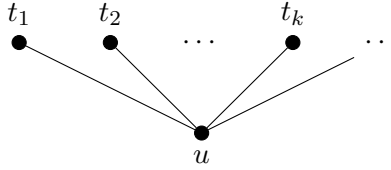
\end{enumerate}

\fdefi

\prop\label{4iojio4ojirjgf}
Every infinite subset of $T$ contains either an infinite chain, or an infinite comb, or an infinite fan.
\fprop
\prue
Fix a sequence $(t_k)_{k\in \om}$ such that $t_k\neq t_l$ for $k\neq l$.  By the Ramsey Theorem, there is $M_0$ such that either $(t_k)_{k\in M_0}$ is a chain or an antichain.

\clam If $(t_k)_{k\in M_0}$ is an antichain, then there is $M_1\con M_0$ such that $t_{k}\wedge t_{l}=t_{k}\wedge t_{m}$ for every $k<l<m$ in $M_1$.
\fclam
\prucl
If $(t_{k})_{k\in M_0}$ is an antichain and since there are no infinite decreasing chains in $T$, by the Ramsey Theorem, there is $M_1\con M_0$ such that
\begin{enumerate}[(a1)]
\item[(a1)] either $t_{k}\wedge t_{l}=t_{k}\wedge t_{m}$ for ever $k<l<m$ in $M_1$,
\item[(b1)] or else $t_{k}\wedge t_{l}<t_{k}\wedge t_{m}$ for ever $k<l<m$ in $M_1$.
\end{enumerate}
Let us see that (b1) cannot happen: Fix $k_0<k_1<k_2<k_3$ in $M_1$. Then
\begin{align*}
t_{k_0}\wedge t_{k_1}\wedge t_{k_i}=(t_{k_0}\wedge t_{k_1})\wedge(t_{k_0}\wedge t_{k_i})=t_{k_0}\wedge t_{k_1},
\end{align*}
for $i=2,3$. On the other hand,
\begin{equation*}
t_{k_0}\wedge t_{k_1}=t_{k_0}\wedge t_{k_1}\wedge t_{k_i}= (t_{k_0}\wedge t_{k_i}) \wedge(t_{k_1}\wedge t_{k_i})\in \{t_{k_0}\wedge t_{k_i}, t_{k_1}\wedge t_{k_i}\},
\end{equation*}
for $i=2,3$; so either $t_{k_0}\wedge t_{k_1}=t_{k_0}\wedge t_{k_i}$ for some $i=2,3$, or else
$t_{k_1}\wedge t_{k_2}=t_{k_0}\wedge t_{k_1}=t_{k_1}\wedge t_{k_3}$. Both cases are impossible since they contradict (b1).
\fprue

For each $k\in M_1$, let $u_k:=t_k\wedge t_l$ for some (all) $l>k$ in $M_1$. Yet again, since there are no infinite decreasing chains in $T$, by the Ramsey Theorem, there is $M_2\con M_1$ such that
\begin{enumerate}[(a2)]
\item[(a2)] either $u_k=u_l=\bar u$ for every $k<l$ in $M_2$,
\item[(b2)] or $u_k<u_l$ for every $k<l$ in $M_2$.
\end{enumerate}
If (a2) holds, then $(t_k)_{k\in M_2}$ is a fan with $\wedge$-root $\bar u$. If (b2) holds, then  $(t_k)_{k\in M_2}$ is a comb with $\wedge$-chain $(u_k)_{k\in M_2}$.
\fprue

We can now use the previous result to describe bettwer the small rank of a family on the tree $T$. 
\cor    \label{ieribgknbvvvv}
Every infinite subtree of $T$ contains either an infinite chain or an infinite fan. Consequently,
\begin{enumerate}[(a)]
\item If $A\con T$ is an infinite accumulation point of a sequence of subtrees of $T$, then $A$ is a  subtree of $T$ that contains either an infinite chain or an infinite fan.
\item  For every family $\mc F$ on $T$  with countable rank one has that
\[\pushQED{\qed}  \mr{srk}(\mc F)=\inf\conj{\mr{rk}(\mc F\rest X)}{X\text{ is an infinite chain, comb or fan}}. \qedhere \popQED \]
\end{enumerate}
\fcor
We can have now characterize the compactness of $\langle\mc F\rangle^{\mr{is}}$ in terms of $<$ and $<_a$-chains. 
\cor	\label{iuyug78687555}
  For a  family $\mc F$  the following are equivalent.
\begin{enumerate}[(a)]
\item $\langle\mc F\rangle^{\mr{is}}$ is compact.
\item $\langle\mc F\rangle^{\mr{is}}\cap \mr{Ch}_c$ and $\mr{Is}(\langle\mc F\rangle^\mr{is})$ are compact. 
\end{enumerate}

\fcor
\prue	 
Suppose that $\langle \mc F\rangle^\mr{is}$ is compact.  The family $\langle \mc F\rangle^\mr{is}$ is by definition hereditary, so   $\langle \mc F\rangle^\mr{is} \cap \mr{Ch}_c$ is closed, and since $\langle \mc F\rangle^\mr{is} \cap \mr{Ch}_c\con \langle \mc F\rangle^\mr{is}$, it is  compact.  Notice that $\mr{Is}(\langle \mc F\rangle^{\mr{is}})$  is hereditary, because $\langle\mc F\rangle^\mr{is}$ is so, by definition. Let $(s_n)_n$ be a sequence in $\mr{Is}(\langle\mc F\rangle^\mr{is})$.  If there is an infinite  subsequence $(s_n)_{n\in M}$ such that $\#s_n\le 1$, then clearly there is a further convergent subsequence of $(s_n)_{n\in M}$  with limit ($\buit$ or  a singleton) in $\mr{Is}(\langle\mc F\rangle^\mr{is})$ because it  is hereditary. Otherwise, from some point on $\#s_n\ge 2$. This means that for such $s_n$ one has that $s_n=\conj{t\wedge_\mr{is} u}{t\neq u\in x_n}$ for some $x_n \in \mc F$, and so    $s_n\in \langle \mc F\rangle^\mr{is}$. Hence, there is a converging subsequence whose limit belongs to $\mr{Is}(\langle\mc F\rangle^\mr{is})$, because $\mr{Is}(\langle\mc F\rangle^\mr{is})$ is hereditary.  Suppose now that {\it(b)} holds, and let $(s_n)_n$ be a sequence in $\langle \mc F\rangle^\mr{is}$. Going towards a contradiction, suppose that $(\tau_n)_n$ is a sequence of $\mr{is}$-subtrees of $T$ converging to an infinite $\mr{is}$-subtree $A$ of $T$. 
From Corollary \ref{ieribgknbvvvv} we know that one of the following holds.   $A$ contains an infinite chain $C$. Then $(\tau_n\cap C)_n$ is a sequence    in  $\langle\mc F\rangle^{\mr{is}}\cap \mr{Ch}_c$ with infinite limit $C$, impossible. The second option is that $A$ contains an infinite fan $X$ with root $t$.  Then the sequence $(\mr{Is}_t''(s_n))_n$ converges to the infinite set $\mr{Is}_t''X$, contradicting the fact that $\mr{Is}(\langle\mc F\rangle^\mr{is})$ is compact. 
\fprue

So, we have the following characterization of the compactness of $\langle\mc F\rangle^{\mr{is}}$ in terms of compact properties of $\langle\mc F\rangle$. 
\cor \label{kjniuhrieuiuhrerer22}
$\langle \mc F\rangle^\mr{is}$ is compact if and only if $\langle \mc F\rangle$ and $\mr{Is}(\langle \mc F\rangle)$ are compact. 
 \fcor
\prue
We know that the compactness of $\langle \mc F\rangle^\mr{is}$ implies the compactess of both $\langle \mc F\rangle$ and $\mr{Is}(\langle \mc F\rangle)=\mr{Is}(\langle \mc F\rangle^\mr{is})$.  Conversely, suppose  that $\langle \mc F\rangle$ and $\mr{Is}(\langle \mc F\rangle)$ are compact. We use the characterization of compactness of $\langle\mc F \rangle^\mr{is}$			in Corollary \ref{iuyug78687555}, so it rests to show that  $\langle \mc F\rangle^\mr{is} \cap \mr{Ch}_c$ is compact: Set $\mc C:=\langle \mc F\rangle \cap \mr{Ch}_c$, $\mc D:=\langle \mc F\rangle^\mr{is}\cap \mr{Ch}_c$. Both are hereditary families, and $\mc C$ is  compact, by hypothesis. 
Let $\mr{Is}$ be the collection of all immediate successors of $T$, and let $\pi: \mr{Is}\to T$ be the mapping assigning to each $t\in \mr{Is}$ its immediate predecessor $\pi(t)$. Then $\pi$ is an adequate mapping, so $\pi^{-1}(\mc C)$ is compact and hereditary.  Since $\mc D\con \mc C\sqcup \pi^{-1}(\mc C)$, $\mc D$ is compact, as desired. 
 \fprue
We analyze now the $<$ and $<_a$-chains of the tree generated by the union of two subtrees of $T$. 
\prop\label{loi21950e}
Suppose that $\tau_0,\tau_1$ are subtrees of $T$, and $t\in T$. Then
\begin{align}
\mr{Is}_t (\langle \tau_0\cup \tau_1\rangle)=&  \mr{Is}_t \tau_0 \cup  \mr{Is}_t \tau_1 ,\label{dfjsfijse32d}\\
\mr{Ch}_c \rest \langle \tau_0\cup \tau_1\rangle \con & (\mr{Ch}_c\rest\tau_0 ) \sqcup_c(  \mr{Ch}_c \rest\tau_1)\sqcup_c [T]^{\le 1}. \label{dfjsfijse32d1}
\end{align}
\fprop
\prue
\eqref{dfjsfijse32d} follows from the fact that $\mr{Is}_t"(\langle s\rangle)=\mr{Is}_t"(s)$.   For \eqref{dfjsfijse32d1}, let $c$ be a chain of $\langle \tau_0\cup \tau_1\rangle$, and suppose that $t_0<t_1$ belong to $c \setminus (\tau_0\cup \tau_1)$. Then $t_0= u_0\wedge v_0$ and $t_1= u_1\wedge v_1$ with $u_0,u_1\in \tau_0$ and $v_0,v_1\in \tau_1$. Then, either $u_0 \wedge u_1= t_0$ or $v_0\wedge v_1=t_0$, and both are impossible.
\fprue

We introduce the   operation $\odot_T$, which is the crucial tool to build bases on a tree $T$ from bases on $<$-chains and bases on $<_a$-chains.  

\defi[The operation $\odot_T$]
  Given  two families $\mc A$ and $\mc C$  on $T$, let $\mc A \odot_T \mc C$ be the family of
all finite subsets $s$ of $T$ such that
\begin{enumerate}[(a)]
\item 
For every $t\in  T $ one has that $\mr{Is}_t"\langle s\rangle \in \mc A$.
\item   
Every chain in $\langle s\rangle$ belongs to $\mc C$.
\end{enumerate}
\fdefi

\nota\label{oiuh324ry32334}
\begin{enumerate}[(i)]
\item The family $\mc A \odot_T \mc C$ is closed under generated subtrees, that is, $\langle s\rangle\in \mc A\odot_T \mc C$ if $s\in \mc A\odot_T \mc C$.
\item   When  $[T]^{\le 1}_a\con \mc A$, the condition {\it(a)} above is equivalent to
\begin{enumerate}[(a')]
\item For every  $t\in \langle s\rangle$, one has that $ \mr{Is}_t"(s)\in \mc A$.
\end{enumerate}
\item When  $[\mr{Is}_t]_a^{\le 1}\con \mc A$ for all $t\in T$, $(\mc A \odot_T \mc C)\cap \mr{Ch}_c = \mc C \cap \mr{Ch}_c$, and when $[T]^{\le 2}_{c}\con \mc C$, then $\mr{Is}(\mc A \odot_T \mc C) = \mc A \cap \mr{Ch}_a$.
\end{enumerate}

\fnota

We introduce some notation. We are going to use $a$ and $c$ to refer to $<_a$ and $<$. For example,
$$\sqcup_a, \, \mr{srk}_a,\, \sqcup_c,\, \mr{srk}_c$$
denote $\sqcup_{(T,<_a)}$, $\mr{srk}_{(T,<_a)}$,  $\sqcup_{(T,<)}$ and $\mr{srk}_{(T,<)}$ respectively. We now compare small ranks of $\mc A\odot_T$ and of $\mc A$ and $\mc C$. 
\prop 
\label{kerhuuerhgf}  Let $\mc A$ and $\mc C$ be two families on chains of $(T,<_a)$ and  of $(T,<)$, respectively. \begin{enumerate}[(a)]
\item  If  $\mc A$ and $\mc C$  are compact, hereditary, then so is $\mc A \odot_T
\mc C$.
\item $\mr{srk}(\mc A\odot_T \mc C)\le
\mr{srk}_a(\mc A)$ if $\mr{Ch}_a$  is non compact, and  $\mr{srk}(\mc A\odot_T \mc C)\le \mr{srk}_c(\mc C) $ if $\mr{Ch}_c$  is non compact.

\item $\mr{srk}(\mc A\odot_T \mc C  )=\mr{srk}_a(\mc A)$ if  $\mr{Ch}_c$ is compact and  $[T]_c^{\le 2}\con \mc C$. 
\item $\mr{srk}_c(\mc C) \le \mr{srk}(\mc A\odot_T (\mc C\sqcup_c[T]^{\le 1} ))\le\mr{srk}_a(\mc C)+1$  if $\mr{Ch}_a$ is compact and     $[T]_a^{\le 2}\con \mc A$. 

\item $\min\{\mr{srk}_a(\mc A),\mr{srk}_c(\mc C)\} \le \mr{srk}(\mc A\odot_T (\mc C\sqcup_c[T]^{\le 1} ))\le \min\{\mr{srk}_a(\mc A),\mr{srk}_c(\mc C)+1\} $  if $\mr{Ch}_a$ and $\mr{Ch}_c$ are not compact, and $[T]_c^{\le 2}\con \mc C$,  $[T]_a^{\le 2}\con \mc A$. 

\end{enumerate}
\fprop
\prue
Set $\mc F:=\mc A \odot_T \mc C$.  {\it(a)}: Hereditariness of  $\mc F$ is trivial. Suppose that $s$ is an
infinite subset of $T$ which is the limit of a sequence $(\tau_k)_k$ in $\mc F$. Since $\mc F$ is, by
definition, $\we$-closed, we may assume without loss of generality that each $\tau_k$ is a subtree.  It
follows that $s$ is a subtree of $T$ as well. Hence, by Corollary \ref{ieribgknbvvvv}, $s$ contains either an
infinite chain $C$, or an infinite fan $F$. In the first case, $\tau_k\cap C\in \mc C$ for every $k$ and
$\tau_k\cap C\to_k C$, which is impossible since $\mc C$ is compact. If $s$ contains an infinite fan $F$ with
$\we$-root $u$, then $\mr{Is}_u"(F)$ is an accumulation point of the sequence $(\mr{Is}_u"(s_k))_k$ in $ \mc
A$, which is impossible by the compactness of $\mc A$.

{\it(b)} follows from the fact that if $X$ is a $<_a$-chain, then $(\mc A\odot_T\mc C)\rest X\con \mc A\rest X$, and if $X$ is a  $<$-chain, then $(\mc A\odot_T\mc C)\rest X\con \mc C\rest X$.  {\it(c)}, {\it(d)} and {\it(e)} follow, by the means of Corollary
\ref{ieribgknbvvvv}, from the following claim.
\clam
\begin{enumerate}[(i)]
\item If $[T]_c^{\le 2}\con \mc C$,  then 
$\mr{rk}((\mc A\odot_T\mc C)\rest X) \ge \mr{srk}_a(\mc A)$ for every infinite fan $X$.
\item If $[T]_a^{\le 1}\con \mc A$, then $\mr{rk}((\mc A\odot_T\mc C)\rest X) \ge \mr{srk}_c(\mc C)$ for every infinite chain $X$.  
\item If $[T]_a^{\le 2}\con \mc A$, then $\mr{rk}((\mc A\odot_T(\mc C\sqcup_c [T]^{\le 1}))\rest X \ge \mr{srk}_c(\mc C)$ for every infinite comb $X$.  

\end{enumerate}

\fclam
\prucl
{\it(i)}: Set $\mc F:= \mc A\odot_T \mc C$.  Suppose that  $X=\{t_n\}_{n<\om}$   is an infinite fan with $\wedge$-root
$u$.  For each $n<\om$, let $v_n:= \mr{Is}_u(t_n)$. Since for every $x\con \om$ one has that $\langle
\{t_n\}_{n\in x}\rangle= \{t_n\}_{n\in x}\cup \{u\}$, it follows that the maximal chains of  $\langle
\{t_n\}_{n\in x}\rangle$ have cardinality 2, so, they belong to $\mc C$, by hypothesis. This means that
$\{t_n\}_{n\in x}\in \mc F$ if and only if $\{v_n\}_{n\in x}\in \mc A$, and consequently $\rk{\mc F\rest X}=
\rk{\mc A \rest \{v_n\}_n}\ge \srk{\mc A}{a}$. {\it (ii)}: Set again $\mc F:= \mc A\odot_T \mc C$.   Suppose that  $X$ is an infinite chain.  Since $[T]_a^{\le 1}\con \mc A$, it follows that $\mc F\rest X=\mc C\rest X$, hence $\mr{rk}(\mc F\rest X)\ge \mr{srk}_a(\mc C)$. {\it(iii)}:  Set now   $\mc F:= \mc A\odot_T( \mc C\sqcup_c [T]^{\le 1}$. 
Suppose that $X=\{t_k\}_k$ is a comb with
$\wedge$-chain  $C:=\{u_k\}_k$. 
Given $\{u_k\}_{k\in x}\in \mc C$ we know that $\{u_k\}_{k\in x}\cup \{t_{p}\}\in
C\sqcup_c [T]^{\le 1}$, where $p:=\max x$. Since
$$\langle \{t_k\}_{k\in x}\rangle = \{t_k,u_k\}_{k\in x} $$
it follows that $\langle \{t_k\}_{k\in x}\rangle\in \mc A \odot_T (\mc C\sqcup_c [T]^{\le 1})$ whenever
$\{u_k\}_{k\in x}\in \mc C$. So, $\{u_k\}_{k\in x}\in \mc C \mapsto \{t_k\}_{k\in x}\in \mc F$ is continuous
and 1-1. Hence $\mr{rk}(\mc F\rest X)\ge \mr{srk}_c(\mc C)$.
 \fprucl
 \fprue

We are ready to define a basis of $T$ from a basis on $<$-chains and a basis on $<_a$-chains. In fact we will define a pseudobasis on $T$ which can be modified to give a basis (Proposition  \ref{iweiijoweee44}). 

\defi[The pseudo-basis on $T$] \label{8747878yyy7w} Let $T$ be an infinite tree, and suppose that for
$(T,<_a)$ and $(T,<)$ either they have a basis of families on their chains or they do not have infinite
chains. Let $(\mk B^a, \times_a)$ be either a basis of families on chains of $(T,<_a)$, or $\mk B^a:=\conj{\mc
A}{\mc A \text{ is compact hereditary of finite rank and } \mc A\con \mr{Ch}_a}$ and $\mc A \times_a \mc H:= \mc A$ for every $\mc
A\in \mk B^a$ and $\mc H\in \mk S$, if there are no infinite $<_a$-chains. We define similarly $(\mk
B^c,\times_c)$. Let $\mk B$ be  the collection of all families $\mc F$ on $T$ such that
\begin{enumerate}[(BT.1)]
\item  $\mc F$ and $\langle \mc F\rangle$   are homogeneous, $\mr{Is}(\mc F)\con \mc F$ and $ \rk{\langle \mc F \rangle}<\iota(\srk{\mc F}{})$.  In addition, if the rank of $\mc F$ is finite, then there is some $n<\om$ such that $\mc F\con [T]^{\le n}$.   
\item  $\mr{Is} \langle\mc F \rangle   \in \mk B^a$ and $\langle \mc F \rangle \cap \mr{Ch}_c \in \mk B^c$.
\end{enumerate}
Given $\mc F\in  \mk B$ and $\mc H\in \mk S$ of infinite rank, let
\begin{enumerate}[(BT.1)] \addtocounter{enumi}{2}
\item  $  \mc F  \times \mc H :=
((\mc A \times_a \mc H )\sqcup_a [T]^{\le 1}) \odot_T ((\mc C\times_c \mc H)\boxtimes_c 5) $  where $\mc
A:= \mr{Is}(\langle \mc  F \rangle)$ and $\mc C:=\langle \mc F \rangle\cap \mr{Ch}_c$.
\end{enumerate}

\fdefi
We are going to see that $\mk B$ is a pseudobasis.  We start by analyzing the small ranks of families of  $\mk B$. 
\prop \label{jkluiui44384} 
Suppose that $\mc F\in \mk B$. Then
\begin{enumerate}[(a)]

\item  $\iota(\mr{srk}(\mc F))=\iota(\mr{srk}(\langle  \mc F \rangle))$.
\item 
$\iota(\mr{srk}(\mc F))=\iota(\mr{srk}(\langle \mc F \rangle))=\iota(\mr{srk}_c( \mc F  \cap
\mr{Ch}_c))=\iota(\mr{srk}_c( \langle \mc F  \rangle \cap \mr{Ch}_c))   $ if   $\mr{Ch}_c$  is not
compact.

\item 
$\iota(\mr{srk}(\mc F))=\iota(\mr{srk}(\langle \mc F \rangle))=\iota(\mr{srk}_a( \mr{Is}(\mc F)))=\iota(\mr{srk}_a( \mr{Is}(\langle \mc F  \rangle)))   $ if   $\mr{Ch}_a$  is not
compact.
\end{enumerate}
\fprop
\prue Fix $\mc F\in \mk B$.   {\it (a)}: 
$$\mr{srk}(\mc F)\le \mr{srk}(\langle \mc F \rangle)\le \mr{rk}(\langle \mc F \rangle)<\iota(\mr{srk}(\mc F)),$$
by property {\it(B.T.1)}. This means that $\iota(\mr{srk}(\mc F))=\iota(\mr{srk}(\langle  \mc F \rangle))$.
{\it (b)}:  Suppose that $\mr{Ch}_c$ is not compact.  $\mc F$ is hereditary, so $\mc F\con \langle \mc F\rangle$. Also,   there are infinite $<$-chains  of $T$ and if $C$ is  a such chain, then $(\mc F\cap \mr{Ch}_c)\rest C= \mc F\rest C$. Then,
$$\mr{srk}(\mc F)\le \mr{srk}_c(\mc F\cap \mr{Ch}_c)\le\mr{srk}_c(\langle\mc F\rangle\cap \mr{Ch}_c)\le \mr{rk}(\langle \mc F \rangle)<\iota(\mr{srk}(\mc F)).$$
This implies that $\iota(\mr{rk}(\mc F))= \iota(\mr{srk}_c(\mc F\cap \mr{Ch}_c))=\iota(\mr{srk}_c(\langle \mc F\rangle\cap \mr{Ch}_c) )$. 

{\it (c)}: Suppose that $\mr{Ch}_a$ is not compact. Similarly than for {\it(b)}, now using that by hypothesis
$\mr{Is}(\langle \mc F\rangle)=\mr{Is}(\mc F)\con \mc F\con \langle \mc F\rangle$, we obtain that 
$$\mr{srk}(\mc F)\le \mr{srk}_a(\mr{Is}(\mc F))=\mr{srk}_a(\mr{Is}(\langle\mc F\rangle))\le \mr{rk}(\langle \mc F \rangle)<\iota(\mr{srk}(\mc F)).$$

 $$\mr{srk}(\mc F)\le \mr{srk}_c(\mc F\cap \mr{Ch}_c)\le \mr{rk}(\mc F)<\iota(\mr{srk}(\mc F)),$$
hence $\iota(\mr{srk}(\mc F))=\iota(\mr{srk}_c(\mc F\cap \mr{Ch}_c))$. 
So, $\iota(\mr{rk}(\mc F))= \iota(\mr{srk}_a(\mr{Is}(\mc F)))=\iota(\mr{srk}_a(\mr{Is}(\langle \mc F\rangle) ))$. 
\fprue

The next two results are the keys to show Theorem \ref{4jtierjteoirtertrtr344}.
\lema\label{oi34r49778576} Suppose that $\mc F$ is a compact family  such that $\mc F=\langle \mc F\rangle$.  Then,
$$\mr{rk}(\mc F) < \max\{ \iota(\mr{rk}(\mr{Is}(\mc F))\cdot \om)),\iota(\mr{rk}(\mc F\cap \mr{Ch}_c)\cdot \om))\}.$$
 
\flema
\lema \label{kjjknkj744}
$\times$ is a multiplication.
\flema
Their proofs are involved, so we reserve specific subsections to present them. In fact, 
 Lemma  \ref{oi34r49778576}  is a consequence of an appropriate   upper bound of the rank of  $\mc A\odot_T \mc C$ in terms of the ranks of $\mc A$ and $\mc C$, and it   is
done  in the next Subsection \ref{iu4787877555}. The  proof of Lemma \ref{kjjknkj744} is mostly combinatorial; It is given in the
Subsection \ref{iu48767457655}, where we find the canonical form of a $\De$-sequence of finite subtrees of
$T$.  The next  gives, Using Lemma \ref{oi34r49778576}, a criteria on $\mc A$ and $\mc B$ to guarantee that $\mc A\odot_T\mc C\in \mk B$.
\prop\label{ioioio32uihihu66} 
 Let $\mc A\in \mk B^a$ and $\mc C\in \mk B^c$ be  such that
 \begin{enumerate}[(a)]
 \item    $[T]_a^{\le 2}\con \mc A$ and
$[T]_c^{\le 2}\con \mc C$.
\item If $\mr{Ch}_c$ is non compact, then  there is some $\mc D\in \mk B^c$ such that $\iota(\mr{srk}(\mc D))= \iota(\mr{srk}(\mc C))$ and $\mc D\sqcup_c [T]^{\le 1}\con \mc C$.   
 \item  $\max\{\mr{srk}_a(\mc A),\mr{srk}_c(\mc C)\}\ge \om$. 
\item   $\iota(\mr{srk}_a(\mc A))=\iota(\mr{srk}_c(\mc C))$ if
$\mr{Ch}_a$ and $\mr{Ch}_c$ are both non compact.  
 \end{enumerate} 
 Then $\mc A\odot_T \mc C  \in \mk B$.

\fprop
\prue
 Fix $\mc A$ and $\mc C$ as in the statement, and set $\mc F:= \mc A\odot_T \mc C$. Then $\langle \mc F \rangle=\mc F$, by definition of $\odot_T$,  $\mr{Is}(\mc F)=\mc A$, because we are assuming $[T]_c^{\le 2}\con \mc C$, and $\mc F\cap \mr{Ch}_c= \mc C$, because $[T]_a^{\le 1}\con \mc A$. This implies that {\it (BT.2)} holds for $\mc F$. Suppose first that $\mr{Ch}_c$ is compact. Then,  by the choice of $\mk B^c$, we know that $\mr{rk}(\mc C)<\om$.    Also,  $\srk{\mc F }{ }=\srk{\mc A}{a}$. Hence, $\la:= \iota(\mr{srk}_a(\mc A))=\iota(\mr{srk}(\mc F))$.  Since $\mc A$ is homogeneous, $\mr{rk}(\mc A)<\la$.  We know also that $\la>\om$, hence it follows that   Lemma   \ref{oi34r49778576} that $\mr{rk}(\mc F) <\la=\mr{srk}(\mc F)$, so {\it(BT.1)} holds for $\mc F$.   Suppose that $\mr{Ch}_a$ is compact.  Then it follows from  condition (b) and Proposition  \ref{kerhuuerhgf} {\it(d)} that $\mr{srk}(\mc D)\le \mr{srk}(\mc A\odot_T (\mc D\sqcup_c [T]^{\le 1}))\le \mr{srk}(\mc F)$. Hence, by hypothesis, $\la=\iota(\mr{srk}(\mc C))\le \iota(\mr{srk}(\mc F))$. Since $\mr{srk}_a(\mc A)=0$, it follows that   $\la>\om$. Since $\mc A$ has finite rank, and $\mc C$ is homogeneous,  it follows from Lemma  \ref{oi34r49778576}  that $\mr{rk}(\mc F)<\la=\iota(\mr{srk}(\mc F))$.   Suppose that $\mr{Ch}_a$ and $\mr{Ch}_c$ are non compact. Set $\la= \iota(\mr{srk}_a(\mc A))= \iota(\mr{srk}_c(\mc C))$.  Then, as before, $\iota(\mr{srk}(\mc F))=\la>\om$. Since both $\mc A$ and $\mc C$ are homogeneous, it follows from Lemma \ref{oi34r49778576}  that 
 $\mr{rk}(\mc F)<\la=\iota(\mr{srk}(\mc F))$. 
\fprue

We are now ready to prove the main result of this section.

\prue[\textbf{Proof of Theorem \ref{4jtierjteoirtertrtr344}}]   

Let us see that $(\mk B, \times)$ defined on
Definition \ref{8747878yyy7w} is a pseudo-basis of families on $T$ (see Proposition \ref{iweiijoweee44}).  {\it(B.1')}: We first see that the cubes belong to $\mk B$:    Notice that if $\tau$ is a finite tree,
then
\begin{equation}\label{uuueuiewe}
\#\tau \le  \frac{a(\tau)^{c(\tau)+1}-1}{a(\tau)-1}
\end{equation}
where $a(\tau)$  and $c(\tau)$ are the maximal cardinality of a $<_a$-chain and a $<$-chain, respectively.
Set  $\mc F=[T]^{\le n}$, $n\ge 1$. Then $\langle \mc F \rangle\con [T]^{\le n^{n+1}}$. Hence, $\mr{Is}(\mc F)=[T]_a^{\le n}\con [T]^{\le n}=\mc F$,  $\langle \mc F
\rangle$ is homogeneous and $\iota(\mr{srk}(\mc F))= \iota(\mr{srk}(\langle \mc F \rangle))=\om$, so $\mc F$ satisfies {\it(BT.1)}.   
Since $\mr{Is}(\mc F)=[T]_a^{\le n}$, $\mr{Is}(\mc F)\in \mk B^a$ whether $\mr{Ch}_a$ is compact or not.  Since $\langle \mc F \rangle \con [T]^{n^{n+1}}$, $\mc G:=\langle \mc F \rangle \cap \mr{Ch}_c\con [T]_c^{\le n^{n+1}}$ is a compact, hereditary family containing $[T]^{\le 1}$, so it belongs to $\mk B^c$ when $\mr{Ch}_c$ is compact.   When $\mr{Ch}_c$ is not compact,   $\iota(\mr{srk}_c(\mc G))=\om=\iota(\mr{srk}_c([T]_c^{\le n^{n+1}}))$ and $[T]_c^{\le n^{n+1}}\in \mk B^c$. Since $\mk B^c$ is a basis, it follows that $\mc G\in \mk B^c$.   

We fix now $\om \le \al<\ou$, and we prove that there is $\mc F\in \mk B$ such that $\al\le \mr{srk}(\mc F)\le \iota(\al)$:     Suppose first that $\mr{Ch}_a$ is compact.   Let $\mc D\in \mk B_\al^c$, and set  $\mc C:=\mc D\sqcup_c[T]^{\le1}$. It follows from Proposition \ref{ioioio32uihihu66} that $\mc F:=[T]_a^{\le 2} \odot_T\mc C\in \mk B$. We check that $\al \le \mr{srk}(\mc F) \le \iota(\al)$.     Proposition   \ref{kerhuuerhgf} {\it(d)}  gives that $\al=\mr{srk}_c(\mc D)\le \mr{srk}(\mc F)\le \mr{srk}_c(\mc D)+1= \al +1<\iota(\al)$, as desired.    Suppose now that $\mr{Ch}_c$ is compact. Let $\mc A\in \mk B^a_\al$ and $\mc C:=[T]_c^{\le 2}$. Then, Proposition \ref{ioioio32uihihu66} gives that  $\mc F:=\mc A\odot_T \mc C \in \mk B$, and Proposition   \ref{kerhuuerhgf} {\it(c)} gives that $\al=\mr{srk}_a(\mc A)=\mr{srk}(\mc F)$.    Finally, if $\mr{Ch}_a$ and $\mr{Ch}_c$ are non compact, we choose $\mc A\in \mk B_\al^a$ and $\mc C\in \mk B_\al^c$. Then, Proposition \ref{ioioio32uihihu66} and Proposition  \ref{kerhuuerhgf} {\it(e)}  give that $\mc F:=(\mc A \sqcup_a [T]_a^{\le 1})  \odot_T (\mc C \sqcup_c [T]^{\le 1} \in \mk B$ and  $\al \le \min\{\mr{srk}_a(\mc A\sqcup_a [T]^{\le 1}), \mr{srk}_c(\mc C)\}\le \mr{srk}(\mc F)\le \min\{\mr{srk}_a(\mc A \sqcup_a [T]^{\le 1}), \mr{srk}_c(\mc C)+1\} = \al+1<\iota(\al)$.  This ends the proof of {\it(B.1')}.   

\noindent {\it(B.2')}: We have to check that $\mk B$ is closed under $\cup$ and $\sqcup$.   Fix $\mc F,\mc G\in \mk B$.  Set $\mc B:= \mc F\cup \mc G$. Then Proposition \ref{uygfuyhgu1} gives that $\mc B$ and $\langle \mc B\rangle=\langle \mc F\rangle \cup \langle \mc G\rangle$   are homogeneous and 
$$\mr{rk}(\langle \mc B\rangle)=\max\{\mr{rk}(\langle \mc F\rangle), \mr{rk}(\langle \mc G\rangle)\}< \max\{\iota(\mr{srk}(  \mc F)), \iota(\mr{srk}( \mc G))\}\le \iota(\mr{srk}(\mc B)).$$ 
Also, $\mr{Is}(\mc B)=\mr{Is}(\mc F)\cup \mr{Is}(\mc G)\con \mc B$, so {\it(BT.1)} holds for $\mc B$.  Since $\mr{Is}(\mc B)=\mr{Is}(\mc F)\cup \mr{Is}(\mc G)$, it follows that $\mr{Is}(\mc B)\in \mk B^a$ whether $\mr{Ch}_a$ is compact or not.  Similarly, $\langle\mc B\rangle\cap \mr{Ch}_c=(\langle \mc F\rangle\cap \mr{Ch}_c)\cup (\langle\mc G\rangle \cap \mr{Ch}_c)\in \mk B^c$. 
 Let us see now 
that $\mc B:=\mc F\sqcup \mc G\in \mk B$. We suppose that $\mc F,\mc G\neq \{\buit\}$, so $[T]^{\le 1}\con \mc F,\mc G$ (as they are homogeneous families).  We know that $\mc B$ is homogeneous and $\mr{rk}(\mc B)=\mr{rk}(\mc F)\dot{+}\mr{rk}(\mc G)$. Suppose first that the ranks of $\mc F$ and $\mc G$ are finite. Then by {\it(BT.1)}, there is some $n$ such that $\mc F, \mc G\con [T]^{\le n}$, hence $\mc B\con [T]^{\le 2n}$.  It follows from \eqref{uuueuiewe} that there is some $m$ such that $\langle \mc B\rangle \con [T]^{\le m}$. This implies  that $\mr{rk}(\langle \mc B\rangle)<\om=\iota(\mr{srk}(\mc B))$, so {\it(BT.1)} holds for $\mc B$.  Now, $\mr{Is}(\langle \mc B\rangle)=\mr{Is}(\mc B)=\mr{Is}(\mc  F) \sqcup_a \mr{Is}(\mc G)\in \mk B^a$ whether $\mr{Ch}_a$ is compact or not. Similarly, $\mc B\cap \mr{Ch}_c\in \mk B^c$.  Suppose now that  $\mc F$ or $\mc G$ has infinite rank.    We have, by Proposition  \ref{loi21950e},  the following inclusions: 
$$\langle \mc F\sqcup \mc G\rangle \con\langle \langle\mc F\rangle \sqcup \langle\mc G\rangle \rangle  \con (\langle \mc F\rangle \sqcup_a \langle \mc G\rangle )\odot_T  (\langle \mc F\rangle \sqcup_c \langle \mc G\rangle \sqcup_c [T]^{\le 1} ).$$
Now,
\begin{align*}
\mr{rk}(\langle F\rangle  \sqcup_a \langle \mc G\rangle) \le&  \mr{rk}(\langle F \rangle)\dot+ \mr{rk}(\langle G \rangle) < \max\{\iota(\mr{srk}(\mc F)),\iota(\mr{srk}(\mc G))\}\le \iota(\mr{srk}(\mc B))\\
\mr{rk}(\langle F\rangle  \sqcup_c \langle \mc G\rangle \sqcup_c [T]^{\le 1}) \le&  \mr{rk}(\langle F \rangle)\dot+ \mr{rk}(\langle G \rangle)+1 <\max\{\iota(\mr{srk}(\mc F)),\iota(\mr{srk}(\mc G))\}\le \iota(\mr{srk}(\mc B)).
\end{align*}
Since $\mc B$ has infinite rank and homogeneous, $\iota(\mr{srk}(\mc B))>\om$, so  it follows from Lemma  \ref{oi34r49778576}   that
\begin{align*}
\mr{rk}(\langle \mc F\sqcup \mc G \rangle) \le \mr{rk}((\langle \mc F\rangle \sqcup_a \langle \mc G\rangle )\odot_T  (\langle \mc F\rangle \sqcup_c \langle \mc G\rangle \sqcup_c [T]^{\le 1} )) < \iota(\mr{srk}(\mc F\sqcup \mc G)).
\end{align*}
On the other hand, by Proposition  \ref{loi21950e}   $\mr{Is}(\langle \mc B \rangle)= \mr{Is}(\mc B)=  \mr{Is}(\mc F)\sqcup_a \mr{Is}(\mc G) \in \mk B_a$, whether $\mk B^a$ is compact or not. And,  
$\langle \mc F \sqcup \mc G\rangle \cap \mr{Ch}_c \con (\langle \mc F\rangle \cap \mr{Ch}_c) \sqcup_c (\langle \mc G\rangle \cap \mr{Ch}_c) \sqcup_c [T]^{\le 1} \in \mk B_c.$     If $\mr{Ch}_c$ is compact, then $\langle \mc F \sqcup \mc G\rangle \cap \mr{Ch}_c \in \mk B^c$. When $\mr{Ch}_c$ is non compact, since
\begin{align*}
   \iota(\mr{srk}(\langle \mc F \sqcup \mc G\rangle \cap \mr{Ch}_c )) =& \max\{\iota( \mr{srk}_c(\mc F \cap \mr{Ch}_c)), \iota( \mr{srk}_c(\mc G \cap \mr{Ch}_c))\} =\\
   =&  \iota(\mr{srk}_c ( (\langle \mc F\rangle \cap \mr{Ch}_c) \sqcup_c (\langle \mc G\rangle \cap \mr{Ch}_c) \sqcup_c [T]^{\le 1}))
\end{align*}
it follows from the property {\it (B.2)} of $\mk B_c$ that $\langle \mc F \sqcup \mc G\rangle \cap \mr{Ch}_c \in
\mk B_c$, hence $\mc B$ satisfies {\it(BT.2)}. This ends the proof of property {\it (B.2')} of $\mk B$.  \textit{(B.3)} is the content of Lemma \ref{kjjknkj744}.
\fprue

\subsection{The operation $\odot_T$ and ranks}\label{iu4787877555}

We  compute an upper bound of the rank of the family $\mc A\odot_T \mc C$ in terms of the ranks of $\mc A$
and $\mc C$, respectively.      Fix a tree $T$, a compact and hereditary family $\mc C$ on chains of $T$  and
a compact and hereditary family on immediate successors of nodes of $T$.   As we have observed in
\eqref{uuueuiewe} for finite trees,  it is natural to expect an upper bound of the rank of $\mc A \odot_T \mc
C$ by an exponential-like function of the rank of $\mc A$ and the rank of $\mc C$.

\defi
Given
  a countable ordinal number $\al$,   we define a function $f_\al: \ou\to\ou$ as follows:
\begin{align*}
f_\al(0):=& 1   \\
f_\al(\xi+1):= & f_\al(\xi) \cdot (\max\{\al,\xi\} \cdot \om )\\
f_\al(\xi):=& \sup_{\eta<\xi} f_\al(\eta), \text{  when $\xi$ is limit.}
\end{align*}
\fdefi
\nota \label{oi43iir44fdgfvcsaa}
\begin{enumerate}[(a)]
\item  $f_\al$ is a continuous  strictly increasing mapping such that $f_\al(\xi)$ is sum- indecomposable \hyphenation{in-de-com-po-s-able}for every $\al$ and $\xi$.
\item $f_\al(\xi)\ge (\al \cdot \om)^\xi$ always, and if $\xi\le \al$ then $f_\al(\xi)=(\al \cdot
\om)^\xi$.

\item Suppose that  $\la>\om$. Then  $f_\al(\xi)<\la$ for every $\al,\xi<\la$ if and only if $\la$ is exp-indecomposable:   Suppose that  $\la$ is closed under $f_{\cdot}(\cdot)$. Let $\al,\xi<\la$. Then $\al^\xi\le  f_\la(\xi) <\al$.
 Suppose that $\la$ is  exp-indecomposable. Since $\la>\om$, this is equivalent to saying that $\om^\la=\la$.   Let $\al,\xi<\la$. Since $f_{\cdot}(\cdot)$ is increasing in both
 variables, and since $\la$ is product-indecomposable, w.l.o.g. we assume that $\al$ is
 sum-indecomposable, i.e. $\al=\om^{\al_0}$, and   $\xi\le \al<\la$. Then,
\begin{equation*}
f_\al(\xi)=(\al\cdot \om)^\xi=\om^{(\al_0 +1)\cdot \xi}<\om^\la=\la.
\end{equation*}

\end{enumerate}

\fnota
Lemma  \ref{oi34r49778576} follows from Remark \ref{oi43iir44fdgfvcsaa} (c) and the following.  
\prop\label{o38989332rr}
Suppose that $\mc F=\langle \mc F\rangle$ and $\mr{Is}(\mc F)$ are compact. Then  
$$\mr{rk}(\mc F)< f_{\mr{rk}(\mr{Is}(\mc F))+1}(\mr{rk}(\mc F\cap \mr{Ch}_c) \cdot 3+2).$$
\fprop
The proof of this result has serveral parts. Recall that given a family $\mc F$,   $\langle \mc F\rangle_T=\conj{\langle s\rangle}{s\in \mc F}$ and $\langle \mc F\rangle_{\mr{is}-T}=\conj{\langle s\rangle_\mr{is}}{s\in \mc F}$.
\lema\label{k2178rhgukt8455}
Suppose that $\langle \mc F\rangle$ is countably ranked. Then
$$\mr{rk}(\langle\mc F\rangle_T)=\mr{rk}(\langle\mc F\rangle)\le \mr{rk}(\langle\mc F\rangle^\mr{is})=\mr{rk}(\langle\mc F\rangle_{\mr{is}-T}).$$
\flema
\prue
This follows from the following. 
\clam
 The following are equivalent for $x\con T$ and $\xi<\ou$:
\begin{enumerate}[(i)]
\item $x\in (\langle \mc F\rangle)^{(\xi)}$.
\item there is some subtree $\tau$ containing $x$  such that $\tau\in (\langle \mc F\rangle_T)^{(\xi)}$.
\item there is some subtree $\tau$ containing $x$ such that $\tau\in (\langle \mc F\rangle_{\mr{is}-T})^{(\xi)}$.
\end{enumerate}
\fclam
\prucl
We only give the proof of the equivalence between (i) and (ii). The equivalence between (i) and (iii) is proved in a similar way. First of all, by definition $x\in \langle \mc F\rangle$ if and only if $\langle x\rangle \in  \langle \mc F\rangle_T$.  Since $\langle \mc F\rangle$ is hereditary, each derivated set $(\langle \mc F\rangle)^{(\eta)}$ is also, so if there is some $\tau$ containing $x$ such that $\tau\in (\langle \mc F\rangle_T)^{(\xi)}\con  (\langle \mc F\rangle)^{(\xi)}$, then $x\in  (\langle \mc F\rangle)^{(\xi)}$. Now suppose that $x\in  (\langle \mc F\rangle)^{(\xi)}$, and let us prove that there is a subtree $\tau$ containing $x$ such that $\tau\in  (\langle \mc F\rangle_T)^{(\xi)}$. The case $\xi=0$ was treated above. 
Suppose that $\xi$  is limit, and let $(\xi_n)_n$ be an increasing sequence with supremum $\xi$. By inductive hypothesis, for every $n$ there is some subtree $\tau_n$ of $T$ such that $x\con \tau_n$ and such that $\tau_n\in (\langle \mc F\rangle_T )^{(\xi_n)}$. By compactness, there is an infinite set $M$ such that $(\tau_n)_{n\in M}$ is a $\De$-sequence with root $\tau$. A limit of subtrees is a subtree, hence $\tau$ is a subtree that contains $x$ and $\tau\in \bigcap_{n\in M} (\langle \mc F\rangle_T)^{(\xi_n)}$, so $\tau\in (\langle \mc F\rangle_T)^{(\xi)}$.  Suppose that $x\in (\langle \mc F\rangle_T)^{(\eta+1)}$. Choose a non-trivial  $\De$-sequence $(x_n)_n$ in   $(\langle \mc F\rangle)^{(\eta)}$ with limit $x$, and for each $n$ choose a subtree $\tau_n$ of $T$ containing $x_n$ and in $ (\langle \mc F\rangle_T)^{(\eta)}$. Now find an infinite subset $M\con \om$  such that $(\tau_n)_n$ is a $\De$-sequence with root $\tau$. Since $x_n\con \tau_n$, it follows that $(\tau_n)_{n\in M}$ is non-trivial, hence $x\con \tau\in (\langle \mc F\rangle_T)^{(\eta+1)}$.
\fprucl
\fprue

\defi
Given a subtree $U$ of $T$ with root $t_0$, let
$$\mr{stem}(U):=\conj{t\in U}{\text{every $u\in U$ is
comparable with $t$}}.$$
Given a chain $c$ and a family $\mc F$, let 
$$\mc F_c:=\conj{x\in \mc F}{ x \text{ is an $\mr{is}$-subtree of $T$ such that $c\ip\mr{stem}(x)$}}.$$
\fdefi
In the definition above, given two chains $c$ and $d$, $c\ip d$ means that $c$ is an initial part of $d$ with respect to the tree ordering $<$. 
Notice that $ \mr{stem}(U)$ is a non-empty chain in $U$, because $t_0\in \mr{stem}(U)$.

%

\lema\label{jki89yhkjmnff}
Suppose that $\mc F=  \langle \mc F\rangle^\mr{is}$ is compact, and suppose that $c$ is a $<$-chain of $T$.   If 
 $\mr{rk}(\mc F_c)\ge f_{\mr{rk}(\mr{Is}(\mc F))+1}(\xi)$, then $c\in (\mc F\cap \mr{Ch}_c)^{(\xi)}$. 
\flema
\prue
Set $\mc A:= \mr{Is}(\mc F)$ and $\mc C:=\mc F\cap \mr{Ch}_c$.  For each countable ordinal $\xi$, set $\be_\xi:=f_{\mr{rk}(\mc A)+1}(\xi)$.
 Fix a chain $c\in \mc C$. and suppose that   $\mr{rk}(\mc F_c)\ge \be_\xi$, and we have to prove that $c\in \mc C^{(\xi)}$.  The proof is by induction on $\xi$. The case $\xi=0$ or limit are trivial. Suppose that $\xi=\eta+1$.
 Since $\mc C$ is compact, we can assume without loss of generality that $c$ is maximal such that $\mr{rk}(\mc F_c)\ge \be_\xi$, i.e.
\begin{equation}\label{uui663333}
\text{if $c\cones c'\in \mc C$, then $\mr{rk}(\mc F_{c'})<\be_\xi$.  }
\end{equation}
If $c\neq \buit$  we define $t_c:=\max c$ and the following families.
\begin{align*}
\mc G:= & \conj{x\in \mc F_c}{\mr{Is}_{t_c}"(x)\con x }\\
\mc H:=& \conj{x\in \mc F_c}{x\cap \mr{Is}_{t_c}=\buit }.
\end{align*}
Both $\mc G$ and $\mc H$ are compact. Now since each $x\in \mc F_c$ is a $\mr{is}$-subtree, it follows that if $\#\mr{Is}_{t_c}'' x\ge 2$, then $\mr{Is}_{t_c}''x= \conj{t\we_\mr{is} u}{t\neq u\in x\setminus \{t_c\}}\con x$. Hence, for $x\in \mc F_c$, either $\mr{Is}_{t_c}''x\con x$ or $\mr{Is}_{t_c}''x\cap x=\buit$. This means that  $\mc F_c = \mc G\cup \mc H$.  Since $\be_\xi$ is sum-indecomposable, it follows that when $c\neq \buit$, then    $\max\{\mr{rk}(\mc
G),\mr{rk}(\mc H)\}\ge \be_\xi$, so we have the following two cases to consider.

 \noindent{\sc Case 1.} $c\neq \buit$ and $\mr{rk}(\mc G)\ge \be_\xi$. Let now
 \begin{align*}
 I:=&\conj{u\in \mr{Is}_{t_c}}{\mr{rk}(\mc F_{c\cup \{u\}})\ge \be_\eta}\\
 J:= & \mr{Is}_{t_c}\setminus I\\
\mc G_I:= &\conj{x\in \mc G}{\mr{Is}_{t_c}"x\in \mc A\rest I}\\
\mc G_J:= &\conj{x\in \mc G}{\mr{Is}_{t_c}"x\in \mc A\rest J}.
 \end{align*}
Clearly $\mc G\con \mc G_I \sqcup \mc G_J$. So, there are two subcases to consider.

\clam\label{claim1}	
$\mr{rk}(\mc G_I)\ge \be_\xi$ and $I$ is infinite.
\fclam
It follows from this that  given 	$u\in I$ we know by inductive hypothesis that $c\cup \{u\}\in \mc
C^{(\eta)}$, so $c\in \mc C^{(\xi)}$, as desired. 
\prue[Proof of Claim \ref{claim1}:]
We argue by contradiction. 
Suppose first that $\mr{rk}(\mc G_I)< \be_\xi$. It follows that $\mr{rk}(\mc G_J)\ge  \be_\xi$.     Let   $\la: \mc G_J\to \mc A\rest J$ be defined by
$\la(x):=\mr{Is}_{t_c}(x)$. This mapping is $\con$-increasing and since $\la(x)= x \cap \mr{Is}_{t_c}$ for
every $x\in \mc G_J$, it follows that $\la$ is continuous. By Proposition \ref{oi43jiro43ro34njrtguyy73} and  the
properties of $f_{\mr{rk}(\mc A)+1}(\cdot)$  we obtain that
$$\be_\eta \cdot (\mr{rk}(\mc A) \cdot \om)\le \be_\xi \le \mr{rk}(\mc G_J) <\sup_{y\in \mc A \rest J} (\mr{rk}\conj{x\in \mc G_J}{\la(x)=y}+1   )\cdot (\mr{rk}(\mc A\rest J)+1).$$
So there must be $y\in \mc A\rest J$ such that $\mr{rk}\conj{x\in \mc G_J}{\la(x)=y}\ge \be_\eta$.
We also have that
\begin{equation} \label{iuo32iuuidfttt}
\conj{x\in \mc G_J}{\la(x)=y}\con \bigsqcup_{u\in y}\mc F_{c\cup \{u\}}.
\end{equation}
Observe that $y\neq \buit$, because $\conj{x\in \mc G_J}{\la(x)=\buit}=\{c\}$ has rank 0. Hence, it follows
from \eqref{iuo32iuuidfttt} that there must be $u\in y$ such that $\mc F_{c\cup \{u\}}$ has rank at least
$\be_\eta$, contradicting the fact that $u\in J$. Finally, suppose that $\mr{rk}(\mc G_I)\ge \be_\xi$ but $I$ is finite. Then, $$\mc G_I \con \bigcup_{K\con I} \bigsqcup_{u\in K} \mc F_{c\cup \{u\}}  $$
it follows that there is some $u\in I$ such that $\mc F_{c\cup \{u\}}$ has rank at least $\be_\xi$,
contradicting \eqref{uui663333}.
\fprue

 \noindent{\sc Case 2.} $c=\buit$, or  $c\neq \buit $  and $\mr{rk}(\mc H)\ge \be_\xi$.   In order to unify the argument,  let 
 $\widetilde{\mc F}= \mc F$ when $c=\buit$, and let   $\widetilde{\mc F}=\mc H$ if $c\neq \buit$.

  Since $\mr{rk}(\widetilde{\mc F})\ge \be_\xi$, given a $\be_\xi$-uniform family $\mc B$ we can use Proposition \ref{njewirjiowejrew} to  can find  $f: \mc B\to \widetilde{\mc
 F}$   continuous, 1-1 and $(\ip,\con)$-increasing.  In particular,  this means that $c\cones f(s)$ for every $s\neq \buit$.  Observe that if $c\cones x\in \widetilde{\mc F}$, then there is $\min (x\setminus c)$: When $c=\buit$, then for such $x$, $\min(x\setminus c)$ is the root of $x$; if $c\neq \buit$ and $c\cones x\in \mc H$, then let $t\in x\setminus c$, and let $u:=\mr{Is}_{\max c}(t)$. Then for any other $t'\in x\setminus c$ we have that $\mr{Is}_{\max c}(t')=u$, since otherwise $u=t\we_\mr{is} t'\in x$, contradicting the fact that $x\in \mc H$.  Hence $\min (x\setminus c)= \bigwedge_{t\in x\setminus c} t$.   So, we can define $\la:\mc B\to \mc C_c$ as follows. Suppose that $c\cones f(\buit)$. We define for $s\in \mc B$ 
 $$\la(s):=\conj{t\in f(s)}{t\le \min (f(\buit)\setminus c)}.$$
Suppose that $f(\buit)=c$. Let $\la(\buit)=c$; and for each $\buit \cones s$, let  
 $$\la(s):=\conj{t\in f(s)}{t\le \min (f(\{\min s\})\setminus c)}.$$ 
Then, $\la$ is well-defined, that is, $\la(s)\in \mc C_c$ for every $s\in \mc B$.  It is also clear that 
 is $(\ip,\con)$-increasing. Moreover $\la$ is continuous: Suppose that $(s_n)_n$ is a $\De$-system with root $s$ such that $s<s_m \setminus s <   s_n\setminus s $ for every $m<n$.  Suppose first that $f(s)\neq c$  It follows that $\la(s_n)=f(s_n)\cap [0, \min(f(s)\setminus c)]$ and $\la(s)=f(s)\cap [0,\min(f(s)\setminus c)]$. Since $f$  and the intersection operation are continuous, it follows that $\la(s_n)\to_n \la(s)$. Suppose that $f(s)=c$. This implies that $s=\buit$, and $\la(s)=c$. Also,  $\la(s_n):=f(s_n)\cap [0,  \min (f(\{\min s_n\})\setminus c)]$. It follows that $\la(s_n)\to_n c=\la(\buit)$

For every $s\in \mc B^\mr{max}$ let $\vphi(s)$ be the
 maximal initial part $u$ of $s$ such that $\la(u)=\la( \buit)$, and let $M\con \om$ be infinite 
 such that $\mc B_0 := \vphi"(\mc B\rest M)$ is a $\ga$-uniform family on $M$ for some $\ga\le \be_\xi$.

\clam	
$\ga<\be_\xi$.
\fclam 
\prucl	
Suppose otherwise that $\ga=\be_\xi$. Since $c\cones \la(s)$ for every $s\neq \buit$ in $\mc B$, and
 since $\mc B_0\neq \{\buit\}$, it follows that $c\cones \la(s)=\la(\buit)=:c'$. By the
definition of $\mc B_0$ it follows that the restriction of $f$ to $\mc B_0$ satisfies that $f(x)\in \mc
F_{c'}$.  Consequently, $\mr{rk}(\mc F_{c'})\ge \ga=\be_\xi$, contradicting \eqref{uui663333}.
\fprucl
So, $\ga<\be_\xi$. Let $\mc B_1$ and $\mc B_2$ be $\be_\eta$-uniform and
$\xi$-uniform families on $M$, respectively. Since by definition $\be_\xi \ge \be_\eta \cdot (\xi \cdot
\om)$, it follows that $\ga,\be_\eta \cdot \xi<\be_\xi$. Since $\be_\xi$ is sum-indecomposable, it follows
that $(\be_\eta \cdot \xi)+\ga<\be_\xi$. Hence, by the properties of the uniform families (Proposition
\ref{ioufhdfhvvbhjss}) we obtain that there is $N\con M$ such that
$$  ((\mc B_1\otimes \mc B_2)\oplus \mc B_0)\rest N\con \mc B.$$
Fix now $s\in (\mc B_0\rest N)^\mr{max}$, set $N_s:= N/s$, and $\la_s: (\mc B_1 \otimes \mc B_2)\rest N_s\to
\mc C_c$, $\la_s(r):= \la(s\cup r)$.
Going towards a contradiction, suppose that  $c\notin \mc C^{(\xi)}$. Then it follows that $\mr{rk}(\mc C_c)<\xi$.
Since $\la_s$ is $(\ip,\con)$-increasing and continuous, it follows from Lemma \ref{778tr5643367898} that
there is some $s<r$ finite and some $r<P$ infinite, $r\cup P\con N_s$, such that $\{r\}\sqcup \mc B_1\rest P \con
(\mc B_1\otimes \mc B_2)\rest N$ and such that $\la_s$ is constant on $\{r\}\sqcup \mc B_1\rest P$ with value
$d:=\la(s\cup r)\in \mc C_c$. Since $\mc B_1\rest P$ contains non-empty elements $q$, it follows that
$s\sqsubset s \cup r \cup q$, so, since  $s\in \mc B_0^\mr{max}$,
\begin{equation*}
c\con \la(\buit)=\la(s)\cones \la(s\cup r\cup q)=\la(s\cup r).
\end{equation*}
On the other hand, the mapping $f_0: \mc B_2\rest P \to \widetilde{\mc F}_d\con \mc F_{d}$, $f_0(q):= f(s\cup r\cup q)$, witnesses
that $\mr{rk}(\mc F_d)\ge \be_\eta$, so by inductive hypothesis, $d\in \mc C^{(\eta)}$.  In this way we can
find $s<r_0<r_1<\cdots <r_n<\cdots$ such that $\la(s\cup r_n)\in \mc C^{(\eta)}$ and
\begin{equation}\label{oi43ieiiowewe}
c\con \la(\buit)=\la(s)\cones \la(s\cup r_n)
\end{equation}
 for every $n$. Since
$\la $ is continuous and $s\cup r_n \to_n s$, it follows that $\la(s\cup r_n)\to_n \la(s)=\la(\buit)$ and
non-trivially, by \eqref{oi43ieiiowewe}.  Hence, $\la(\buit)\in \mc C^{(\xi)}$, and so $c\in \mc C^{(\xi)}$,
because $\mc C$ is hereditary and $c\con \la(\buit)$, contradicting our hypothesis. 
\fprue

\prue[{\sc Proof of Proposition \ref{o38989332rr}}] Set $\mc C:=\mc F\cap \mr{Ch}_c$, $\mc D:=\langle\mc F\rangle^\mr{is}\cap \mr{Ch}_c$ and $\mc A:=\mr{Is}(\mc F)=\mr{Is}(\langle \mc F\rangle)= \mr{Is}(\langle \mc F\rangle^\mr{is})$. 
It follows from  Lemma \ref{jki89yhkjmnff} that 
 $\mr{rk}(\langle \mc F\rangle^\mr{is}) < f_{\mr{rk}(\mc A)+1}(\mr{rk}(\mc D) +1))$,
 so the proof will be finished once we have the following.
\clam
$\mr{rk}(\mc D) \le \mr{rk}(\mc C) \cdot 3 +1$.
\fclam
\prucl
As in the proof of the Corollary \ref{kjniuhrieuiuhrerer22}, let   $\mr{Is}$ be the collection of all immediate successors of $T$, and let $\pi: \mr{Is}\to T$ be the mapping assigning to each $t\in \mr{Is}$ its immediate predecessor $\pi(t)$. Then $\mc D\con \mc C\sqcup \pi^{-1}(\mc C)$.    It follows from Proposition \ref{oi43jiro43ro34njrtguyy73}  that
\begin{equation}\label{lkmklmoij78755}
\mr{rk}(\la^{-1}(\mc C))< \sup_{x\in \mc C}(\mr{rk}(\conj{s\in \la^{-1}(\mc C)  }{\la(s)\con x})+1)   \cdot (\mr{rk}(\mc C)+1).
\end{equation}
Now, given $x\in \mc C$, let us see that $\mr{rk}(\conj{s\in \la^{-1}(\mc C)  }{\la(s) \con x})\le 1$: Set $\mc E=  \conj{s\in \la^{-1}(\mc C)  }{\la(s)\con x}$. Every $s\in \mr{Is}$ such that $\pi(s)\con x$ can be written as  $s= \{\mr{Is}_t(\max x)\}_{t\in \pi(s)\setminus \{\max x\}}\cup \bar s$ with $\bar s\con \mr{Is}_{\max x}$ such that $\# \bar s\le 1$.  We see that $\mc E' \con \conj{ \{ \mr{Is}_t(\max x)\}_{t\in y}}{y\con x\setminus \{\max x\}}$, so $\mc E'$ is finite and $\mr{rk}(\mc E)\le 1$:  Otherwise, suppose that $s\in \mc E'$ is such that  $s\cap \mr{Is}_{\max x}=\{u\}$. Let $(s_n)_n$ be a non-trivial $\De$-sequence in $\mc E$ with root $s$. It follows that $s_n\cap \mr{Is}_{\max x}=\{u\}$ for every $n$. Since $x$ is finite, we may assume that $\la(s_n)= y\con x$ for every $n$. Hence
$s_n=\{u\}\cup \{\mr{Is}_t(\max x)\}_{t\in y\setminus \{\max x\}}$, so $(s_n)_n$ is a constant sequence, a contradiction. 
 It follows from \eqref{lkmklmoij78755} that $\mr{rk}(\la^{-1}(\mc C))< 2 (\mr{rk}(\mc C)+1)=2\mr{rk}(\mc C)+2$, so  $\mr{rk}(\mc D) \le \mr{rk}(\mc C) \dot{+} 2\mr{rk}(\mc C) +1 \le \mr{rk}(\mc C) \cdot 3 +1$, as desired. 
\fprucl
\fprue

\subsection{Canonical form of sequences of finite subtrees}\label{iu48767457655}
We prove here Lemma \ref{kjjknkj744}, that is if $\mc F\in \mk B$ and $\mc H\in \mk S$, then for every
sequence $(s_n)_{n<\om}$  in $\mc F$ there is an infinite subset $M\con \om$ such that $\bigcup_{n\in x}
s_n\in \mc F\times \mc H$ for every $x\in \mc H\rest M$.  The proof is based on a combinatorial analysis of
sequences of finite subtrees of $T$, done in the next Lemma \ref{ijiwurijjsfsd},  that  uses crucially the
Ramsey property. This relation  between the Ramsey theory, uniform fronts and BQO-WQO theory of trees is well
studied and has produced fundamental results like Kruskal Theorem \cite{Kr} (see also Nash-Williams paper
\cite{Na1})  and Laver Theorem \cite{La}.

 We start with some simple analysis of the tree generated by two finite subtrees $\tau_0$ and $\tau_1$.

\defi
Given $\tau_0$, $\tau_1$ two finite subtrees of $T$ and $t\in \tau_0\cup \tau_1$, let
\begin{align*}
i(t):=&\min\conj{i\in 2}{t\in \tau_i},\\
\pi(\tau_0,\tau_1):=&\conj{w\in  \langle \tau_0\cup \tau_1\rangle}{(\tau_0\setminus \tau_1)_{\ge w}\neq \buit\text{ and } (\tau_1\setminus \tau_0)_{\ge w}\neq \buit},\\
\sig(\tau_0,\tau_1):=&\pi(\tau_0,\tau_1)_{\mr{max}},\\
\bar{\sig}(\tau_0,\tau_1):=& \conj{t_0\we t_1}{t_0\perp t_1\text{ are in }\tau_0\cup \tau_1\text{ and }t_0\we t_1\notin \tau_0\cup \tau_1}.
\end{align*}
 \fdefi

\defi
For each $w\in \sig(\tau_0,\tau_1)$, fix $t^0(w)\in (\tau_0\setminus \tau_1)_{\geq w}$ and $t^1(w)\in
(\tau_1\setminus \tau_0))_{\geq w}$ such that $w=t^0(w)\wedge t^1(w) \notin \tau_0\cup \tau_1$ and whenever
$w\in \bar\sig(\tau_0,\tau_1)$, then $t^0(w)\perp t^1(w)$
\fdefi

\prop\label{jweofnow23}
$\bar{\sig}(\tau_0, \tau_1)\con \sig(\tau_0,\tau_1)$.
\fprop
\prue
Clearly $\bar{\sig}(\tau_0, \tau_1)\con \pi(\tau_0,\tau_1)$, so given $w \in \bar\sig(\tau_0,\tau_1)$, let us
prove that $w$ is maximal there, so that $w \in \sig(\tau_0,\tau_1)$. Suppose otherwise that there is $w'\in
\sig(\tau_0,\tau_1)$ such that $w<w'$ and let us get a contradiction. If $\IS{w}{t^0(w)}=\IS{w}{w'}$, then
$w=t^0(w) \we t^1(w) =t^1(w')\we t^1(w)\in \tau_1$, a contradiction. Otherwise, $w=t^0(w) \we t^1(w)
 = t^0(w)\we t_0(w')\in \tau_0$,   contradicting the hypothesis.
\fprue
%
%
%
%
%

\defi
Given two finite subtrees $\tau_0,\tau_1$, let
\begin{align*}
\ro: = & \ \ro(\tau_0,\tau_1) =  \tau_0 \cap \tau_1\\
\bar{\ro}: = & \ \bar{\ro}(\tau_0,\tau_1):=\langle\tau_0\rangle_\mr{is}\cap \langle\tau_1\rangle_\mr{is}\\
\ro_0:= & \ \bar{\ro}\cup \{0\}\\
(\tau_0,\tau_1)_\infty:= & \ \conj{u\in (\ro_0)_\mr{max}}{ (\pi(\tau_0,\tau_1))_{\ge u}\neq \buit}.
\end{align*}
\fdefi

\prop\label{weewfdf}
For every $u\in (\tau_0,\tau_1)_\infty$ one has that $\#(\sig(\tau_0,\tau_1))_{\ge u}=1$.
\fprop
\prue
If $u\in (\tau_0,\tau_1)_\infty$, then clearly $(\sig(\tau_0,\tau_1))_{\ge u}\neq \emptyset$. Suppose there
are $w_0\neq w_1\in (\sig(\tau_0,\tau_1))_{\ge u}$ for some $u\in (\tau_0,\tau_1)_\infty$. Observe that
$w_0\perp w_1$, since both of them are maximal in $\pi(\tau_0,\tau_1)$. Hence,
$$ u\leq t^0(w_0)\we t^0(w_1)= t^1(w_0)\we t^1(w_1), $$
so that
$$ u< t^0(w_0)\we_\mr{is} t^0(w_1)= t^1(w_0)\we_\mr{is} t^1(w_1)\in \langle\tau_0\rangle_\mr{is}\cap \langle\tau_1\rangle_\mr{is} = \bar{\ro} \subseteq \ro_0, $$
contradicting the maximality of $u$ in $\ro_0$.
\fprue

\defi
For every $u\in (\tau_0,\tau_1)_\infty$, let $\varpi_{\tau_0,\tau_1}(u)$ be the unique element of
$(\sig(\tau_0,\tau_1))_{\ge u}$.
\fdefi

\prop\label{iwjriwejiwerewe}
For every $w\in \pi(\tau_0,\tau_1)$, either there is $u\in (\tau_0,\tau_1)_\infty$ such that $w\le
\varpi_{\tau_0,\tau_1}(u)$, or else there is $u\in (\ro_0)_\mr{max}$ such that $w< u$. Consequently,
$\overline{\sig}(\tau_0,\tau_1)\con \ran (\varpi_{\tau_0,\tau_1})$.
\fprop
\prue
Given $w \in \pi(\tau_0,\tau_1)$, suppose there is no $u\in (\ro_0)_\mr{max}$ such that $w < u$ and let
$$u:=\max\{v \in \ro_0: v \leq w \}.$$
Let us prove that $u$ is maximal in $\ro_0$ so that $w$ witnesses that $u \in (\tau_0,\tau_1)_\infty$.
Suppose by contradiction that there is $v \in (\ro_0)_\mr{max}$ such that $u<v$ and in particular, $v \in
\bar{\ro}$. Notice that the definition of $u$ implies that $v \perp w$. Hence, $u \le w\wedge v < w$, so that
$u < w \we_\mr{is} v \leq w$. But $w \we_\mr{is} v = t^0(w) \wedge_\mr{is} v = t^1(w) \wedge_\mr{is} v$, so
that $w \we_\mr{is} v \in \bar{\ro}$ and we get a contradiction with the maximality of $u$ below $w$. It
follows that $u\in (\tau_0,\tau_1)_\infty$ and $w\le \varpi_{\tau_0,\tau_1}(u)$.

Finally, suppose that $w\in \overline{\sig}(\tau_0,\tau_1)$. Then, $w\in \sig(\tau_0,\tau_1)$, by Proposition
\ref{jweofnow23}. It is easy to see from the definition of  $ \overline{\sig}(\tau_0,\tau_1)$ that there is
no  $u\in (\ro_0)_\mr{max}$ such that $w< u$. Hence $w\le \varpi_{\tau_0,\tau_1}(u)$ for some $u\in
(\tau_0,\tau_1)_\infty$ and it follows from the maximality of $w$ that $w= \varpi_{\tau_0,\tau_1}(u)$.
\fprue

The following result guarantees that the new points of the tree generated by two finite subtrees $\tau_0$ and
$\tau_1$ are given by the function $\varpi_{\tau_0,\tau_1}$ and hence, they are controlled by the maximal
elements of $(\tau_0,\tau_1)_\infty$.

\cor    \label{oi438y747y44}
$\langle \tau_0\cup \tau_1\rangle=\tau_0\cup \tau_1\cup \conj{\varpi_{\tau_0,\tau_1}(u)
}{u\in(\tau_0,\tau_1)_\infty}$.
\fcor
\prue
If $w\in \langle \tau_0\cup \tau_1\rangle\setminus (\tau_0\cup \tau_1)$, then there are $t_0 \in \tau_0
\setminus \tau_1$ and $t_1 \in \tau_1 \setminus \tau_0$ such that $w = t_0 \we t_1$ and notice that $t_0 \we
t_1 \in \bar{\sig}(\tau_0, \tau_1)$. Then, by Proposition \ref{iwjriwejiwerewe}, there is $u\in
(\tau_0,\tau_1)_\infty$ such that $t_0 \we t_1 = \varpi_{\tau_0,\tau_1}(u)$. The other inclusion follows
directly from the definitions.
\fprue

To prove that $\times$ is a multiplication we have to deal with the tree generated by a sequence of finite
subtrees. Given a sequence $(\tau_k)_k$ of finite subtrees and $M \subseteq \om$, let $\tau_M$ be the
subtree generated by $\bigcup_{k \in M} \tau_k$. In order to be able to control the chains and the immediate
successors of some $\tau_M$, we will first find some suitable infinite $M$ such that the subsequence
$(\tau_k)_{k \in M}$ has some uniformity respective to the new points of $\tau_M$. This is the content of the
next result, that can be seen as a generalization of Proposition \ref{4iojio4ojirjgf}, which guarantees the
existence of an infinite fan, chain or comb inside any infinite subset of a tree.

If we assume each $\tau_k$ to be a singleton $\{t_k\}$ and apply Proposition \ref{4iojio4ojirjgf} to get an
infinite $M$ such that $\{t_k: k \in M\}$ is a fan, a chain or a comb, then the corresponding tree $\tau_M$
is given by $\{t_k: k \in M\} \cup \{w\}$, $\{t_k: k \in M\}$ or $\{t_k: k \in M\} \cup \{w_k: k \in M\}$,
respectively. The case (2.1) corresponds to $\{t_k: k \in M\}$ being a comb, so that the new points
$\{\varpi_k: k \in M\}$ form a chain; case (2.2) corresponds to $\{t_k: k \in M\}$ being a fan with root $w$
which is the only new point; case (2.3) corresponds to $\{t_k: k \in M\}$ being a chain and no new points
($\varpi_k=t_k$); and case (2.4) corresponds to $t_k = t_{k'}$ for all $k, k' \in M$.

In the next, after refining the sequence to get a fixed $\tau_\infty$ and $\varpi_{i}(u):=
\varpi_{i,j}(u) = \varpi_{i,k}(u)$ for $i<j<k$, each of these four cases might happen for each of the
sequences of  points $(\varpi_k(u))_{k \in M}$.

\teor[Canonical form of sequences of subtrees]\label{ijiwurijjsfsd}
Suppose that $(\tau_k)_k$ is a sequence of finite subtrees of $T$ forming a $\Delta$-system with root $\ro$
and such that $(\langle \tau_k\rangle_\mr{is})_k$ forms a $\De$-system with finite root $\bar{\ro}$. Then
there is a subsequence $(\tau_k)_{k\in M}$ such that
\begin{enumerate}[(1)]
\item  For every $i\neq j$ and $k\neq l$ in $M$ one has that
\begin{align}
\tau_\infty:=&(\tau_i,\tau_j)_\infty=(\tau_k,\tau_l)_\infty
\end{align}
\item Let $u\in \tau_\infty$. For each $i<j$ write $\varpi_{i,j}(u):=\varpi_{\tau_i,\tau_j}(u)$.  Then $\varpi_i(u):=\varpi_{i,j}(u)=\varpi_{i,k}(u)$ for every $i<j<k$, and $\varpi_i(u)\le \varpi_j(u)$ for every $i\le j$.

\noindent Moreover, one of the following holds.
\begin{enumerate}[(2.1)]
\item   $\varpi_i(u)<\varpi_j(u)$ for every $i<j$ and $\varpi_i(u)\notin \bigcup_{k} \tau_k$ for every $i<j$.
\item  $w(u):=\varpi_i(u)=\varpi_j(u)\notin \bigcup_k \tau_k$ for every $i$.
\item   $\varpi_i(u)<\varpi_j(u)$ and $\varpi_i(u)\in \tau_i\setminus \ro$ for every $i<j$.
\item $w(u):= u = \varpi_i(u)=\varpi_j(u)\in \ro$ for every  $i<j$.
 \end{enumerate}
  \end{enumerate}
\fteor

\defi
We call a $\De$-sequence $(\tau_k)_k$ of root $\ro$ such that $(\langle\tau_k\rangle_\mr{is})_k$ is a
$\Delta$-sequence of root $\bar{\ro}$ satisfying (1) and (2) of Theorem \ref{ijiwurijjsfsd} above a \emph{well-placed} sequence.
In this case, let $\tau_\infty^0$ be the set of those $u\in \tau_\infty$ such
that $w(u)=\varpi_i(u)=\varpi_j(u)$ for every $i<j$ in $\om$ and let
$\tau_\infty^1:=\tau_\infty\setminus \tau_\infty^0$. For each $u\in \tau_\infty^1$, let $z(u):=\sup_{i\in
\om} \varpi_i(u)$.

Given $I\con \om$, we use the terminology $\tau_I$ to denote $\langle \bigcup_{k\in I}\tau_k \rangle$.
\fdefi

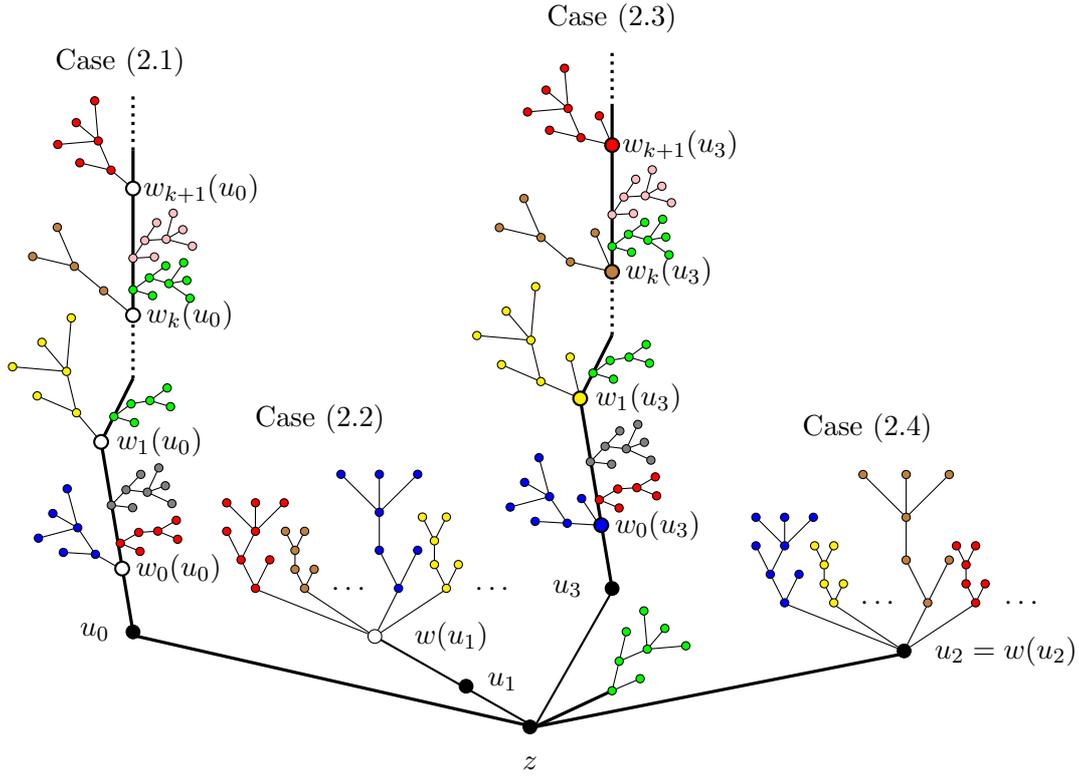
\begin{figure}[ht]

 \def\bla#1#2#3#4#5{
\begin{scope}[shift={(#1,#2)}]
\draw[rotate around={#4:(0,0)},scale=#5] node[circle, draw, thin,fill=#3, scale=0.3]  {}
    child { node[circle, draw, thin,fill=#3, scale=0.3] {} }
    child { node[circle, draw, thin,fill=#3, scale=0.3] {}
      child { node[circle, draw, thin,fill=#3, scale=0.3] {}
        child { node[circle, draw, thin,fill=#3, scale=0.3] {} }
        child { node[circle, draw, thin,fill=#3, scale=0.3] {} }
        child { node[circle, draw, thin,fill=#3, scale=0.3] {} } }
      child { node[circle, draw, thin,fill=#3, scale=0.3] {} } };
\end{scope}
}

 \def\blano#1#2#3#4#5{
\begin{scope}[shift={(#1,#2)}]
\draw[rotate around={#4:(0,0)},scale=#5] node[circle, draw, thin,fill=#3, scale=0.3]  {}
    child { node[circle, draw, thin,fill=#3, scale=0.3] {}
      child { node[circle, draw, thin,fill=#3, scale=0.3] {}
        child { node[circle, draw, thin,fill=#3, scale=0.3] {} }
        child { node[circle, draw, thin,fill=#3, scale=0.3] {} }
        child { node[circle, draw, thin,fill=#3, scale=0.3] {} } }
      child { node[circle, draw, thin,fill=#3, scale=0.3] {} } };
\end{scope}
}

 \def\blaa#1#2#3#4#5{
\begin{scope}[shift={(#1,#2)}]
\draw[fill=yellow,rotate around={#4:(0,0)},scale=#5] node[circle, draw, thin,fill=#3, scale=0.3]  {}
    child { node[circle, draw, thin,fill=#3, scale=0.3] {} }
    child { node[circle, draw, thin,fill=#3, scale=0.3] {}
      child { node[circle, draw, thin,fill=#3, scale=0.3] {}
        child { node[circle, draw, thin,fill=#3, scale=0.3] {} }
        child { node[circle, draw, thin,fill=#3, scale=0.3] {} } }
      };
\end{scope}
}
 \def\blaano#1#2#3#4#5{
\begin{scope}[shift={(#1,#2)}]
\draw[fill=yellow,rotate around={#4:(0,0)},scale=#5] node[circle, draw, thin,fill=#3, scale=0.3]  {}
    child { node[circle, draw, thin,fill=#3, scale=0.3] {}
      child { node[circle, draw, thin,fill=#3, scale=0.3] {}
        child { node[circle, draw, thin,fill=#3, scale=0.3] {} }
        child { node[circle, draw, thin,fill=#3, scale=0.3] {} } }
      };
\end{scope}
}

\def\blaaa#1#2#3#4#5{
\begin{scope}[shift={(#1,#2)}]
\draw[rotate around={#4:(0,0)},scale=#5] node[circle, draw, thin,fill=#3, scale=0.3]  {}
    child { node[circle, draw, thin,fill=#3, scale=0.3] {} }
    child { node[circle, draw, thin,fill=#3, scale=0.3] {}
      child { node[circle, draw, thin,fill=#3, scale=0.3] {}
        child { node[circle, draw, thin,fill=#3, scale=0.3] {} }
        child { node[circle, draw, thin,fill=#3, scale=0.3] {} }
        child { node[circle, draw, thin,fill=#3, scale=0.3] {} } }
       };
\end{scope}
}

 \begin{tikzpicture}[grow=up,scale=1.2]

 \begin{scope}[shift={(-1.15,-0.25)},scale=0.7] (comb1)
 \node[circle, draw, thin,fill=black, scale=.5]   at (5,-1) {};
\node (v30) at (4.4,-1) {$u_0$}; \node (v333) at (5,0) { \hspace{1.2cm}$w_0(u_0)$}; \node (v4) at (5,2) {
\hspace{0.7cm}$w_1(u_0)$}; \node (v70) at (4.8,8) {Case (2.1)};

\node (v6) at (3,0) {}; \node (v7) at (3,1) {}; \node (v8) at (3,2) {}; \node (v9)  at (3,3) {}; \node (v10)
at (3,4) {}; \node (v11) at (3,5) {}; \node (v11) at (3,6) {};

\draw[very thick] (5,-1) edge (4.5,2); \draw[very thick] (4.5,2) edge (5,3 ); \node (v20) at (5,4) {
\hspace{1.5cm}$w_k(u_0)$};

\draw[very thick] (5,4) edge (5,6.6); \draw[very thick,dotted] (5,3) edge (5,4); \node (v21) at (5,6) {
\hspace{1.8cm}$w_{k+1}(u_0)$}; \draw[very thick,dotted] (5,6.6) edge (5,7.5);

\blano{4.5}{2}{yellow}{40}{.4} \blano{4.8}{0}{blue}{60}{.3} \blaano{5}{4}{brown}{50}{.4} \blano{5}{6}{red}{50}{.3}
\bla{4.66}{1.0}{gray}{-70}{.2} \blaa{4.8}{.4}{red}{-86}{.2} \bla{5}{4.4}{green}{-80}{.19}
\bla{5}{4.9}{pink}{-60}{.2} \blaa{4.7}{2.4}{green}{-80}{.2} \node[circle, draw, thick, fill=white,
scale=0.5,overlay]  (v334) at (4.83,0) {}; \node[circle, draw, thick, fill=white, scale=0.5,overlay]  (v422)
at (4.5,2) {}; \node[circle, draw, thick, fill=white, scale=0.5,overlay]  (v202) at (5,4) {}; \node[circle,
draw, thick, fill=white, scale=0.5,overlay]  (v221) at (5,6) {};
 \end{scope}


\node (v430) at (10.4,1.3) {Case (2.4)};
\begin{scope}[shift={(10.8,-1.16)},scale=0.7] (fan2)
\node (v320) at (1.6,-0) {$u_2=w(u_2)$}; \draw[rotate around={0:(1,2)},scale=0.5] node[circle, draw,
thin,fill=black, scale=0.5]  {} child { node {$\cdots$} edge from parent[draw=none] }
    child { node[circle, draw, thin,fill=red, scale=0.1] {} }
      child { node[circle, draw, thin,fill=brown, scale=0.1] {} }
    child { node {$\cdots$} edge from parent[draw=none] }
     child { node[circle, draw, thin,fill=yellow, scale=0.1] {} }
     child { node[circle, draw, thin,fill=blue, scale=0.1] {} }
   ;
   \bla{-1.87}{0.76}{blue}{0}{.3}
\blaa{-1.1}{0.76}{yellow}{0}{.2} \blaaa{0.37}{0.76}{brown}{0}{.45} \blaa{1.12}{0.76}{red}{0}{.2}
\end{scope}

\draw[thick] (5,-1) edge (6.75,-2.0); \draw[thick]  (6.75,-2.0) edge (7.6,-0.5);
\node (v330) at (6.7,-2.4) {$z$}; \node[circle, draw, thin,fill=black, scale=0.5]  at (6.7,-2.) {};
\draw[very thick]  (6.75,-2.0) edge (7.6,-1.6); \bla{7.6}{-1.6}{green}{-40}{.2}

\draw[very thick] (2.4,-1) edge (6.75,-2.0); \node (v330) at (6.4,-1.5) {$u_1$}; \node[circle, draw,
thin,fill=black, scale=0.5]  at (6,-1.55) {}; \node (v330) at (4.4,1.4) {Case (2.2)};

\begin{scope}[shift={(5.,-1)},scale=0.7,rotate=0] (fan1)
\node (v320) at (1.2,0.0) {$w(u_1)$}; \draw[rotate around={0:(1,2)},scale=0.5] node[circle, draw,
thin,fill=white, scale=0.5]  {} child { node {$\cdots$} edge from parent[draw=none] }
    child { node[circle, draw, thin,fill=red, scale=0.1] {} }
      child { node[circle, draw, thin,fill=brown, scale=0.1] {} }
    child { node {$\cdots$} edge from parent[draw=none] }
     child { node[circle, draw, thin,fill=yellow, scale=0.1] {} }
     child { node[circle, draw, thin,fill=blue, scale=0.1] {} }
   ;
   \bla{-1.87}{0.76}{red}{0}{.3}
\blaa{-1.1}{0.76}{brown}{0}{.2} \blaaa{0.37}{0.76}{blue}{0}{.4} \blaa{1.12}{0.76}{yellow}{0}{.25}
\end{scope}
\draw[very thick]  (6.75,-2.0) edge (10.75,-1.2);

\begin{scope}[shift={(4.1,0.23)},scale=0.7] (comb1)
 \node[circle, draw, thin,fill=black, scale=.5]   at (5,-1) {};
\node (v30) at (4.3,-1) {$u_3$}; \node (v333) at (5,0) { \hspace{1.2cm}$w_0(u_3)$}; \node (v4) at (5,2) {
\hspace{0.7cm}$w_1(u_3)$}; \node (v80) at (5,8) {Case (2.3)};

\node (v6) at (3,0) {}; \node (v7) at (3,1) {}; \node (v8) at (3,2) {}; \node (v9)  at (3,3) {}; \node (v10)
at (3,4) {}; \node (v11) at (3,5) {}; \node (v11) at (3,6) {};

\draw[very thick] (5,-1) edge (4.5,2); \draw[very thick] (4.5,2) edge (5,3 ); \node (v20) at (5,4) {
\hspace{1.5cm}$w_k(u_3)$};

\draw[very thick] (5,4) edge (5,6.6); \draw[very thick,dotted] (5,3) edge (5,4); \node (v21) at (5,6) {
\hspace{1.8cm}$w_{k+1}(u_3)$}; \draw[very thick,dotted] (5,6.6) edge (5,7.5);

\bla{4.5}{2}{yellow}{40}{.4} \bla{4.8}{0}{blue}{60}{.3} \blaa{5}{4}{brown}{50}{.4} \bla{5}{6}{red}{50}{.3}
\bla{4.66}{1.0}{gray}{-70}{.2} \blaa{4.8}{.4}{red}{-86}{.2} \bla{5}{4.4}{green}{-80}{.19}
\bla{5}{4.9}{pink}{-60}{.2} \blaa{4.7}{2.4}{green}{-80}{.2} \node[circle, draw, thick, fill=blue,
scale=0.5,overlay]  (v334) at (4.83,0) {}; \node[circle, draw, thick, fill=yellow, scale=0.5,overlay]  (v422)
at (4.5,2) {}; \node[circle, draw, thick, fill=brown, scale=0.5,overlay]  (v202) at (5,4) {}; \node[circle,
draw, thick, fill=red, scale=0.5,overlay]  (v221) at (5,6) {};
 \end{scope}

\end{tikzpicture}
\caption{A well placed sequence $(\tau_k)_k$.}
\end{figure}
Each color in the figure corresponds to one of the elements of the sequence: blue nodes belong to the subtree
$\tau_0$, yellow nodes to $\tau_1$, and so on. Black is used to denote elements of the extended root
$\bar{\ro}$ and white, to nodes not belonging to any of the subtrees.

\prue[Proof of Theorem \ref{ijiwurijjsfsd}]
We will apply the Ramsey Theorem and refine the sequence $(\tau_k)_{k \in \om}$ finitely many times in order
to get the desired subsequence $(\tau_k)_{k \in M}$.

First, since $\bar{\ro}$ is finite and each $(\tau_i, \tau_j)_\infty \subseteq \ro_0 = \bar{\ro}\cup\{0\}$,
it follows from the Ramsey Theorem that we may assume, by passing to a subsequence $(\tau_k)_{k\in M}$, that
(1) holds. Now fix $u\in \tau_\infty$. Applying the Ramsey Theorem and passing again to a subsequence, we may assume
that exactly one of the following holds:
\begin{enumerate}[(a1)]
\item[(a1)] $\varpi_{i,j}(u) \notin \ro$ for every $i < j$ in $M$.
\item[(b1)] $\varpi_{i,j}(u) \in \ro$ for every $i < j$ in $M$.
\end{enumerate}

If (b1) holds, since $\ro$ is finite, we may pass to a further subsequence and get that $\varpi_{i,j}(u) =
\varpi_{k,l}(u)$ for every $i < j$ and $k<l$ in $M$. In particular we get (2.4). From now on we will assume that (a1) holds and prove that we have one of the other three cases (2.1), (2.2)
or (2.3).

\clam\label{kjasdnk}
For every $i<j$ in $M$ and $k \in M \setminus \{i,j\}$, one has that $\varpi_{i,j}(u)\notin \tau_k$.
\fclam
 \prucl
 We color a triple $i<j<k$ by $0$ if $\varpi_{i,j}(u)\in \tau_k$, by $1$ if $\varpi_{i,j}(u)\notin \tau_k$ and $\varpi_{i,k}(u)\in \tau_j$, by $2$ if $\varpi_{i,j}(u)\notin \tau_k$, $\varpi_{i,k}(u)\notin \tau_j$ and $\varpi_{j,k}(u)\in \tau_i$ and by $3$ otherwise. By the Ramsey Theorem, we may assume that all triples in $M$ are monochromatic. We prove that its color is $3$. In the other two cases, there are $i< j$ and $k\neq l$ such that
$$\varpi_{i,j}(u)\in \tau_k\cap \tau_l=\ro$$
which contradicts (a1).
 \fprue

For each $i<j$, let $t_{i,j}^i(u)\in \tau_i \setminus \ro$ and $t_{i,j}^j(u)\in \tau_j \setminus \ro$ be such
that
$$\varpi_{i,j}(u)=t_{i,j}^i(u)\we t_{i,j}^j(u).$$
Since each $\tau_i$ is finite and $t_{i,j}^i(u) \in \tau_i$, we may assume by the Ramsey Theorem that
$$t_i(u):=t_{i,j}^i(u)=t_{i,k}^i(u) \text{ for every }i<j<k.$$
Hence, for each $i<j<k$ in $M$, $\varpi_{i,j}(u)$ and $\varpi_{i,k}(u)$ are comparable since they are both
below $t_i(u)$.

\clam
By passing to an infinite subset of $M$, we may assume that $\varpi_{i,j}(u)=\varpi_{i,k}(u)$ for every
$i<j<k$ in $M$.
\fclam
\prucl
By the Ramsey Theorem, we may assume that one of the following holds:
 \begin{enumerate}[(a2)]
 \item[(a2)] $\varpi_{i,j}(u)=\varpi_{i,k}(u)$ for every $i<j<k$ in $M$.
 \item[(b2)] $\varpi_{i,j}(u)<\varpi_{i,k}(u)$ for every $i<j<k$ in $M$.
 \item[(c2)] $\varpi_{i,j}(u)>\varpi_{i,k}(u)$ for every $i<j<k$ in $M$.
 \end{enumerate}
Notice that (c2) is impossible, since trees have no infinite strictly decreasing chains. We claim that (b2) is also impossible and therefore, (a2) holds. Given $i<j<k$, if (b2) holds, then $\varpi_{i,j}(u) \leq t^j_{i,j}(u), t^k_{i,k}(u)$, so that $\varpi_{i,j}(u) \in \pi(\tau_j,\tau_k)$. By the maximality of $\varpi_{j,k}(u)$ in $\pi(\tau_j, \tau_k)$, we get that $\varpi_{i,j}(u) \leq \varpi_{j,k}(u)$. Then, $\varpi_{i,j}(u) \leq \varpi_{i,k}(u) \we \varpi_{j,k}(u)$. 

If $\varpi_{i,j}(u) < \varpi_{i,k}(u) \we \varpi_{j, k}(u)$, then the fact that $\varpi_{i,k}(u) \we
\varpi_{j, k}(u) \in \pi(\tau_i,\tau_j)$ contradicts the fact that $\varpi_{i,j}(u)$ is maximal in
$\pi(\tau_i,\tau_j)$. If $\varpi_{i,j}(u) = \varpi_{i,k}(u) \we \varpi_{j, k}(u)$, then we get that
$\varpi_{i,j} \in \tau_k$, which is a contradiction with Claim \ref{kjasdnk}. Therefore, (b2) cannot be true
and we conclude that (a2) holds.
\fprue

Let now $\varpi_i(u):=\varpi_{i,j}(u)$ for every $i<j$.
\clam
By passing to an infinite subset of $M$, we may assume that either $\varpi_i(u)=\varpi_j(u)$ for every $i<j$
in $M$, or $\varpi_i(u)<\varpi_j(u)$ for every $i<j$ in $M$.
\fclam
\prucl
By the Ramsey Theorem, we may assume that one of the following holds:
\begin{enumerate}[(a3)]
\item[(a3)] $\varpi_i(u)=\varpi_j(u)$ for every $i<j$ in $M$.
\item[(b3)] $\varpi_i(u)<\varpi_j(u)$ for every $i<j$ in $M$.
\item[(c3)] $\varpi_i(u)>\varpi_j(u)$ for every $i<j$ in $M$.
\item[(d3)] $\varpi_i(u)$ and $\varpi_j(u)$ are incompatible for every $i<j$ in $M$.
\end{enumerate}

Again (c3) is impossible, since trees have no infinite strictly decreasing chains. We claim that (d3) is also
impossible and therefore, either (a3) or (b3) holds. If (d3) holds, then we have that for $i<j<k<l$, $u \leq
\varpi_i(u) = t_{i,k}^k(u) \we t_{j,k}^k(u) <  t_{i,k}^k(u) \we_\mr{is} t_{j,k}^k(u) \in
\langle\tau_k\rangle_\mr{is}$ and $t_{i,k}^k(u) \we_\mr{is} t_{j,k}^k(u) = t_{i,l}^l(u) \we_\mr{is}(u)
t_{j,l}^l \in \langle\tau_l\rangle_\mr{is}$, contradicting the maximality of $u$ in $\ro_0$. Hence, either
(a3) or (b3) holds.
\fprue

In any case, by Claim \ref{kjasdnk} we may assume that $\varpi_i(u) \notin \tau_k$ for $i \neq k$. Hence, by
the Ramsey Theorem, we may assume that one of the following holds:
\begin{enumerate}[(a4)]
\item[(a4)] $\varpi_i(u) \notin \bigcup_k \tau_k$ for every $i$ in $M$.
\item[(b4)] $\varpi_i(u) \in \tau_i \setminus \ro$ for every $i$ in $M$.
\end{enumerate}

Now, if (a3) holds, (b4) cannot hold and we get that $u$ satisfies (2.2). If (b3) and (a4) hold, we get (2.1)
and if (b3) and (b4) hold, then we get (2.3). This concludes the proof of the theorem.
\fprue

\cor \label{iudfhakjds}
Given a well-placed sequence $(\tau_k)_{k < \om}$, $I \subseteq \om$ and $u \in \tau_\infty^1$, we have that:
\begin{enumerate}[(i)]
\item For every $t \in [u,z(u)]$, there is $i \in I$ such that if $t'\in (\tau_I)_{>t}$ with $\mr{Is}_t(t') \neq \mr{Is}_t(z(u))$, then $t'\in \tau_{i}$.
\item  $(\tau_I)_{< z(u)} \subseteq \bigcup_{k \in I} (\tau_k\cup\{\varpi_k(u)\})$.
\end{enumerate}
\fcor
\prue
(i): Let $u \leq t \leq z(u)$ and suppose there are $i_0<i_1$ in $I$ and $t_j \in (\tau_{i_j} \setminus
\ro)_{> t}$ such that $\mr{Is}_t(t_j) \neq \mr{Is}_t(z(u))$, $j=0,1$. Let $w = t_0 \we t_1$ and notice that
$w \in \pi(\tau_{i_0},\tau_{i_1})$. By Proposition \ref{iwjriwejiwerewe}, either there is $v \in \tau_\infty$
such that $w \leq \varpi_{i_0,i_1}(v)$ or there is $v \in (\ro_0)_\mr{max}$ such that $w < v$. Since $u \in
(\ro)_\mr{max}$ and $u \leq t \leq w$, the second alternative cannot hold and the first alternative holds
with $u=v$. Since $w\leq \varpi_{i_0,i_1}(u) = \varpi_{i_0}(u) <z(u)$, it follows that $w=t$. But then, $t=
\varpi_{i_1}(u) \we t_1 \in \tau_{i_1}$, which is impossible both in cases (2.1) and (2.3). Finally, notice
that if $t' \in (\tau_I)_{>t} \setminus \bigcup_{k \in I} \tau_k$ is such that $\mr{Is}_t(t') \neq
\mr{Is}_t(z(u))$, then there are $i_0<i_1$ in $I$ and $t_j \in (\tau_{i_j} \setminus \ro)_{> t'}$, so that
that $\mr{Is}_t(t_j) = \mr{Is}_t(t') \neq \mr{Is}_t(z(u))$, $j=0,1$, which we just proved that cannot happen.

(ii): If $t \in \tau_I \setminus \bigcup_{k \in I} \tau_k$, by Corollary \ref{oi438y747y44}, there are
$i_0<i_1$ in $I$ such that $t =\varpi_{i_0,i_1}(v)$ for some $v \in \tau_\infty$. If $t \leq z(u)$, the
maximality of $u$ and $v$ guarantee that $u=v$. Hence, $t=\varpi_{i_0,i_1}(u)=\varpi_{i_0}(u)$
\fprue

\cor \label{iudfhakjds2}
Given a well-placed sequence $(\tau_k)_{k < \om}$, $I \subseteq \om$ and $u \in \tau_\infty^0$, we have that:
\begin{enumerate}[(i)]
\item For every $t \in \mr{Is}_{w(u)}$ there is a $i \in I$ such that $(\tau_I)_{\geq t} \subseteq \tau_{i}$.
\item For every $t \in [u,w(u)[$ there is $i \in I$ such that if $t'\in (\tau_\om)_{>t}$ is such that  $\mr{Is}_t(t') \neq \mr{Is}_t(w(u))$, then $t'\in \tau_{i}$.
\item $(\tau_I)_{<w(u)} \subseteq \bigcup_{k \in I} \tau_k$.
\end{enumerate}
\fcor
\prue
(i): Let $t \in \mr{Is}_{w(u)}$ and suppose there are $i_0<i_1$ in $I$ and $t_j \in (\tau_{i_j} \setminus
\ro)_{\geq t}$, $j=0,1$. Let $w = t_0 \we t_1$ and notice that $w \in \pi(\tau_{i_0},\tau_{i_1})$. By
Proposition \ref{iwjriwejiwerewe}, either there is $v \in \tau_\infty$ such that $w \leq \varpi_{i_0,i_1}(v)$
or there is $v \in (\ro_0)_\mr{max}$ such that $w < v$. Since $u \in (\ro)_\mr{max}$ and $u \leq w(u) < t
\leq w$, the second alternative cannot hold and the first alternative holds with $u=v$, which cannot hold as
well, since $\varpi_{i_0}(u) = w(u) < t \leq w$. Finally, notice that if $t' \in (\tau_I)_{>t} \setminus
\bigcup_{k < \om} \tau_k$, then there are $i_0<i_1$ in $I$ and $t_j \in (\tau_{i_j} \setminus \ro)_{> t'}$,
so that that $t_j> t$, $j=0,1$, which we just proved that cannot happen.

(ii): Given $t \in [u, w(u)[$, suppose there are $i_0<i_1$ in $I$ and $t_j \in (\tau_{i_j} \setminus \ro)_{>
t}$ such that $\mr{Is}_t(t_j) \neq \mr{Is}_t(w(u))$, $j=0,1$. Let $w = t_0 \we t_1$ and notice that $w \in
\pi(\tau_{i_0},\tau_{i_1})$. By Proposition \ref{iwjriwejiwerewe}, either there is $v \in \tau_\infty$ such
that $w \leq \varpi_{i_0,i_1}(v)$ or there is $v \in (\ro_0)_\mr{max}$ such that $w < v$. Since $u \in
(\ro)_\mr{max}$ and $u \leq t \leq w$, the second alternative cannot hold and the first alternative holds
with $u=v$. Since $w\leq \varpi_{i_0,i_1}(u) = \varpi_{i_0}(u) = w(u)$, it follows that $w=t=w(u)$, a
contradiction. Finally, notice that if $t' \in (\tau_I)_{>t} \setminus \bigcup_{k \in I} \tau_k$ is such that
$\mr{Is}_t(t') \neq \mr{Is}_t(w(u))$, then there are $i_0<i_1$ in $I$ and $t_j \in (\tau_{i_j} \setminus
\ro)_{> t'}$, so that that $\mr{Is}_t(t_j) = \mr{Is}_t(t') \neq \mr{Is}_t(w(u))$, $j=0,1$, which we just
proved that cannot happen.

(iii): If $t \in \tau_I\setminus \bigcup_{k \in I} \tau_k$, by Corollary \ref{oi438y747y44}, there are
$i_0<i_1$ in $I$ such that $t =\varpi_{i_0,i_1}(v)$ for some $v \in \tau_\infty$. Hence, if $t < w(u)$, the
maximality of $u$ and $v$ guarantee that $u=v$. Hence, $t=\varpi_{i_0,i_1}(u)=\varpi_{i_0}(u) = w(u)$, a
contradiction.
\fprue

\cor \label{io4jioji34t345thh2}
Let $(\tau_k)_{k\in \om}$ be a well-placed sequence of finite subtrees of $T$. There are finite sets $\mk a,
\mk f\con T$ and, for each $z \in \mk f$, there is a chain $\{\varpi_i(z): i \in \om\}$ such that for
any finite $\emptyset \neq I \subseteq \om$, the following hold:
\begin{enumerate}[(i)]
\item For every $t\in \tau_I$, there are $z\in \mk f$ and $i\in I$ such that
\begin{equation*}
\#\left(  [t\we z, t] \cap (\tau_I \setminus \tau_i) \right)\le 1.
\end{equation*}
\item  For every $z\in \mk f$,
 \begin{equation*}
[0,z] \cap (\tau_I \setminus \bigcup_{i\in I}\tau_{i}) \subseteq \{\varpi_i(z): i \in I\}.
 \end{equation*}
\item  For every $t\in \tau_I \setminus \mk a$, there is $i\in I$ such that
\begin{equation}\label{jiuiu34iu4uyuty765}
\#(\IS{t}{\tau_I} \setminus \IS{t}{\tau_i})\le 1.
\end{equation} 
\end{enumerate}
\fcor
\prue
Let
\begin{align*}
\mk f:= &((\ro_0)_\mr{max}\setminus \tau_\infty)\cup \conj{w(u)}{u\in \tau_\infty^0}\cup \conj{z(u)}{u\in \tau_\infty^1}
\end{align*}
and let us prove  {\it (i)}: Given $t \in \tau_I$, if there is $u \in (\ro_0)_{\mr{max}}$ such that $t \leq u$, then either $u \in \mk f$,
or $w(u) \in \mk f$ or $z(u) \in \mk f$, depending on whether $u \in (\ro_0)_\mr{max} \setminus \tau_\infty$,
$u \in \tau_\infty^0$ or $u \in \tau_\infty^1$, respectively. In any case, there is $z \in \mk f$ such that
$t \we z = t$ so that property 1 holds trivially. Otherwise, there is a unique $u \in (\ro_0)_{\mr{max}}$
such that $u < t$. In case $u \notin \tau_\infty$, it follows from the definition of $\tau_\infty$ that $(\pi(\tau_i,\tau_j))_{\geq
u}=\emptyset$ for every $i \neq j$ in $I$, which implies property 1. If $u \in \tau_\infty^1$, then $z(u) \in \mk f$ and $t \we z(u) \in [u,z(u)]$. Then, Corollary
\ref{iudfhakjds}.(i) guarantees that there is $i \in \om$ such that $]t\we z(u), t] \subseteq \tau_i$, so that
property 1 holds. Finally, if $u \in \tau_\infty^0$, then $w(u) \in \mk f$ and $u \leq t \we w(u) \leq w(u)$. If $t \we w(u) < w(u)$,
then Corollary \ref{iudfhakjds2}.(ii) guarantees that there is $i \in \om$ such that $]t\we w(u), t] \subseteq
\tau_i$, so that condition 1 holds. If $t \we w(u) = w(u)$, it follows that $w(u) \leq t$ and the case when
$t=w(u)$ is trivial. If $w(u) < t$, Corollary \ref{iudfhakjds2}.(i) applied to $\mr{Is}_{w(u)}(t)$
guarantees that there is $i \in \om$ such that $]t\we w(u), t] \subseteq \tau_i$, so that property 1 holds and
this concludes the proof of {\it (i)}.

To prove {\it (ii)}, given $z \in \mk f$, let us consider three different cases. If $z = z(u)$ for some $u \in \tau_\infty^1$, let $(\varpi_i(z))_{i \in \om}$ be the sequence $(\varpi_i(u))_{i \in \om}$ and if $z = w(u)$ for some $u  \in \tau_\infty^0$, let $(\varpi_i(z))_{i \in \om}$ be the constant sequence equal
to $\varpi_i(u) = w(u)$. Finally, if $z \in (\ro_0)_\mr{max} \setminus \tau_\infty$, let $(\varpi_i(z))_{i \in \om}$ be the constant sequence equal
to $z$

Now, given $t \in \tau_I \setminus \bigcup_{i \in I} \tau_i$, by
Corollary \ref{oi438y747y44} there are $i_0<i_1$ in $I$ such that $t = \varpi_{i_0,i_1}(v)= \varpi_{i_0}(v)$ for some $v \in
\tau_\infty$ since the sequence is well-placed. Then, if $t \leq z$, the maximality of $v$ implies that $z \in \tau_\infty$ and since $u$ is also maximal, $u=v$ and $t = \varpi_{i_0}(v)= \varpi_{i_0}(u)$, which concludes the proof of {\it (ii)}

It remains to prove {\it (iii)}. Let
\begin{align*}
\mk a:=& \ro_0 \cup \conj{w(u)}{u\in \tau_\infty^0}
\end{align*}
and fix $t\in \tau_I \setminus \mk a$. Let $u \in (\ro_0)_\mr{max}$ be such that $u < t$ or $t < u$
(the equality cannot hold since $\ro_0 \subseteq \mk a$ and $t \notin \mk a$).

Suppose that $t < u$.  We see that there is at most one index $i$ such that $\mr{Is}_t"(\tau_{i}) \setminus \{\mr{Is}_t(u)\}\neq\buit$, that proves \eqref{jiuiu34iu4uyuty765}.  We suppose otherwise that there are $i_0<i_1$ such that $\mr{Is}_t"(\tau_{i_j}) \setminus \{\mr{Is}_t(u)\}\neq\buit$ for
$j =0,1$. Choose $t_j\in \mr{Is}_t(\tau_{i_j}) \setminus \{\mr{Is}_t(u)\}$, $j=0,1$. Then, $t_0,t_1\perp
u$, hence $\tau_{i_0}\ni t_0\wedge u =t=t_1\wedge u\in \tau_{i_1}$, and so $t\in \ro$, which is impossible.

Suppose now that $u < t$. If $u \notin \tau_\infty$, then \eqref{jiuiu34iu4uyuty765} holds trivially. If $u < t$ and $u \in
\tau_\infty^0$, we have to consider that cases $t < w(u)$ and $w(u) < t$ (again the equality cannot hold
since $w(u) \in \mk a$ and $t \notin \mk a$). If $t < w(u)$, then Corollary \ref{iudfhakjds2}.(ii) implies
that $\mr{Is}_t(\tau_I) = \mr{Is}_t"(\tau_i) \cup \{\mr{Is}_t(w(u))\}$ for some $i \in I$, so  \eqref{jiuiu34iu4uyuty765} 
holds. If $w(u)<t$, then Corollary \ref{iudfhakjds2}.(i) applied to $\mr{Is}_{w(u)}(t)$ implies that
$\mr{Is}_t"(\tau_I) = \mr{Is}_t(\tau_i)$ for some $i \in I$, so  \eqref{jiuiu34iu4uyuty765} holds. If $u < t$ and $u \in
\tau_\infty^1$, then Corollary \ref{iudfhakjds}.(i) implies that $\mr{Is}_t"(\tau_I) = \mr{Is}_t(\tau_i) \cup
\{\mr{Is}_t(z(u))\}$ for some $i \in I$, so  \eqref{jiuiu34iu4uyuty765} holds, including in case $z(u) < t$.
\fprue

\prue[{\sc Proof of Lemma \ref{kjjknkj744}}]
Fix $\mc F\in \mk B$ and $\mc H\in \mk S$ of infinite rank. We recall that  
$$ \mc G:= \mc F  \times \mc H =
((\mc A \times_a \mc H )\sqcup_a [T]^{\le 1}) \odot_T ((\mc C\times_c \mc H)\boxtimes_c 5)$$  where $\mc
A:= \mr{Is}(\langle \mc  F \rangle)$ and $\mc C:=\langle \mc F \rangle\cap \mr{Ch}_c$. We know from Proposition \ref{jkluiui44384} that $\iota(\mr{srk}(\mc F))=\iota(\mr{srk}_a(\mc A))$ when $\mr{Ch}_a$ is not compact and $\iota(\mr{srk}(\mc F))=\iota(\mr{srk}_c(\mc C))$, when $\mr{Ch}_c$ is not compact.  We check first that $\mc G\in \mk B$,  and for this purpose we use Proposition  \ref{ioioio32uihihu66}. Set 
$$\bar{ \mc A}:=(\mc A \times_a \mc H )\sqcup_a [T]^{\le 1}\in \mk B^a, \, \bar{\mc C}:=(\mc C\times_c \mc H) \boxtimes_c 5\in \mk B^c.$$ Then when $\mr{Ch}_a$ or $\mr{Ch}_c$ is not compact, then $\bar{\mc A}$ and $\bar{\mc C}$ fulfill the conditions {\it(a)}, {\it(b)} and {\it(c)} of that Proposition, so $\mc G\in \mk B$.  Suppose that $\mr{Ch}_a$ and $\mr{Ch}_c$ are both non compact. Then {\it(d)} of Proposition  \ref{ioioio32uihihu66}, i.e., $\iota(\mr{srk}_a(\bar{\mc A}))=\iota(\mr{srk}_c(\bar{\mc C}))$ holds because $\iota(\mr{srk}_a(\bar{\mc A}))=\iota(\mr{srk}(\mc F))=\iota(\mr{srk}_c(\bar{\mc C}))$.  Now we verify that $\mc G$ satisfies {\it(M.1)}: Since $\mc G\in \mk B$, it follows from Proposition \ref{jkluiui44384} that if $\mr{Ch}_a$ is non compact then 
\begin{align*}
 \iota(\mr{srk}(\mc G))=& \iota(\mr{srk}_a(\bar{\mc A}))=\iota(\mr{srk}_a(\mc A\times_a \mc H))=\max\{\iota(\mr{srk}_a(\mc A)),\iota(\mr{srk}(\mc H))\}=\\
 =&\max\{\iota(\mr{srk}(\mc F)),\iota(\mr{srk}(\mc H))\},
\end{align*} 
and if  $\mr{Ch}_c$ is non compact then 
\begin{align*}
 \iota(\mr{srk}(\mc G))=& \iota(\mr{srk}_c(\bar{\mc C}))=\iota(\mr{srk}_c(\mc C\times_c \mc H))=\max\{\iota(\mr{srk}_c(\mc C)),\iota(\mr{srk}(\mc H))\}=\\
 =&\max\{\iota(\mr{srk}(\mc F)),\iota(\mr{srk}(\mc H))\}.
\end{align*} 
This ends the proof of
property {\it(M.1)} for $\times$. Let us prove now {\it(M.2)} for $\times$.  Let $(\tau_k)_k$ be a sequence in $\langle \mc F \rangle$. Since $\langle\mc F\rangle\in \mk B$,  it follows that $\mr{Is}(\langle \mc F\rangle)\con \langle \mc F\rangle$, so $\langle \mc F\rangle^\mr{is}$ is compact, by Corollary \ref{kjniuhrieuiuhrerer22}. Hence we may assume  that each
$\tau_k$ is a tree and that the sequence $(\langle \tau_k \rangle_\mr{is})_k$ is a $\De$-sequence.  By
Corollary \ref{io4jioji34t345thh2} there is a subsequence $(\sig_k)_{k<\om}$ of $(\tau_k)_k$ and there are
$\mk a, \mk f$ finite subsets of $T$ such that 1., 2. and 3. there hold. By refining the subsequence
$(\nu_k)_{k <\om}$ finitely many times, we may assume that for every $I\in \mc H$, we have that:
\begin{enumerate}[(i)]
\item For all $z \in \mathfrak{f}$, $\{w_i(z): i\in I\} \in \mc C \times_c \mc H$.
\item  For all $z \in \mathfrak{f}$, $\bigcup_{i\in I}(\nu_i \cap [0,z])\in \mc C \times_c \mc H$.
\item  For all $t \in \mathfrak{a}$, $\bigcup_{i\in I}\mr{Is}_{t}"\nu_i\in  \mc A \times_a \mc H $.
\end{enumerate}
For (i) we use that $[T]^{\le 1}\con \mc F$. Fix $I\in \mc H\rest M$ and let us prove that $\nu_I \in \mc
F\times \mc H$, which is enough to guarantee that $\bigcup_{i\in I}\nu_i\in \mc F\times \mc H$,  since  this
family is hereditary.

\clam
For every $t\in \nu_I$ one has that $(\nu_I)_{\le t} \in  \bar{\mc C}$.
\fclam
\prucl
Given $t\in \nu_I$, by property 1. of the sequence $(\nu_k)_{k\in M}$, there are $z \in \mk f$, $\bar{t} \in
\nu_M$ and $i \in M$ such that
$$(\nu_M)\cap [t\we z, t]\subseteq \nu_i \cup \{\bar{t}\}.$$
Then, the property 2. of $(\nu_k)_{k\in M}$ implies that
$$\nu_I \cap [0,z] \con (\bigcup_{i \in I} (\nu_i \cap[0,z])) \cup  \{w_i(z): i\in I\}   \cup \{z\}.$$
Hence,
\begin{align*}
\nu_I \cap [0,t]\con & \nu_I\cap[t \we z,t] \cup (\nu_I \cap [0,z]) \con (\nu_i\cap [t \we z,t]) \cup  (\bigcup_{i \in I} (\nu_i  \cap[0,z])) \cup  \{w_i(z): i\in I\}  \cup \{\bar{t},z\}.
\end{align*}

Now notice that
\begin{enumerate}[$\bullet$]
\item $\nu_i\cap [t \we z,t] \in \mc C\con \mc C \times_c \mc H$;
\item $\bigcup_{i \in I} (\nu_i \cap[0,z]) \in \mc C \times_c \mc H$ by (ii) above;
\item $\{w_i(z): i\in I\} \in \mc C \times_c \mc H$ by (i) above;
\item $\{\bar{t},z \} \in [T]^{\le 2 }\con \mc C\boxtimes 2$.

\end{enumerate}
Putting all together, we obtain that
\[\nu_I\cap [0,t]\in(\mc C\times_c \mc H)\boxtimes_c 5=\bar{\mc C}.  \qedhere \popQED
\]
\let\qed\relax\fprue
\clam
For every $t \in \nu_I$ one has that $\mr{Is}_t"(\nu_I)\in \bar{\mc A}$.
\fclam
\prucl
Given $t \in \nu_I$, if $t\notin \mk a$, the property 2. of $(\nu_k)_{k\in M}$ implies that there are $j\in
I$ and $\bar t\in \mr{Is}_t$ such that
$$
\mr{Is}_t"(\nu_I)\con \mr{Is}_t"(\nu_j) \cup \{\bar t\}\in \mc A\sqcup_a [T]^{\le 1}.
$$
If $t\in \mk a$, it follows from (iii) above that
$$
\mr{Is}_t"(\nu_I) \con \bigcup_{i\in I} \mr{Is}"(\nu_i) \in \mc A \times_a \mc H
$$
In any case, we have that $\mr{Is}_t"(\nu_I)\in  \bar{\mc A}$.
\fprue
These two claims imply that $\nu_I \in\mc G$ for every $I\in \mc H$.
\fprue

\section{Bases of families on (not so) large cardinals}\label{bskdjgbsekjbgcnamnc}

The purpose of the section is to use Theorem \ref{4jtierjteoirtertrtr344} to prove Theorem \ref{firstMahlo}, that is, to show the existence of bases of families on every cardinal below the first Mahlo cardinal. 


\subsection{Binary trees}
We start by analyzing the case of binary trees. 
\teor\label{oiiojoui434344}
Suppose that $\ka$ has a basis. Then $2^\ka$ also has a basis.
\fteor

 \prue 
Suppose that $\ka$ has a basis. Let $T$ be the complete binary tree $2^{\le \ka}$.   The height mapping
$\mr{ht}:T\to \ka+1$ is strictly monotone, so
  it follows from Theorem  \ref{ui32iuui32t744} that there is a basis of families  on chains of $T$. Each set $\mr{Is}_t$ has size 2,
  so it follows from Theorem \ref{4jtierjteoirtertrtr344} that  $T$ has a basis. Since $T$ has cardinality $2^\ka$, we are done. 
\fprue

\nota
Let $\vep$ be the first exp-indecomposable ordinal $>\om$, that, is $\vep$ is the first ordinal such that $\om^\vep=\vep$. Then it is easy to see that $f_1(\al)=(\al\cdot
\om)^\al$ for every  $\al<\vep$.  Using this, and the construction of the basis on $2^\ka$ from the one in $\ka$ we can give upper
bounds of the ranks of $\om^\al$-homogeneous families in small exponential cardinals.  Let $\mr{ht}: 2^{\le \aleph_0}\to \aleph_0+1$  be the height function. We know that $\mr{ht}$ is adequate. For each $1\le \al<\vep$, let $\mc H_\al$ be a $\om^\al$-homogeneous family of rank exactly $\om^\al$ (e.g. the Schreier families on the index set $\aleph_0+1\sim \om$).  Let $\mc C_\al:= \mr{ht}^{-1}(\mc H_\al)$. It follows from Lemma \ref{kloi89787855} that $\mc C_\al$ is a $\be$-homogeneous family on chains of $T$ such that $\om^\al\le \be <\om(\om^\al+1)$ and such that $\mr{rk}(\mc C_\al)<\om(\om^\al+1)$. In fact, we can have a better estimation. From Proposition \ref{oi43jiro43ro34njrtguyy73}  it follows that 
$$\mr{rk}(\mc C_\al) <    \sup_{x\in \mc H_\al}(\mr{rk}(\conj{s\in \mc C_\al}{\mr{ht}(s)\con x})+1)   \cdot (\mr{rk}(\mc H_\al)+1).$$   
Let  $x\in \mc H_\al$, and set $\mc E:=\conj{s\in \mc C_\al}{\mr{ht}(s)\con x}$. We claim that $\mr{rk}(\mc E)\le 1$:  Suppose that $s\in \mc E'$ and let $(s_n)_n$ be a non trivial $\De$-sequence in $\mc E$ with root $s$. We claim that  $\aleph_0\notin \mr{ht}''(s)$ (recall that $\aleph_0+1=\aleph_0\cup \{\aleph_0\}$).  Otherwise, let $f\in s$ be the unique node such that $\mr{ht}(f)=\aleph_0$. Since each $s_n$ is a chain, it  follows that  $\conj{f\in s_n}{\mr{ht}(f)=\aleph_0}=\{f\}$, so
$$s_n \con \{ f \rest m \}_{m\in x\setminus \{\aleph_0\}}\cup \{f\}.$$ 
Then there is an infinite set $M\con \om$ such that   $s_n=s$  for every $n\in M$, a contradiction.  Let us see that $\mc E''=\buit$. Suppose otherwise, and let   $s\in \mc E''$ and suppose that  $(s_n)_n$ is a non-trivial $\De$-sequence in $\mc E'$ with root $s$.  Since $\conj{r}{ \mr{ht}''(r)\con x\setminus \{\aleph_0\}}$ is a finite set, it follows that    $M=\conj{n\in \om}{\aleph_0\in \mr{ht}''(s_n)}$ is co-finite, in particular non-empty. Let $n\in M$. Then $s_n\in \mc E'$ and $\aleph_0\in \mr{ht}''s_n$, contradicting the previous fact.   We have just proved that 
$$\mr{rk}(\mc C_\al) < 2(\om^\al+1)=\om^\al+2.$$
 Let $\mc F_\al:=[2^{\le \aleph_0}]_a^{\le 2}\odot_T (\mc C_\al\sqcup_c [2^{\le \aleph_0}]^{\le 1})$. We know that $\mc F_\al$ is an homogeneous family and $\om^{\al}+1\le \mr{srk}(\mc F_\al) \le \om^{\al}+2$. On the other hand, it follows from Proposition \ref{o38989332rr} that 
 $$\om^\al+1\le \mr{srk}(\mc F_\al)\le \mr{rk}(\mc F_\al)< f_{1}( (\om^\al+2)\cdot 3+2)= f_{1}( \om^\al\cdot 3+8)\le \om^{\om^\al\cdot 3 +\al \cdot 8 + 8 }.$$
One step up further, we use $\mc F_\al$ (on $2^{\aleph_0}\sim 2^{\le \aleph_0}$) and the height function  to find homogeneous families $\mc D_\al$  on chains of   $2^{\le 2^{\aleph_0}}$ with
$$\om^{\al}+1\le \mr{srk}_c(\mc D_\al) \le \mr{rk}(\mc D_\al)< \om  \cdot(\om^{\om^\al\cdot 3 +\al \cdot 8 + 8 }+1)= \om^{\om^\al\cdot 3 +\al \cdot 8 + 8 } +\om.$$ 
Set $T:= 2^{\le 2^{\aleph_0}}$, and let $\mc G_\al:= [T]_a^{\le 2} \odot_T \mc D_\al$. We know that this is an homogeneous family with 
$$\om^\al+1\le \mr{srk}(\mc G_\al)\le \mr{rk}(\mc G_\al)< \om^{\om^{\om^\al \cdot 3+ \al\cdot 8 +8}  +\om^{\al+1}  }$$ 
  \fnota

\subsection{Trees from walks on ordinals} \label{trees_walks} We pass now to study certain trees on inaccessible cardinal numbers. They are produced  using the method of walks on ordinals. We introduce some basic notions of this. For more details we refer the reader to the monograph \cite{To}.
\defi
A  \emph {$C$-sequence}  $\bar C:=(C_\al)_{\al<\theta}$ is a sequence such that $C_\al\con \al$ is a closed and unbounded subset of $\al$ with $\mr{otp}(C_\al)=\mr{cof}(\al)$. The \emph{$\bar C$-walk} from $\be$ to $\al<\be$ is the finite sequence of ordinals defined recursively by
\begin{align*}
\mr{Tr}(\al,\be):=&(\be)\conc \pi(\al,\min (C_\be\setminus \al))\\
\mr{Tr}(\al,\al):=&(\al).
\end{align*}
We write then the $\bar C$-walk as $\be=\pi_0(\al,\be)>\cdots>\pi_{l}(\al,\be)=\al$, where
$l+1=\mr{ht}(\pi(\al,\be)$, and for each $i\le l$, $\pi_i(\al,\be)$ is the $i^{\mr{th}}$ term of
$\pi(\al,\be)$. Let
$$\ro_2(\al,\be):= \mr{ht}(\mr{Tr}(\al,\be))-1.$$
We now define the mapping $\ro_0:[\theta]^2\to (\mc P(\theta))^{<\om}$ for $\al\le \be$ recursively by
\begin{align*}
\ro_0(\al,\be):=&(C_\be\cap \al)\conc \ro_0(\al,\min (C_\be\setminus \al))\\
\ro_0(\al,\al):=&\buit.
\end{align*}
Let $T=T(\ro_0)$ be the tree whose nodes are $\ro(\cdot,\be)\rest \al$, $ \al\le \be $, ordered by end-extension as functions.
Given $t\in T(\ro_0)$, let $\al_t\le \be_t$ be such that $t=\ro_0(\cdot,\be_t)\rest \al_t$. We say that $(\al_t,\be_t)$ represents $t$.
\fdefi
\prop\label{oii32985555}
$T$ has size $\theta$, and  if $\theta$ is strong limit, then for every $t\in T$ one has that $\#\mr{Is}_t(T)<\theta$.
\fprop
\prue
This is a tree on a quotient of $[\theta]^{\le 2}$, so it has cardinality $\theta$.   Now observe that the immediate successors of $t=\ro_0(\cdot,\be)\rest \al$ are  extensions $u$ of $t$ whose support is $\al+1$. It follows that the number of them is at most $(2^{\al})^{<\om}<\theta$, when we assume that $\theta$ is strong limit.
\fprue

In other words, the partial ordering $<_\mr{a}$ is the disjoint union of small partial orderings.

Observe that if $t<u$ in $T(\ro_0)$, then we can take $(\al_t,\be_t),(\al_u,\be_u)$ representing $t$ and $u$ respectively such that $\al_t<\al_u$ and $\be_t\le \be_u$: Take representatives $(\al_t,\be_t)$, $(\al_u,\be_u)$ of $t$ and $u$ respectively. Then $\al_t<\al_u$ and  $\ro_0(\cdot, \be_t)\rest \al_t=\ro_0(\cdot,\be_u)\rest \al_t$, hence $(\al_t,\min\{\be_t,\be_u\})$ is a representative of $t$ and satisfies the required condition together with $(\al_u,\be_u)$. The following is well-known.
\prop
$t<u$ if and only if $\ro_0(\al_t,\be_t)=\ro_0(\al_t,\be_u)$.\qed
\fprop

\defi
Given a C-sequence $\bar C$ on $\theta$, let
$$\mc I(\bar C):=\conj{C\con \theta}{  C\ip C_\al \text{ for some $\al<\theta$} }.$$
We consider $\mc I(\bar C)$ ordered by $\sqsubset$.
\fdefi

\prop\label{sijeirjewirowfsdfd}
$\ro_0: (T,<)\to \mc I(\bar C)_\mr{lex}^{<\om}$ is strictly monotone, and consequently $\ro_0: (T,<)\to \mc
I(\bar C)_\mr{qlex}^{<\om}$ is adequate.
\fprop
\prue
 Suppose that $t<u$ in $T(\ro_0)$.

\clam
Suppose that $\ro_0^i(t) = \ro_0^i(u)$ for every $i\le k$. Then $\pi_{i}(\al_t,\be_u)=\pi_{i}(\al_u,\be_u)$ for every $i\le k+1$  and $\ro_0^{k+1}(t)\ip \ro_0^{k+1}(u)$.
\fclam
 \prucl
 Induction on $k \ge 0$.  Suppose is true for $k-1$. Then $\pi_{i}(\al_t,\be_u)=\pi_{i}(\al_u,\be_u)$ for every $i\le k$.  It follows that
 \begin{align*}
 \pi_{k+1}(\al_t,\be_u)= &\min (C_{\pi_k(\al_t,\be_u)}\setminus \al_t)= \min (C_{\pi_k(\al_u,\be_u)}\setminus \al_t)\\
 C_{\pi_k(\al_t,\be_t)}\cap \al_t = & \ro_0^k(t)= \ro_0^k(u)=C_{\pi_k(\al_u,\be_u)}\cap \al_u.
 \end{align*}
In particular,
$ C_{\pi_k(\al_u,\be_u)}\cap [\al_t,\al_u[=\buit $ hence $ \min (C_{\pi_k(\al_u,\be_u)}\setminus \al_t)= \min (C_{\pi_k(\al_u,\be_u)}\setminus \al_u)$, so
$$ \pi_{k+1}(\al_t,\be_u)=\pi_{k+1}(\al_u,\be_u).$$
Finally,
\[\pushQED{\qed} 
\ro_0^{k+1}(t)=\ro_0^{k+1}(\al_t,\be_t)=\ro_0^{k+1}(\al_t,\be_u)= C_{\pi_{k+1}(\al_t,\be_u)}\cap \al_t =
C_{\pi_{k+1}(\al_u,\be_u)}\cap \al_t \ip \ro_0^{k+1}(u). \qedhere
\popQED
\] \let\qed\relax
 \fprucl
 It follows that $\ro_0^0(t)\ip \ro_0^0(u)$, so there must be $k<\ro_2(\al_u,\be_u)$ such that $\ro_0^{k}(t)\sqsubset \ro_0^k(u)$, since otherwise for every $k$ $\pi_k(\al_t,\be_u)=\pi_k(\al_u,\be_u)$ would imply that $\al_t=\al_u$.
 \fprue

So, the mapping $\ro_0$ is adequate, hence if  $\mc I(\bar C)_\mr{lex}^{<\om}$ has a basis of families on
chains, then $\ro_0$ will transfer it to a basis on $<$-chains.

\prue[{\sc Proof of Theorem \ref{ui32iuui32t7442212}}]

(a): Suppose that $(\mk B^\mc P,\times_\mc P)$ and $(\mk B^\mc Q, \times_\mc Q)$ are bases of families on
chains of $\mc P$ and $\mc Q$ respectively. Let $\mk B$ be the collection of all $\mc P \times_\mr{lex} \mc
Q$-homogeneous families $\mc F$ such that
\begin{enumerate}[(i)]
\item $\pi_\mc P"(\mc F):=\conj{\pi_\mc P"(x)}{x\in \mc F}\in \mk B_\mc P$.
\item  $(\mc F)_\mc P:=\conj{(x)_p}{x\in \mc F, \, p\in P}\in \mk B_\mc Q$.
 \end{enumerate}
Given $\mc F\in \mk B$ and $\mc H\in \mk S$, let
$$\mc F\times \mc H:= (\pi_\mc P"( \mc F )\times_\mc P \mc H) \circledast_\mr{F} ((\mc F)_\mc P \times_\mc Q\mc H).$$
We verify that $(\mk B, \times)$ is a pseudo-basis.  First of all, each family on $\mk B$ is $\mc P \times_\mr{lex} \mc Q$-homogeneous. Next, given $n<\om$, $\pi_\mc P"([P\times Q]^{\le n}_\mr{lex})=[P]_\mc P^{\le n}$, and $([P\times Q]^{\le n}_\mr{lex})_\mc P=[Q]_\mc Q^{\le n}$, so $[P\times Q]^{\le n}_\mr{lex}\in \mk B$. Now given $\al$ infinite, we choose $\mc F\in \mk B^\mc P_\al$ and $\mc G\in \mk B^\mc Q_\al$. Then $\mc Z:=\mc F\circledast_\mr{F}\mc G\in \mk B_\al$:  We know from Lemma \ref{kloi8978785522222} that $\mc Z$ is $\al$-homogeneous. Since $\pi_\mc P"(\mc Z)=\mc F$ and $(\mc Z)_\mc P=\mc G$, we have that $\mc Z\in \mk B$. We check now (B.2'): Let $\mc F_0,\mc F_1\in \mk B$. Since $\pi_\mc P"(\mc F_0\cup \mc F_1)=\pi_\mc P"(\mc F_0)\cup \pi_\mc P"(\mc F_1)$ and $(\mc F_0\cup \mc F_1)_\mc P=(\mc F_0)_\mc P \cup (\mc F_1)_\mc P$, we obtain that $\mk B$ is closed under $\cup$.  Secondly, $\pi_\mc P"(\mc F_0\sqcup_{\mr{lex}} \mc F_1)=\pi_\mc P"(\mc F_0)\sqcup_\mc P \pi_\mc P"(\mc F_1)$, and $(\mc F_0)_\mc P ,(\mc F_1)_\mc P \con (\mc F_0\sqcup_\mr{lex} \mc F_1)_\mc P \con (\mc F_0)_\mc P \sqcup_\mc Q (\mc F_1)_\mc P$.  This means that
$$\iota(\mr{srk}_\mc Q((\mc F_0\sqcup_\mr{lex} \mc F_1)_\mc P) )=\max\{\iota(\mr{srk}_\mc Q((\mc F_0)_\mc Q)),\iota(\mr{srk}_\mc Q((\mc F_1)_\mc Q))\}=\iota(\mr{srk}_\mc Q((\mc F_0)_\mc P \sqcup_\mc Q (\mc F_1)_\mc P)),$$
so $(\mc F_0\sqcup_\mr{lex} \mc F_1)_\mc P\in \mk B_\mc Q$. Since in addition $\mc F_0\sqcup_\mr{lex}\mc F_1$ is $\mr{lex}$-homogeneous, we obtain that $\mc F_0\sqcup_\mr{lex}\mc F_1\in \mk B$, by definition.  Finally we check that $\times$ is a multiplication. The property (M.1) of $\times$ follows from Lemma \ref{kloi8978785522222} (c). Now  suppose that $(s_n)_n$ is a sequence in $\mc F\in \mk B$. Let $(t_n)_n$ be a  subsequence of $(s_n)_n$ such that
\begin{enumerate}[(1)]
\item  $(\pi_\mc P"t_n)_n$ is a $\De$-sequence with root $y$.
\item For every $x\in  \mc H$ one has that $\bigcup_{n\in x} \pi_\mc P"(t_n)\in (\pi_\mc P"\mc F) \times_\mc P \mc H$.
\item For every $x\in  \mc H$  and every $p\in y$ one has that $\bigcup_{n\in x} (t_n)_p \in (\mc F)_\mc P \times_\mc Q \mc H$.
\end{enumerate}
Since $(\mc F)_\mc P \con (\mc F)_\mc P \times \mc H$,  the conditions  above imply that given $x\in \mc H$ one has that $\bigcup_{n\in x}t_n\in \mc F\times \mc H$.

(b)  Finite lexicographical powers have bases of families on chains by (a).  For each $n<\om$, let $(\mk B_n, \times_n)$ be a basis on chains of $\mc P^n_\mr{lex}$. Let $\mk B$ be the collection of all $\mr{qlex}$-homogeneous families on $\mc P^{<\om}_\mr{qlex}$ such that
 \begin{enumerate}[(i)]
 \item  $\mc F\rest [P]^n \in \mk B_n$ for every $1\le n<\om$.
 \item  $\mr{lh}"(\mc F):=\conj{\mr{lh}"(s)}{s\in \mc F}\in \mc H$.
 \end{enumerate}
 Given $\mc F\in \mk B$ and $\mc H\in \mk S$, let
 $$\mc F\times \mc H=  (( (\mc F\rest [P]^n) \times_n \mc H  )_n)^{\mr{lh}"(\mc F)\times_\om  \mc H}.$$
 We check that $(\mk B,\times)$ is a pseudo-basis. Given $1\le k<\om$, set $\mc F:=[P^{<\om}]^{\le k}_\mr{qlex}$. Then for each $n$ one has that $\mc F\rest [P]^{n}=[ [P]^n]^{\le k}_\mr{lex}\in \mk B_n$, and
 $\mr{lh}"(\mc F)=[\om]^{\le k}\in \mk B_\om$, so $\mc F\in \mk B$.  One shows as in (a) that $\mk B$ is closed under $\cup$ and $\sqcup_\mr{qlex}$. Finally, we check that $\times$ is a multiplication. That $\times$ satisfies (M.1) it follows from Lemma \ref{kloi8978785522222} (c).  Let now  $\mc F\in \mk B$, $\mc H\in \mk B_\om$, and let $(s_k)_k\in \mc F$.  Let $(t_k)_k$ be a subsequence of $(s_k)_k$ such that
 \begin{enumerate}[(1)]
 \item  $(\mr{lh}"t_k)_k$ is a $\De$-sequence with root $y\con \om$.

 \item  For every $x\in \mc H$ and every $n\in y$ one has that $\bigcup_{k\in x} (s_k \cap [P]^n)\in (\mc F\rest [P]^n)\times_n \mc H$.
 \item   For every $x\in \mc H$ one has that $\bigcup_{k\in x}\mr{lh}"(y_k) \in (\mr{lh}"\mc F \times_\om \mc H)$.
   \end{enumerate}
 It is easy to verify that $\bigcup_{k\in x} t_k\in \mc F\times \mc H$ for every $x\in\mc H$.
\fprue

\subsection{Cardinals smaller than the first Mahlo have a basis}

\defi
A $C$-sequence on $\theta$ is \emph{small} when there is a function $f:\theta\to \theta$ such that $\mr{otp}(C_\al)<f(\min C_\al)$ for every $\al<\theta$.
\fdefi
\prop
A strong limit   cardinal $\theta$ has a small $C$-sequence if and only if $\theta$ is smaller than the first Mahlo cardinal.
\fprop
\prue
Suppose that $\theta$ is smaller than the first Mahlo cardinal. Choose a closed and unbounded set $D\con
\theta$ consisting of non-inaccessible cardinals.  For each $\al<\theta$ let $\la(\al)\in D$ be the maximal
element of $D$ smaller or equal than $\al$, and for each $\la$ in $D$ let $\la_D^+$ be the first element of
$D$ bigger than $\la$. Let $f:\theta\to \theta$, $f(\al)=2^{\la(\al)_D^+}+1$. $f(\al)<\theta$ because we are
assuming that $\theta$ is strong limit.  Observe also that $2^\al<f(\al)$. We define now the $\bar C$
sequence. Fix $\al\in \theta$. Suppose first that $\al\notin D$. Write $\al=\la(\al)+\be$. Since
$\mr{cof}(\al)=\mr{cof}(\be)$, we can choose a club $C_\al\con[\la(\al),\al[$ with
$\mr{otp}(C_\al)=\mr{cof}(\be)$. It follows that $\mr{otp}(C_\al)=\mr{cof} (\be)
<\la(\al)^+_D<f(\la(\al))=f(\min C_\al)$.  Suppose that $\al\in D$. If $\al$ is singular, then we choose
$C_\al$ in a way that $\mr{otp}(C_\al)=\mr{cof}<\min C_\al$. Observe that then $\mr{otp}(C_\al)<\min C_\al <
f(\min C_\al)$. Finally, if $\al$ is not strong limit, then let $\be<\al$ be such that $2^\be \ge \al$. Let
now $C_\al\con [\be,\al[$. It follows that $\mr{otp}(C_\al)=\mr{cof}(\al)\le \al \le 2^\be\le 2^{\min C_\al}
< f(\min C_\al)$.

Suppose now that $\theta$ is bigger or equal to a Mahlo cardinal $\ka$. By pressing down Lemma there is a
stationary set $S$ of inaccessible cardinals of $\ka$ and $\ga<\ka$ such that $\min C_\al=\ga$ for every
$\al\in S$.  Hence $\al=\mr{cof}(\al)=\mr{otp}(C_\al)<f(\ga)$ for every $\al$, and this is of course
impossible.
\fprue

\prue[{\sc Proof of Theorem \ref{firstMahlo}}]
 The proof is by induction on $\theta<\mu_0$.  We see that there is a tree $T=(T,<)$ on $\theta$ that has bases of families on chains of $(T,<)$ and of $(T,<_a)$.  Suppose that $\theta$ is not strong limit.
 Then there is  $\ka<\theta$ such that $\theta\le 2^\ka$.  By Theorem \ref{oiiojoui434344}, $2^\ka$ has a
 basis,  so $\theta$ does.

 Suppose   that $\theta$ is strong limit. Let $T:=T(\ro_0)$ on $\theta$. Let $\bar C$ be a small C-sequence on $\theta$, and let $f:\theta\to
 \theta$ be a witness of it. By inductive hypothesis, and by Proposition \ref{hjgdty42932833} we know that
 the disjoint union $\biguplus_{\xi<\theta}f(\xi)$    of $(f(\xi))_{\xi<\theta}$  has a basis of families on chains. Let now $\la: \mc I(\bar C)\to
 \biguplus_{\xi<\theta}f(\xi)=\bigcup_{\xi<\theta} f(\xi)\times \{\xi\}$, $\la(C):= (\mr{otp}(C),\min C)$.
 Then $\la: (\mc I(\bar C), \sqsubset)\to (\biguplus_{\xi<\theta}f(\xi),<) $ is chain preserving and 1-1 on
 chains. Hence, by Theorem \ref{ui32iuui32t744}, there is a basis of families  on  chains of $\mc I(\bar
 C)_{\mr{qlex}}^{<\om}$.   Since $T$ has infinite chains and Proposition \ref{sijeirjewirowfsdfd} tells that $\ro_0$ is strictly monotone,
 it follows again by  Theorem \ref{ui32iuui32t744} that there is a basis of families on chains of $(T,<)$.

 Observe that $\bigcup_{t\in T} \mr{Is}_t$ is a disjoint union,  $\#T=\theta$.  Since we are assuming that $\theta$ is inaccessible, it follows that $T$ is
 $<\theta$-branching, hence for every $t\in T$ there is a basis of families on $\mr{Is}_t$. Hence, by Proposition
 \ref{hjgdty42932833}, there is a basis of families  on chains of $(T,<_a)$.
\fprue

\section{Subsymmetric sequences and $\ell_1^\al$-spreading models}
We present now new examples of Banach spaces without subsymmetric basic sequences of density $\ka$. Their
construction  uses bases of families on $\ka$. On one side, the multiplication of families will imply the non
existence of subsymmetric basic sequences. On the other side,  the fact that the families are homogeneous
will allow to bound the complexity of finite subsymmetric basic sequences. In fact, we will give examples of
spaces such that every non-trivial sequence on it has a subsequence such that a large family of finite
further subsequences behave like $\ell_1^n$. Since in addition the spaces are reflexive, we will have, as for
the Tsirelson space, that there are no subsymmetric basic sequences.

\defi
Recall that a non-constant  sequence $(x_n)_n$ in a Banach space $X=(X,\nrm{\cdot})$  is called \emph{subsymmetric}
when there is a constant $C\ge 1$ such that
\begin{equation}
\nrm{\sum_{i=1}^n a_i x_{l_i}}\le C\nrm{\sum_{i=1}^n a_i x_{k_i}}
\end{equation}
for every $n$, scalars $(a_i)_{i=1}^n$, $l_1<l_2<\cdots <l_n$ and $k_1<k_2<\cdots<k_n$.
\fdefi
Sometimes it is assumed, not in here, that subsymmetric sequences are  unconditional basic sequences. Notice
that by Rosenthal's $\ell_1$-Theorem and Odell's partial unconditionality, it follows that if $(x_n)_n$ is a
subsymmetric basic sequence, then either is equivalent to the $\ell_1$ unit basis, or its difference sequence
$(x_{2n}-x_{2n-1})_n$ is subsymmetric and unconditional. This is sharp, as it is shown by the summing basis
of $c_0$.

\defi Let $\mc S_\al$ be an $\al$-Schreier family on $\om$. A bounded sequence $(x_n)_n$ in a normed space
$\mk X$ is called an \emph{$\ell_1^\al$-spreading model} when there is a constant $C>0$ such that
\begin{equation*}
\nrm{\sum_{n\in s} a_n x_n}\ge C \sum_{n\in s} |a_n| \text{ for every $s\in \mc S_\al$.}
\end{equation*}
\fdefi
Let us say that a sequence   in a Banach space is \emph{non-trivial} when it does not have norm-convergent
subsequences.
\nota
\begin{enumerate}[(a)]
\item   Suppose that $(\mc S_\al)_{\al<\ou}$ is a generalized Schreier sequence.   If $(x_n)_n$ is a $\ell_1^\al$-spreading model, and $\be\le \al$, then $(x_n)_n$ is a $\ell_1^\be$-spreading model: This is a consequence of the fact that for every $\be<\al$ there is some integer $n$ such that $\mc S_\be \rest \om \setminus n \con \mc S_\al$.

\item Suppose that a space $X$  does not contain $\ell_1$ and it is such that every non-trivial sequence has a  $\ell_1$-spreading model subsequence.
 Then $X$ does not have  subsymmetric sequences. If in addition $X$ have an unconditional basis, then $X$
  is in addition reflexive:  Suppose otherwise that $(x_n)_n$ is a subsymmetric sequence $(x_n)_n$. It
  follows that  $(x_n)_n$ is bounded and that $(x_n)_n$ does not have
   norm-convergent subsequences. So, by hypothesis, there is a $\ell_1$-spreading model subsequence
   $(y_n)_n$. This implies that $(y_n)_n$ is equivalent to the unit basis of $\ell_1$, and this is
   impossible.  The latter condition follows from the James criterion of reflexivity.

\end{enumerate}
\fnota

\defi
Recall that given a family $\mc F$ on $\ka$, we define the corresponding \emph{generalized Schreier space}
$X_\mc F$ as the completion of $c_{00}(\ka)$ with respect to the norm
\begin{equation*}
\nrm{x}_{\mc F}:=\max\{\nrm{x}_\infty,\max_{s\in \mc F}\sum_{\xi\in s}|(x)_{\xi}|  \}.
\end{equation*}
\fdefi
It is easy to see that the unit basis of $c_{00}(\kappa)$ is  a 1-unconditional basis of $X_\mc F$, and that
$X_\mc F$ is $c_0$-saturated if $\mc F$ is compact, and contains  a copy of $\ell_1$ otherwise. When the
family $\mc F$ is compact, hereditary and   $\al$-homogeneous with $\al$ infinite, then every subsequence of
the unit basis of $X_\mc F$ has a $\ell_1$-spreading model subsequence, consequently, no subsequence of the
unit basis is subsymmetric. These families $\mc F$ exist on cardinal numbers not being $\om$-Erd\H{o}s (see
\cite{LoTo}).

\teor
Suppose that $\theta$ is smaller than the first Mahlo cardinal number. Then for every $\al<\ou$ there is a
Banach space $X$ of density $\theta$ with a long 1-unconditional basis $(u_\xi)_{\xi<\theta}$ such  that
every subsequence of $(u_\xi)_{\xi<\theta}$ has a further $\ell_1^\al$-spreading model subsequence, and no
subsequence of $(u_\xi)_{\xi<\theta}$ is a $\ell_1^{\iota(\om^\al)}$-spreading model.
\fteor
\prue
 Fix a basis $(\mk B,\times)$ on $\theta$, let $\mc F$ be an $\om^\al+1$-homogeneous family in $\mk B$  and let $X:= X_\mc F$. Let $(u_\xi)_{\xi\in M}$ be an infinite subsequence of
 the unit basis $(u_\xi)_{\xi<\theta}$ of $X_\mc F$.  Since $\mr{rk}(\mc F\rest M)>\om^\al$, and  $\mr{rk}(\mc S_\al)=\om^\al$ there is some infinite subset $N$ of $M$ such that
  $\{\xi_n\}_{n\in x}\in \mc F$ for every $x\in \mc S_\al\rest N$. Let $N=\{n_k\}_k$ be the increasing enumeration of $N$, and set $x_k:= u_{\xi_{n_k}}$ for every $k<\om$. We claim that
 $$\nrm{\sum_{k\in x} a_k x_k}_\mc F=\sum_{k\in x}|a_k|$$
 for every $x\in \mc S_\al$:  Fix $x\in \mc S_\al$. Then $\{n_k\}_{k\in x}\in \mc S_\al\rest N$, because $\mc S_\al$ is spreading. This means that $\{\xi_{n_k}\}_{k\in x}\in \mc F$, so
 $$\nrm{\sum_{k\in x} a_k x_k}_\mc F =  \nrm{\sum_{k\in x} a_k u_{\xi_{n_k}}}_\mc F\ge \sum_{k\in x} |a_k|.$$
On the other hand, let given a subsequence $(x_n)_{n<\om}$ be a subsequence of $(u_\xi)_{\xi<\theta}$,
$x_n:=u_{\xi_n}$. we assume that $\xi_n<\xi_{n+1}$ for every $n$. Let $\mc G:= \conj{x\con
\om}{\{\xi_n\}_{n\in x}\in \mc F}$. Then the mapping $x\in \mc G\mapsto  \{\xi_n\}_{n\in x}\in \mc F$ is
continuous and 1-1, hence $\mr{rk}(\mc G)\le \mr{rk}(\mc F)<\iota(\om^\al)$.  Since $\iota(\om^\al)$ is
exp-indecomposable, $\mr{rk}(\mc S_{\iota(\om^\al)})=\om^{\iota(\om^\al)}=\iota(\om^\al)$.  It follows by the
quantitative version of Ptak's Lemma  (see for example \cite[Lemma 4.7]{LoTo3}) that for every $\vep>0$
there is some convex combination $(a_n)_{n\in x}$ supported in $x\in \mc S_{\iota(\om^\al)}$ such that
$$\sup_{x\in \mc G}\sum_{n\in x} |a_n|<\vep.$$
This means that $\nrm{\sum_{n\in x}a_n x_n}_\mc F<\vep\sum_{n\in x}|a_n|$.
\fprue
\defi
Recall that given   an $\al$-Schreier family $\mc S_\al$, let $T_\al:=T_{\mc S_\al}$ be the $\al$-Tsirelson
space defined as the completion of $c_{00}$ under the norm
\begin{equation*}
\nrm{x}_\al:= \max\{\nrm{x}_\infty, \sup_{(E_i)_i} \frac{1}{2}\sum_{i} \nrm{E_i x}_\al\}
\end{equation*}
where the $\sup$ above runs over all sequences of sets $(E_i)_i$ such that $E_i<E_{i+1}$ and $\{\min
E_i\}_i\in \mc S_\al$, and where $Ex =\sum_{n\in E} (x)_n$.
\fdefi
An equivalent way of defining $\nrm{\cdot}_\al$ is as follows. Let $K_0:=\{\pm u_n\}_{n<\om}$ and let
$$K_{n+1}:=K_n\cup \conj{\frac12\sum_{i<k}\vphi_i}{\{\vphi_i\}_{i<k}\con K_n \text{, $\vphi_i<\vphi_{i+1}$,
$i<k-1$, and $\{\min \supp\vphi_i\}_{i<k}\in \mc S_\al$}}.$$ Let $K:=\bigcup_n K_n$. Then
$\nrm{x}_\al=\sup_{\vphi\in K}\langle\vphi,x\rangle$.  It is easy to see that each $\vphi\in K$ has a
decomposition
$$\vphi=\sum_{i}\frac{1}{2^i} \vphi_i$$
where $\vphi_i$ is a vector with coordinates $-1$ or 1, supported in $\mc S_{\al \cdot i}$ and $(f_i)_i$
pairwise disjointly supported. It is well known that every normalized block subsequence of the unit basis
$(u_n)_n$ of $T_\al$ is equivalent to a subsequence of the unit basis. Since clearly from the definition
every subsequence  of $(u_n)_n$ is a $\ell_1^\al$-spreading model, it follows that every non-trivial sequence
in $T_\al$ has a $\ell_1^\al$-spreading model subsequence. Now suppose that $(x_n)_n$ is a non-trivial
sequence in $T_\al$. W.l.o.g. we assume that $(x_n)_n$ is a subsequence of the unit basis, $x_n=u_{k_n}$. Let
$\vep>0$, and let $n$ be such that $\vep 2^n>1$. By the quantitative version of Ptak's Lemma, there is some
convex combination $(a_n)_{n\in x}$ supported in $x\in \mc S_{\al\cdot\om }$  such that $\sup_{x\in
\bigcup_{i<n}\mc S_{\al\cdot n}} \sum_{k\in x}|a_k|<\vep$.  We claim that $\nrm{\sum_{n\in x} a_n x_n}_\al\le
3\vep$: Fix $\vphi\in K$,  $\vphi=\sum_i 2^{-i}\vphi_i$  decomposed as above. Then,
\begin{align*}
|\langle \vphi, \sum_{j\in x} a_j u_{k_j}\rangle| \le & \sum_{i<n} \frac{1}{2^i}|\langle \vphi_{i}, \sum_{j\in x}a_j u_{n_j}\rangle| +\frac{1}{2^n}\sum_{j\in x}|a_j| \le  \sum_{i<n} \frac{1}{2^i}\sum_{j\in x\cap \supp \vphi_j}|a_j| +\frac{1}{2^n}\sum_{j\in x}|a_j| \le 3\vep.
\end{align*}
We have just proved the following.
\teor\label{ioj43iotio777}
$T_\al$ is a reflexive Banach space whose unit basis is 1-unconditional and such that every non-trivial
sequence has a  $\ell_1^\al$-spreading model subsequence but it does not have $\ell_1^{\al\cdot
\om}$-spreading models.  Consequently, $T_\al$ does not  have subsymmetric  basic sequences. \qed
\fteor

\subsection{The interpolation method}\label{ffsdfd}

We recall the following well-known construction, presented in a general, not necessarily separable, context:
fix an infinite cardinal number $\ka$, let $(\nrm{\cdot}_n)_{n\in \om}$ be a sequence of norms in $c_{00}(\ka)$ and
$\nrm{\cdot}_X$ be a norm on $c_{00}(\N)$ such that $(e_n)_n$ is a 1-unconditional basic sequence of the completion
$X$ of $(c_{00}(\N),\nrm{\cdot}_X)$. Let $X_n$, $n\in \N$, be the completion of $(c_{00}(\ka),\nrm{\cdot}_n)$. For $x\in
c_{00}(\ka)$, define
$$\tnrm{x}:=\nrm{\sum_{n}\frac{\nrm{x}_n}{2^{n+1}}e_n}_X.$$
It is not difficult to see that $\tnrm{\cdot}$ is a norm on $c_{00}(\ka)$ (the fact that $(u_n)_n$ is a
1-unconditional basic sequence of $(c_{00}(\N),\nrm{\cdot}_X)$ is crucial to prove the triangle inequality). Let
$\mk X$ be the completion of $(c_{00},\tnrm{\cdot})$.

\nota\label{n4hui3t8767rrreer}
Observe that the dual unit ball of $\eqs$ is closed under the following operation. Given  $x_i^*\in
B_{X_i^*}$ for $i=1,\dots,n$ and $\sum_{i=1}^n b_i e_i^*\in B_{X^*}$, then
$$\sum_{i=1}^n \frac{b_i}{2^{i+1}} x_i^*\in B_{\eqs^*}.$$
To see this, there is a simple computation. Given $x\in c_{00}$ we have that
\begin{align*}
 |(\sum_{i=1}^n \frac{b_i}{2^{i+1}} x_i^*)(x)|\le & \sum_{i=1}^n \frac{|b_i|}{2^{i+1}} \nrm{x}_i =|(\sum_{i=1}^n |b_i| e_i^*)(\sum_{i=1}^n \frac{1}{2^{i+1}}\nrm{x}_i e_i)   | \le \\
 \le & \left\|\sum_{i=1}^n |b_i| e_i^*\right\|_{X^*}  \left\|\sum_{i=1}^n \frac{1}{2^{i+1}}\nrm{x}_i e_i \right\|_X \le\nrm{x}_{\eqs}
\end{align*}
\fnota
The following follows easily from the definition.
\prop
Suppose that $(x_\xi)_{\xi<\la}$ is a $C$-unconditional basic sequence of each $X_n$. Then
$(x_\xi)_{\xi<\la}$ is a $C$-unconditional basic sequence of $\mk X$. \qed
\fprop

In our construction, this will be the case, so that we will be able to apply the following result.

\prop \label{iijgihgnjkbngnjkui74}
Suppose that $X$ is a space with an unconditional basis and without isomorphic copies of $\ell_1$. Then the
following are equivalent.
\begin{enumerate}[(a)]
\item  Every non trivial bounded sequence in $X$ has an   $\ell_1^\al$-spreading model subsequence.
\item  Every non-trivial weakly-convergent sequence in $X$ has $\ell_1^\al$-spreading model subsequence.
\item  Every non-trivial weakly-null sequence in $X$ has $\ell_1^\al$-spreading model subsequence.
\end{enumerate}
\fprop
\prue
Suppose that (b) holds. It follows that $c_0$ does not embed into $X$. Hence, by James' criteria of
reflexivity for spaces with an unconditional basis, $X$ is reflexive, whence every bounded sequence has a
weakly-convergent subsequence and now (a) follows directly from (b). Now suppose that (c) holds and let us
prove (b): Suppose that $(x_k)_k$ is a non-trivial weakly-convergent sequence with limit $x$.  Let $x^*\in
S_{X^*}$ be such that $x^*(x)= \nrm{x}$.     Let $y_k:=x_k-x$ for every $k$. By hypothesis, we can find
$\vep>0$ and a subsequence $(z_n)_{n}$ of $(y_n)_{n}$ such that $\nrm{\sum_{n} a_n y_n} \ge \vep \sum_n
|a_n|$ for every sequence of scalars $(a_n)_n$ supported in $\mc S_\al $.   Let $(v_n)_n$ be a further
subsequence of $(z_n)_n$ such that $|x^*(v_n)|\le \vep/2$ for every $n$.  We claim that $(v_n +x)_n$ is a
subsequence of $(x_n)_n$  which is a $\ell_1^\al$-spreading model: Fix a sequence $(a_n)_n$ supported in $\mc
S_\al$, and let $y^*\in B_{X^*}$ be such that $y^*(\sum_n a_n z_n) \ge \vep \sum_n |a_n| $. Let $z^*:= y^* -
\la x^* \in \ker(x)\cap 2 B_{X^*}$, where $\la:= y^*(x)/\nrm{x}$. Then,
\begin{align*}
\nrm{\sum_n a_n (v_n+x)} \ge & \frac{1}{2} z^*\left(\sum_n a_n (v_n+x)\right)= \frac{1}{2} z^*\left(\sum_n a_n v_n\right) \ge \frac12 y^*\left(\sum_n a_n v_n\right)- \\ \pushQED{\qed}
- &  \la x^*\left(\sum_n a_n v_n\right)   \ge \frac{\vep}{4} \sum_n |a_n|.   
\qedhere
\popQED
\end{align*}\let\qed\relax
\fprue
Finally, one of the interesting features of the resulting space of the interpolation is given by the
following proposition.

\prop\label{lkjjrjrfffff333}
If $Y$ can be isomorphically embedded into $\mk X$, then a subspace $Z$ of $Y$ can be isomorphically
embedded  into some $X_n$ or into $X$.
\fprop
\prue
Let $Y_0\con \mk X$ be isomorphic to $Y$. If there is some $n$ such that $I_n: Y_0 \to X_n$ is not strictly
singular, then choosing $Y_1\con Y_0$ such that $I_n\rest Y_1$ is an isomorphism we obtain that $Y_1$ embeds
into $X_n$.  Suppose that all $I_n\rest Y_0$ are strictly singular.

Fix a strictly positive summable sequence $(\vep_n)_n$, $\sum_n \vep_n<1/4$.  Let $y_0\in S_{Y_0}$.    Let
$n_0\in \N$ be such that
\begin{equation*}
\nrm{\sum_{n>n_0}\frac{1}{2^{n+1}} \nrm{y_0}_n e_n}_X \le \vep_0.
\end{equation*}
Let $y_1\in S_{Y_0}$  be such that
\begin{equation*}
\max_{n\le n_0} \nrm{y_1}_n\le \vep_1
\end{equation*}
Let $n_1>n_0$ be such that
\begin{equation*}
\nrm{\sum_{n>n_1}\frac{1}{2^{n+1}} \nrm{y_1}_n e_n}_X \le \vep_1.
\end{equation*}
Let now $y_2\in S_{Y_0}$ be such that
\begin{equation*}
\max_{n\le n_1} \nrm{y_2}_n\le \vep_2.
\end{equation*}
In this way, we can find an strictly increasing sequence $(n_k)_{k\in \N}$ of integers normalized vectors
$y_n\in Y_0$ such that for every $k$ one has that
\begin{align*}
\nrm{\sum_{n>n_k}\frac{1}{2^{n+1}} \nrm{y_k}_n e_n}_X & \le \vep_k   \\
\max_{n\le n_{k-1}} \nrm{y_{k}}_n & \le \vep_k.
\end{align*}
Set   $w_0:=\sum_{n\le n_0} (\nrm{y_0}_{n}/2^{n+1})e_n \in X$, and for each $k\ge 1$, let
\begin{equation*}
w_k:=\sum_{n=n_{k-1}+1}^{n_k} \frac{\nrm{{y_k}_n}}{2^{n+1}} e_n \in X.
\end{equation*}
Let $(a_k)_k$ be a sequence of scalars with $\max_k |a_k|=1$. Then, after some computations, one can show
that
\begin{align*}
\left| \nrm{\sum_{k} a_k y_k}_\mk X - \nrm{\sum_{k} a_k w_k}_X   \right| \le &
|a_0|
\nrm{\sum_{n\le n_0} \frac{1}{2^{n+1}} \nrm{y_0}_n e_n}_X + \\
+&  \sum_{k\ge 1} |a_k|\left(  \nrm{\sum_{n\le n_{k-1} } \frac1{2^{n+1}} \nrm{y_k}_n e_n   }_X +    \nrm{\sum_{n>n_{k} } \frac1{2^{n+1}} \nrm{y_k}_n e_n   }_X  \right)< \\
< &\frac12.
\end{align*}
Since $\nrm{w_k}_X\ge 3/4$, and $(w_k)_k$ is a block subsequence of the basis $(e_n)_n$, it follows that
$\nrm{\sum_k a_k w_k}_X\ge 3/4 \max_k |a_k|$, we obtain that
\begin{equation*}
 \left| \nrm{\sum_{k} a_k y_k}_\mk X - \nrm{\sum_{k} a_k w_k}_X   \right| <\frac{2}{3}\nrm{\sum_k a_k w_k}_X.
\end{equation*}
Hence,
\begin{equation*}\pushQED{\qed}
\frac13 \nrm{\sum_k a_k w_k}_X\le \nrm{\sum_{k}a_k y_k}_\mk X\le \frac43\nrm{\sum_k a_k w_k}_X. \qedhere
\popQED
\end{equation*}\let\qed\relax
\fprue

\subsection{The Banach space $\mathfrak{X}$}
The next result gives the existence of the desired Banach space $\mk X$ of density $\ka$, subject to the
existence of bases of families   on  $\ka$.

\teor\label{lkklklklerefvvvvmjk}
Suppose that $\ka$ has a basis of families. Then  for every $1\le \al<\ou$   there is a reflexive Banach
space $\mk X_\al$ of density $\ka$ with a long unconditional basis such that
\begin{enumerate}[(a)]
\item   every bounded sequence without norm convergent subsequences has a $\ell_1^\al$-spreading model subsequence, and
\item  $\mk X_\al$ does not have $\ell_1^{\iota(\om^\al) + \al \cdot \om }$-spreading models.
\end{enumerate}
Consequently,
\begin{enumerate}[(a)]\addtocounter{enumi}{2}
\item $\mk X_\al$ does not have subsymmetric sequences;
\item  if $\iota(\om^\al) +\al \cdot \om \le \be$, then $\mk X_\al$ and $\mk X_\be$ are totally incomparable, i.e. there is no infinite dimensional subspace of $\mk X_\al$ isomorphic to a subspace of $\mk X_\be$.
\end{enumerate}

\fteor

\prue  Let $(\mk B,\times)$ be a basis of families on $\ka$. Let $\mc F_0:=[\ka]^{\le 1}$, and for $n<\om$ let $\mc F_{n+1}:= \mc F_n \times  \mc S_\al$.  Notice that $\mr{rk}(\mc F_n)<\iota(\om^\al)$ for every $n$.
Let $\mk X$ be the interpolation space  from an $\al$-Tsirelson space  $T_\al$ and the sequence of
generalized Schreier spaces $X_{\mc F_n}$, $n\in \N$. Since each $X_{\mc F_n}$ is $c_0$-saturated, and
$T_\al$ is reflexive,  it follows from Proposition \ref{lkjjrjrfffff333} that  $\mk X$ does not have
isomorphic copies of $\ell_1$.
\clam
Every non-trivial bounded sequence  has an  $\ell_1^\al$-spreading model.
\fclam
From this claim, and the fact that the unit basis $(u_\ga)_{\ga<\ka}$ is unconditional, we obtain that $\mk
X$ is reflexive.  We pass now to prove that previous claim.
\prucl
Fix such sequence $(x_k)_k$.   Since $\mk X$ does not have isomorphic copies of $\ell_1$,  we may assume, by
Proposition \ref{iijgihgnjkbngnjkui74}, that   $(x_k)_k$ is  a non-trivial weakly-null sequence.  Since
$(u_\xi)_{\xi<\ka}$ is a Schauder basis of $\eqs$, by going to a subsequence if needed, we assume that
$(x_k)_k$ is disjointly supported, and
\begin{equation*}
\nrm{x_k}\ge \ga>0 \text{ for every $k$.}
\end{equation*}
 The proof is now rather similar to that of Proposition \ref{lkjjrjrfffff333}.

\caso{1} There is $\vep>0$, $n\in \N$ and an infinite subsequence $(y_k)_k$ of $(x_k)_k$ such that
\begin{equation*}
\nrm{y_k}_n\ge \vep \text{ for every $k$.}
\end{equation*}
For each $k$ choose $s_k\in \mc F_n\rest  \supp y_{k}$ such that
\begin{equation*}
\sum_{\xi\in s_k} |u_{\xi}^*(y_{k}) |\ge \vep.
\end{equation*}
By hypothesis, $\mc F_{n+1}=\mc F_n \times \mc S_\al$, so  there is a subsequence $(t_k)$ of $(s_k)_k$ such
that $\bigcup_{k\in v} t_k\in \mc F_{n+1}$ for every $v\in \mc S_\al$. Let $(z_k)_k$ be the subsequence of
$(y_k)_k$ such that $\sum_{\xi\in t_k} |u_\xi^*(z_k)|\ge \vep$ for every $k$. We claim that $(z_k)_k$ is a
$\ell_1^\al$-spreading model.  So, fix a sequence of scalars $(a_k)_{k\in s}$ indexed by  $s\in \mc S_\al$.
Then $t:=\bigcup_{k\in s} t_k\in \mc F_{n+1}$, hence,
\begin{align*}
\nrm{\sum_{k\in s} a_k z_k}\ge & \frac{1}{2^{n+2}} \nrm{\sum_{k\in s} a_k z_k}_{n+1}\ge \frac{1}{2^{n+2}}\sum_{\xi\in t} \left| u_\xi^*( \sum_{k\in s} a_k z_k)    \right|= \frac{1}{2^{n+2}}\sum_{k\in s} |a_k|\sum_{\xi\in t_k} \left| u_\xi^*(  z_k)    \right| \ge \\
\ge & \frac{\vep}{2^{n+2}} \sum_{k\in s} |a_k|.
\end{align*}

\caso{2} For every $\vep>0$ and every $n$ the set
\begin{equation*}
\conj{k\in \N}{ \nrm{x_k}_n\ge \vep}\text{ is finite.}
\end{equation*}
So, let $(y_k)_k$ be a subsequence of $(x_k)_k$ such that
\begin{equation*}
 \nrm{y_k}_{k-1}\le \frac{1}{2^{k+1}} \text{ for every $k$}
\end{equation*}
Now we can find $(n_k)_k$ such that for every $k$ one has that
\begin{align*}
\nrm{\sum_{n=n_k}^{n_{k+1}-1}\frac{1}{2^{n+1} }\nrm{y_k}_n t_n}_T\ge & \frac{\ga}2 \\
\nrm{\sum_{n\notin[n_k,n_{k+1}[}\frac{1}{2^{n+1} }\nrm{y_k}_n t_n}_T\le & \frac{\ga}{2^{k+4}}
\end{align*}
For every $k$ choose  $\sum_{n=n_k}^{n_{k+1}-1}b_n t_n^*\in S_{T_\al^*}$ and $s_n\in \mc{F}_n$ with $n\in
[n_k,{n_{k+1}}[$  such that
\begin{equation*}
\sum_{n=n_k}^{n_{k+1}-1}\frac{1}{2^{n+1} } b_n \psi_n(y_n)  \ge \frac{\ga}2 \text{ for every $k$,}
\end{equation*}
where $\psi_n:= \sum_{\xi \in s_n } \vep_\xi u_\xi^*\in \mc B_{(X_{\mc F_n})^*}$ and $\vep_\xi$ is the sign
of $|u_\xi^*(z_k)|$ for every $k$, every $n\in [n_k,n_{k+1}[$ and every $\xi\in s_n$. For every $k$ let
\begin{equation*}
\vphi_k:=\sum_{n=n_k}^{n_{k+1}-1} \frac{b_n}{ 2^{n+1}} \psi_n\in B_{\eqs^*},
\end{equation*}
because of Remark \ref{n4hui3t8767rrreer}. Notice that $\vphi_k (y_k)\ge \ga/2$ for every $k$.  Given $s\in
\mc S_\al$ we have that  $\{n_k\}_{k\in s}\in \mc S_\al$, because $\mc S_\al$ is a spreading family. Hence,
\begin{equation*}
\frac12\sum_{k\in s} \sum_{n=n_k}^{n_{k+1}-1}b_n t_n^*\in  B_{T_\al^*}.
\end{equation*}
 It follows from this and   Remark \ref{n4hui3t8767rrreer} that
\begin{equation*}
\frac12\sum_{k\in s} \vphi_k= \frac12\sum_{k\in s} \sum_{n=n_k}^{n_{k+1}-1}\frac{1}{2^{n+1}}b_n \psi_n\in  B_{\mk X^*}.
\end{equation*}
Now, fix $s\in \mc S_\al$ and scalars $(a_k)_{k\in s}$, and for each $k\in s$, let $\sig_k$ be the sign of
$a_k$. Then,
\begin{align*}
\nrm{\sum_{k\in s} a_k y_k}\ge & |\langle \frac{1}{2}\sum_{k\in s}\sig_k\vphi_k ,\sum_{k\in s}a_k y_k\rangle|
\ge  \frac{\ga}{4}\sum_{k\in s}|a_k| -\sum_{k\in s} |\langle\frac12 \sum_{j\in s\setminus\{k\}} \vphi_j,  a_ky_k\rangle| \ge  \nonumber\\
\ge & \frac{\ga}{4}\sum_{k\in s}|a_k| -\sum_{k\in s} |a_k|\left\|\sum_{n\notin[n_k,n_{k+1}[}\frac{1}{2^{n+1} }\nrm{y_k}_n t_n\right\|_{T_\al} \ge\\
\ge & \frac{\ga}{4}\sum_{k\in s}|a_k| -\sum_{k\in s} |a_k|\frac{\ga}{2^{k+4}}   \ge \frac{\ga}{8}\sum_{k\in s}|a_k| .
\end{align*}
Hence, $(y_k)_k$ is an $\ell_1^\al$-spreading  model.
 \fprucl
(b):  Suppose otherwise  that $(x_k)_k$ is a weakly-null sequence such that $\nrm{\sum_{k\in s} a_k x_k}_\mk
X \ge \vep \sum_{k\in s} |a_k|$ for every sequence of scalars $(a_k)_{k\in s}$ supported in $s\in \mc
S_{\be}$, for $\be:=\iota(\om^\al)+ \al\cdot \om$.  There are two cases to consider.

\noindent{\sc Case 1.} There is some subsequence $(y_k)_k$ of $(x_k)_k$ some $n$ and some $C>0$ such that
\begin{equation*}
\frac{1}{2^{n+1}} \nrm{\sum_{k \in s} a_k y_k}_n \le \nrm{\sum_{k\in s} a_k y_k}_\mk X \le C\nrm{\sum_{k\in s}a_ky_k}_n
\end{equation*}
for every $(a_k)_{k\in s}$ supported in $s\in \mc S_{\iota(\om^\al)}$.   Let $\theta: \mc F_n \to \fin$,
$\theta(s):=\conj{k\in \om}{\supp y_k \cap s\neq \buit}$. This is a continuous mapping, so $\mc C:= \theta"
\mc F_n$ is compact, and $\ga:=\mr{rk}(\mc C)\le \mr{rk}(\mc F_n)<\iota(\om^\al)$, by the choice of $\mc
F_n$. By the quantitative version of  by Ptak's Lemma we can find a convex combination $(a_k)_{k\in s}$
supported in $\mc S_{\iota(\om^\al)}$ such that $\sum_{k\in t} |a_k|< \vep /(C \sup_k \nrm{y_k}_n)$ for every
$t\in \mc C$.  Let $v\in \mc F_n$ be such that
$$\nrm{\sum_{k\in s} a_k y_k}_n = \sum_{\xi\in v} |u_\xi^*(\sum_{k\in s} a_k y_k)|.$$ Then
\begin{align*}
\nrm{\sum_{k\in s} a_k y_k}_\mk X \le C \nrm{\sum_{k\in s}a_k y_k}_n  \le \sum_{k\in \theta(s)} |a_k| \nrm{y_k}_n<\vep,
\end{align*}
and this is impossible.

\noindent {\sc Case 2.} There is a normalized block subsequence $(y_k)_k$ of $(x_k)_k$ with $y_k:= \sum_{i\in
s_k} b_i x_i$  and  $s_k\in \mc S_{\iota(\om^\al)}$ and that is 2-equivalent to a block subsequence $(w_k)_k$
of the basis $(t_n)_n$ of $T_\al$. Since for every $\de,\ga$ there is an integer $n$ such that $\mc
S_\de\otimes \mc S_\ga \rest (\om\setminus n)\con \mc S_{\de+\ga}$, we assume without loss of generality
$\bigcup_{k\in x} s_k \in \mc S_{\iota(\om^\al) +\al\cdot \om}$ for every $x\in \mc S_{\al\cdot \om}$.   Let
$K:= \sup_k \nrm{x_k}_\mk X$.  By Theorem \ref{ioj43iotio777}, let $(a_k)_{k\in v}$ be a sequence supported
in $v\in \mc S_{\al\cdot \om}$ such that $\nrm{\sum_{k\in v}a_k w_k}_{T_\al}<(\vep/2 K) \sum_{k\in v} |a_k|$.
Hence,
\begin{equation*}
\nrm{\sum_{k\in v} a_k \sum_{i\in s_k}b_i x_i}_\mk X \le 2 \nrm{\sum_{k\in v} a_k w_k}_{T_\al} <\frac{\vep}K \sum_{k\in v}| a_k| \le \vep \sum_{k\in v} |a_k|\sum_{i\in s_k}  |b_i|,
\end{equation*}
and this is impossible because $\bigcup_{k\in v} s_k \in \mc  S_{\iota(\om^\al)+ \al\cdot \om}$.
 \fprue


\begin{thebibliography}{99}

\bibitem[ArGoRo]{ArGoRo} S. A. Argyros, G. Godefroy and H. P. Rosenthal, \emph{Descriptive set theory and Banach spaces.} Handbook of the geometry of Banach spaces, Vol. 2, 1007--1069, North-Holland, Amsterdam, 2003.

\bibitem[ArLoTo]{ArLoTo} S. A. Argyros, J. Lopez-Abad and S. Todorcevic, \emph{ A class of Banach spaces
with few non-strictly singular operators}, {J. Funct. Anal.} \textbf{222} (2005), no. 2, 306--384.

\bibitem[ArTod]{ArTo} S. A. Argyros and S. Todorcevic, \emph{Ramsey methods in analysis.} Advanced Courses in Mathematics. CRM Barcelona. Birkh\"{a}user Verlag, Basel, 2005

\bibitem[ArMo]{ArMo} S. A. Argyros and P.  Motakis,  \emph{$\al$-Large families and applications to Banach space theory}, Topology and its Application \textbf{72} (2014),  47--67.

\bibitem[ArTo]{ArTol} S.A. Argyros, A. Tolias, \emph{Methods in the theory of hereditarily indecomposable Banach spaces,} Mem. Amer. Math. Soc. \textbf{70} (806) (2004).

\bibitem[Be]{Be} S. F. Bellenot, \emph{Tsirelson superspaces and $\ell_p$}, J. Funct. Anal. \textbf{66}
(1986), 2, 207--228.



\bibitem[DoLoTo]{DoLoTo} P. Dodos, J. Lopez-Abad and S. Todorcevic,  \emph{Unconditional basic sequences in spaces of
large density},  {Advances in Mathematics}, \textbf{226} (2011), 3297--3308.



\bibitem[ErHa]{ErHa} P. Erd\H{o}s and A. Hajnal, \emph{On a problem  of B. Jonsson,} Bull. de l'Acad. Polon. Sci. 14 (1966) 19--23.

\bibitem[GoMa]{GoMa} W.T. Gowers, B. Maurey, \emph{The unconditional basic sequence problem,} J. Amer. Math. Soc. 6 (1993) 851--874.



\bibitem[LiTz]{LiTz}J. Lindenstrauss and L. Tzafriri,
\emph{ Classical Banach spaces I}, Springer-Verlag  Vol. 92, (1977).

\bibitem[Lo]{Lo}  J. Lopez-Abad, \emph{Families of finite subsets of $\N$}, Zb. Rad. (Beogr.) 17 (25) (2015), Selected topics in combinatorial analysis, 145--169.
\bibitem[LoTo1]{LoTo3}  J. Lopez-Abad and S. Todorcevic \emph{Pre-compact families of finite sets of integers and weakly null sequences
in Banach spaces}, Topology and its Applications \textbf{156} (2009) 1396--1411.
\bibitem[LoTo2]{LoTo} J. Lopez-Abad and S. Todorcevic \emph{Positional graphs and conditional structure of weakly null sequences,} Advances in Mathematics \textbf{242} (2013) 163-186.



\bibitem[Ke]{Ke} J. Ketonen, \emph{Banach spaces and large
cardinals}, Fund. Math., \textbf{81} (1974), 291--303.

\bibitem[Kr]{Kr}
J. B. Kruskal, \emph{Well-quasi-ordering, the tree theorem, and Vazsonyi's conjecture}, Transactions of the
American Mathematical Society {\bf 95} (2): 210--225.

\bibitem[La]{La} R. Laver, \emph{On Fra\"{\i}ss\'{e}'s order type conjecture}, Ann. Math. {\bf 93} (1971), 89--111.

\bibitem[MaRo]{MaRo} B. Maurey and H. P.  Rosenthal,
\emph{Normalized weakly null sequence with no unconditional subsequence}, Studia Math. \textbf{61} (1977),
no. 1, 77--98.

\bibitem[Na1]{Na1}  C. St.J. A. Nash-Williams,  \emph{On well-quasi-ordering finite trees}, Proc. Of the Cambridge Phil. Soc. {\bf 59} (04): 833--835.
\bibitem[Na2]{NW} C. St. J. A. Nash-Williams, \emph{On well-quasi-ordering transfinite sequences,}
 {Proc. Cambridge Philos. Soc.} {\bf 61} (1965), 33--39.

\bibitem[Od1]{Od} E. Odell, \emph{A nonseparable Banach space not containing a spreading basic sequence}, Israel J. Math. \textbf{52} (1985), no. 1-2, 97--109.

\bibitem[Od2]{Od1} E. Odell, On Schreier unconditional sequences, Contemp. Math. \textbf{144} (1993), 197-201.


\bibitem[Pu-Ro]{PR}P. Pudlak and V. R\"{o}dl, \textit{Partition theorems
for systems of finite subsets of integers}, Discrete Math., \textbf{39}, (1982), 67-73.


\bibitem[Schr]{Schr}
    {J. Schreier,}
    \emph{Ein Gegenbeispiel zur Theorie der schwachen Konvergenz.}
    Studia Math. \textbf{2} (1930), 58--62.



\bibitem[Si]{Si} W. Sierpinski, \emph{Cardinal and ordinal numbers}. Second revised edition. Monografie Matematyczne, Vol. 34
Panstowe Wydawnictwo Naukowe, Warsaw 1965.

\bibitem[To]{To} S. Todorcevic, \emph{Walks on ordinals and their characteristics.}
Progress in Mathematics, \textbf{263}. Birkh\"{a}user Verlag, Basel, 2007.
\bibitem[Tsi]{Tsi} B.S. Tsirelson, \emph{Not every Banach space contains $l_p$ or $c_0$,} Funct. Anal. Appl. \textbf{8} (1974) 138--141.



\end{thebibliography}
\end{document}